\documentclass[11pt]{amsart}
\usepackage{amssymb, adjustbox, enumerate, amsbsy, stmaryrd}
\usepackage{array, longtable}
\usepackage{geometry}
\usepackage{todonotes}

\usepackage{amsfonts, amssymb, amscd}
\numberwithin{equation}{section}

\usepackage[symbol]{footmisc}

\usepackage{bm}
\usepackage{verbatim}
\usepackage{mathrsfs}
\usepackage{graphicx}
\usepackage{tikz-cd}
\usepackage{subcaption}
\usepackage{listings}
\usepackage{subfiles}
\usepackage[toc,page]{appendix}
\usepackage{mathtools}
\usepackage{comment}
\usepackage{enumerate}
\usepackage{enumitem}
\usepackage[all]{xy}

\usepackage{graphicx}
\graphicspath{{images/}}

\usepackage{appendix}
\usepackage{hyperref}
\hypersetup{
    colorlinks=true,
    citecolor=red,
    linkcolor=blue,
    filecolor=magenta,      
    urlcolor=red,
}
\newcommand{\Dd}{{\mathcal{D}}}

\newcommand{\Qq}{\mathbb{Q}}

\newcommand{\Rr}{\mathbb{R}}

\newcommand{\vol}{\operatorname{vol}}
\newcommand{\Center}{\operatorname{center}}

\newcommand{\glct}{\operatorname{glct}}

\newcommand{\mld}{\operatorname{mld}}
\newcommand{\tmld}{\operatorname{tmld}}

\newcommand{\lct}{\operatorname{lct}}
\newcommand{\act}{\operatorname{act}}

\newcommand{\pet}{\operatorname{pet}}

\newcommand{\Supp}{\operatorname{Supp}}

\newcommand{\mult}{\operatorname{mult}}
\newcommand{\Rct}{\operatorname{Rct}}

\newcommand{\Ft}{\operatorname{Ft}}
\newcommand{\CY}{\operatorname{CY}}

\newcommand{\Ii}{\Gamma}

\newcommand{\Sing}{\mathrm{Sing}}

\newtheorem{thm}{Theorem}[section]

\newtheorem{cor}[thm]{Corollary}
\newtheorem{lem}[thm]{Lemma}
\newtheorem{prop}[thm]{Proposition}

\theoremstyle{definition}
\newtheorem{defn}[thm]{Definition}
\newtheorem{agree}[thm]{Agreement}
\newtheorem{ques}[thm]{Question}
\theoremstyle{definition}
\newtheorem{rem}[thm]{Remark}

\newtheorem{notalem}[thm]{Notation-Lemma}
\newtheorem{ex}[thm]{Example}
\newtheorem{nota}[thm]{Notation}
\newtheorem{exlem}[thm]{Example-Lemma}
\newtheorem{cons}[thm]{Construction}
\newtheorem{setup}[thm]{Set-up}

\theoremstyle{definition}

\begin{document}

\title{Optimal bounds on surfaces}

\author{Jihao Liu}
\address{Department of Mathematics, Northwestern University, 2033 Sheridan Rd, Evanston, IL 60208, USA}
\email{jliu@northwestern.edu}

\author{Vyacheslav V. Shokurov}
\address{Department of Mathematics, Johns Hopkins University, 3400 N. Charles Street, Baltimore, MD 21218, USA}
\email{vshokur1@jhu.edu}
\address{Steklov Mathematics Institute, Russian Academy of Sciences, Cubkina Str. 8, 119991, Moscow, Russia}
\email{shokurov@mi-ras.ru}

\date{\today}

\begin{abstract}
We prove that the first gap of $\mathbb R$-complementary thresholds of surfaces is $\frac{1}{13}$. More precisely, the largest $\mathbb R$-complementary threshold for surfaces that is strictly less than $1$ is $\frac{12}{13}$. This result has many applications in explicit birational geometry of surfaces and threefolds and allows us to find several other optimal bounds on surfaces.

We show that the first gap of global log canonical threshold for surfaces is $\frac{1}{13}$, answering a question of V. Alexeev and W. Liu. We show that the minimal volume of log surfaces with reduced boundary and ample log canonical divisor is $\frac{1}{462}$, answering a question of J. Koll\'ar. We show that the smallest minimal log discrepancy (mld) of exceptional surfaces is $\frac{1}{13}$. As a special case, we show that the smallest mld of klt Calabi-Yau surfaces is $\frac{1}{13}$, reproving a recent result of L. Esser, B. Totaro, and C. Wang. After a more detailed classification, we classify all exceptional del Pezzo surfaces that are not $\frac{1}{11}$-lt, and show that the smallest mld of exceptional del Pezzo surfaces is $\frac{3}{35}$.

We also get better upper bounds of $n$-complements and Tian's $\alpha$-invariants for surfaces. Finally, as an analogue of our main theorem in high dimensions, we propose a question associating the gaps of $\mathbb R$-complementary thresholds with the gaps of mld's and study some special cases of this question.
\end{abstract}

\maketitle
\tableofcontents

\section{Introduction}

We work over the field of complex numbers $\mathbb C$. Nevertheless, since this paper mainly focuses on surfaces, many results in this paper should also hold for varieties over algebraically closed fields or non-closed fields with arbitrary characteristics by lifting (cf. \cite{BBKW22}).

\medskip

\noindent\textbf{First and second gaps of global log canonical thresholds of surfaces}. In this paper, global lc thresholds stand for values $b$, such that there exist lc Calabi-Yau pairs $(X,bS)$ for some non-zero effective Weil divisor $S$. (We remark that, in some references, Tian's $\alpha$-invariant is also called the ``global lc threshold". We will not use this notation in this paper.)

The main theorem of this paper shows that the first gap of global lc thresholds for surfaces is $\frac{1}{13}$. More precisely, the largest global lc threshold for surfaces that is strictly less than $1$ is $\frac{12}{13}$. This result answers a question of V. Alexeev and W. Liu \cite[Paragraph after Definition 1.5]{AL19a}. Moreover, we also show that the second gap of global lc thresholds for surfaces is $\frac{1}{11}-\frac{1}{13}$. That is, the largest global lc threshold for surfaces that is strictly less than $\frac{12}{13}$ is $\frac{10}{11}$.

\begin{thm}\label{thm: global lct 12/13 intro}
Let $b$ be a real number and $(X,bS)$ a klt log Calabi-Yau surface pair such that $S$ is a non-zero effective Weil divisor. Then $b\leq\frac{10}{11}$ or $b=\frac{12}{13}$.
\end{thm}

The following example is complementary to Theorem \ref{thm: global lct 12/13 intro}, which shows that the value $\frac{12}{13}$ is attainable. The value $\frac{10}{11}$ is also attainable by \cite[Table, Case 10]{Abe97}. We will prove later in Theorem  \ref{thm: global lct 12/13} that, if $X$ has Picard number $1$, then Example \ref{ex: 12/13} is the only case when $b=\frac{12}{13}$. 

\begin{ex}[{\cite[40]{Kol13},~\cite[Section 4]{AL19a}}]\label{ex: 12/13}
Consider the surface log pair $$\left(X,B=\frac{12}{13}S\right)=\left(\mathbb P(3,4,5),\frac{12}{13}(x^3y+y^2z+z^2x=0)\right).$$ Then $(X,B)$ is klt and $K_X+B\sim_{\mathbb Q}0$. 
\end{ex}
Example is also related to some special K3 surfaces. See Remark \ref{rem: k3 and example 12/13} below.

Theorem \ref{thm: global lct 12/13 intro} is a special case of the following theorem on the first and second gaps of $\Rr$-complementary thresholds for surfaces.

\medskip

\noindent\textbf{First and second gaps of $\mathbb R$-complementary thresholds of surfaces}. Recall that a projective pair $(X,B)$ is called \emph{$\Rr$-complementary} if there exists $B^+\geq B$ such that $(X,B^+)$ is lc and $K_X+B\sim_{\mathbb R}0$. The $\Rr$-complementary threshold ($\Rct$ for short) measures how far a pair is away from being not $\Rr$-complementary. $\Rr$-complementary thresholds can be viewed as the combination of two well-known thresholds, the lc thresholds and the anti-canonical thresholds \cite[Lemma 14]{Sho20}. It is known that the set of $\Rr$-complementary thresholds for (weak) Fano type pairs satisfies the ACC (\cite[Theorem 8.20]{HLS19}, \cite[Theorem 21]{Sho20}), and it is conjectured that the set of $\Rr$-complementary thresholds satisfies the ACC always (cf. \cite[Conjecture 8.19]{HLS19}). In particular, there should be a gap of (non-trivial) $\Rr$-complementary thresholds and $1$ in any dimension.

In this paper, we show that the first and the second gaps of $\Rr$-complementary thresholds for surfaces are equal to $\frac{1}{13}$ and $\frac{1}{11}$ respectively.

\begin{thm}[First and second gaps of $\Rr$-complementary thresholds]\label{thm: 12/13 rct intro}
Let $X$ be an $\Rr$-complementary projective surface and let $S$ be a non-zero effective Weil divisor on $X$. Let
$$b:=\Rct(X,0;S):=\sup\{t\geq 0\mid (X,tS)\text{ is }\Rr\text{-complementary}\}.$$
Then either $b=1$, or $b=\frac{12}{13}$, or $b\leq\frac{10}{11}$.
\end{thm}

Example \ref{ex: 12/13} and \cite[Table, Case10]{Abe97} also imply that the values $\frac{12}{13}$ and $\frac{10}{11}$ are attainable in Theorem \ref{thm: 12/13 rct intro}, even for rational del Pezzo surfaces. We refer the reader to Theorem \ref{thm: 12/13 rct} for a more precise version of Theorem \ref{thm: 12/13 rct intro}.

\smallskip

 Theorems \ref{thm: global lct 12/13 intro} and \ref{thm: 12/13 rct intro} have many important applications in explicit geometry of surfaces.

\medskip

\noindent\textbf{Log surface with the smallest volume}. By applying Theorem \ref{thm: global lct 12/13 intro}, we find the optimal lower bound of volumes of log surfaces with non-empty reduced boundary and ample log canonical divisor. The bound is $\frac{1}{462}$. Together with previously known examples, we answer a question of J. Koll\'ar \cite{Kol94,Kol13}. This optimal bound is known to have applications in moduli theory of surfaces: when considering the stable degeneration of surfaces of general type, lc surface pairs with reduced boundaries and ample log canonical divisors naturally appear as components of the normalization of special fibers with non-normal singularities. See Remark \ref{rem: connection with stable degeneration} (cf. \cite[Proposition 4.2]{Kol94}, \cite[40]{Kol13}) for details.

\begin{thm}\label{thm: vol 1/462}
Let $(X,S)$ be a projective lc surface pair such that $S$ is a non-zero effective Weil divisor and $K_X+S$ is ample. Then $\vol(X,S):=(K_X+S)^2\geq\frac{1}{462}$. Moreover, if $(K_X+S)^2=\frac{1}{462}$, then $S$ a non-singular rational curve. In particular, $S$ is a prime divisor.
\end{thm}

The following example shows that the bound in Theorem \ref{thm: vol 1/462} is optimal.

\begin{ex}[{\cite[Section 5]{AL19a}}]\label{ex: volume 1/462}
There exists a projective lc surface pair $(X,S)$ such that $S$ is a prime divisor, $K_X+S$ is ample, and $(K_X+S)^2=\frac{1}{462}$.
\end{ex}

It remains interesting to ask what is the smallest positive volume of lc (or klt) projective surfaces without boundary. The smallest positive volume for klt projective surfaces is expected to be $\frac{1}{48983}$ which was first given by V. Alexeev and W. Liu \cite[Theorem 1.4]{AL19a}. It is shown very recently by B. Totaro that the surface is a degree $438$ hypersurface in $\mathbb P(219,146,61,11)$ that is not quasi-smooth \cite[Introduction]{Tot23}. However, we don't know whether this value is optimal or not. There is no estimation on the smallest volume for lc but not klt projective surfaces, but we know the value is $\geq\frac{1}{86436}$ (\cite[Theorem 1.7(iii)]{AL19b}).

\medskip

\noindent\textbf{Smallest minimal log discrepancy (mld) of exceptional and Calabi-Yau surfaces}. We discover a correspondence between exceptional (Fano) varieties and R-complementary thresholds; see Section \ref{sec: 1/11 exceptional fano} for details. Using Theorem \ref{thm: 12/13 rct intro}, we can find the optimal lower bound of mld's of exceptional surfaces. Recall that a variety $X$ is called \emph{exceptional} if $|-K_X|_{\mathbb R}\not=\emptyset$ and $(X,D)$ is klt for any $D\in |-K_X|_{\mathbb R}$. By the theory of lc rational polytopes, $X$ is exceptional if and only if  $|-K_X|_{\mathbb Q}\not=\emptyset$ and $(X,D)$ is klt for any $D\in |-K_X|_{\mathbb Q}$.

\begin{thm}\label{thm: 1/13 exceptional}
Exceptional surfaces are $\frac{1}{13}$-lc.
\end{thm}

It is clear that any klt Calabi-Yau variety is exceptional. Therefore, as an immediate corollary of Theorem \ref{thm: 1/13 exceptional}, we reprove the following recent result of L. Esser, B. Totaro, and C. Wang.

\begin{cor}[{\cite[Proposition 6.1]{ETW22}}]\label{cor: calabi-yau surface mld}
Klt Calabi-Yau surfaces are $\frac{1}{13}$-lc.
\end{cor}

\begin{rem}\label{rem: k3 and example 12/13}
The bound $\frac{1}{13}$ in Corollary \ref{cor: calabi-yau surface mld} is optimal. Let $(X,B=\frac{12}{13}S)$ be the pair as in Example \ref{ex: 12/13}, $W\rightarrow X$ the extraction of all divisors with log discrepancy $1$ with respect to $(X,B)$, and $W\rightarrow X'$ the contraction of the strict transform of $S$ on $W$. Then $X'$ is a klt Calabi-Yau surface such that $\mld(X')=\frac{1}{13}$. As a consequence, the bound $\frac{1}{13}$ in Theorem \ref{thm: 1/13 exceptional} is optimal.

It is also interesting to notice that $W$ has a degree $13$ cover $\pi: \tilde W\rightarrow W$ ramified along the strict transform of $S$ on $W$, such that $\tilde W$ is a non-singular K3 surface associated with a non-symplectic automorphism action of order 13. The variety $\tilde W$ and its non-symplectic automorphism action is unique by \cite[Theorem 8.4]{AST11}. This matches our observation (cf. Theorem \ref{thm: global lct 12/13}) that Example \ref{ex: 12/13} is the only example when the global log canonical threshold is equal to $\frac{12}{13}$.
\end{rem}

Theorem \ref{thm: 1/13 exceptional} immediately implies that exceptional del Pezzo surfaces are $\frac{1}{13}$-lt. Moreover, since we are able to get a detailed classification of all cases when $\Rct=\frac{12}{13}$ is achieved (cf. Theorem \ref{thm: global lct 12/13}), we can classify all exceptional del Pezzo surfaces that are not $\frac{1}{11}$-lt. We have the following theorem:

\begin{thm}\label{thm: 1/13 exceptional Fano}
An exceptional del Pezzo surface is not $\frac{1}{11}$-lt if and only if it belongs to 25 types of Table \ref{table: 1/11-lt fano} below. In particular, we have the following.
\begin{enumerate}
    \item Exceptional del Pezzo surfaces that are not $\frac{1}{11}$-lt have isolated moduli, i.e. there are only 25 such surfaces up to isomorphism.
    \item  The smallest mld of exceptional del Pezzo surfaces is $\frac{3}{35}$.
    \item The smallest mld of exceptional del Pezzo surfaces of Picard number $1$ is $\frac{7}{79}$.
    \item Exceptional del Pezzo surfaces that are not $\frac{1}{11}$-lt are $13$-complementary.
\end{enumerate}
\end{thm}

\begin{rem}
We remark that Theorem \ref{thm: 1/13 exceptional Fano} also implies a classification of exceptional Fano type surfaces that are not $\frac{1}{11}$-lt. In particular:
\begin{enumerate}
    \item The smallest mld of exceptional Fano type surfaces is $\frac{3}{35}$ and is attained at a del Pezzo surface.
    \item Exceptional Fano type surfaces that are not $\frac{1}{11}$-lt are $13$-complementary.
\end{enumerate}
We refer the reader to Corollary \ref{cor: classification non-1/11-lt exceptional fano type} below for details.
\end{rem}

Theorem \ref{thm: 1/13 exceptional Fano} can be viewed as one key step towards a general classification of exceptional del Pezzo surfaces, which are known to be bounded \cite{Ale94,Bir19}. Moreover, since exceptional del Pezzo surfaces are always K-stable \cite{Tia87,OS12}, we expect Theorem \ref{thm: 1/13 exceptional Fano} to be useful in the moduli theory and the classification of K-stable surfaces. We remark that I. Cheltsov, J. Park, and C. Shramov have a classification of exceptional del Pezzo surfaces that are quasi-smooth well-formed hypersurfaces in weighted projective varieties \cite{CPS10}, but no general classification of exceptional surfaces is known in literature.

The exceptional del Pezzo surfaces constructed in Theorem \ref{thm: 1/13 exceptional Fano} seem to be new, except for some inexplicit examples in \cite{Lac20}. The old record of an exceptional del Pezzo surface with the smallest known mld was constructed in \cite[Lemma 3.1.17]{CPS10}, where the authors show that a degree $256$ hypersurface in $\mathbb P(11, 49, 69, 128)$ is an exceptional del Pezzo surface with $\mld=\frac{5}{49}>\frac{1}{11}$. Therefore, all exceptional del Pezzo surfaces constructed in Table \ref{table: 1/11-lt fano} are new as none of them is $\frac{1}{11}$-lt. In fact, the exceptional del Pezzo surfaces constructed in Theorem \ref{thm: 1/13 exceptional Fano} are related to varieties with $\Rr$-complementary thresholds equal to $\frac{12}{13}$, while a degree $256$ hypersurface in $\mathbb P(11,49,69,128)$ is related to varieties with $\Rr$-complementary thresholds equal to $\frac{10}{11}$. We also note that our construction can be used to construct more exceptional del Pezzo surfaces that are not necessarily $\frac{1}{11}$-lt. See Remark \ref{rem: some non-1/11-lt exc del pezzos} below for more details.

We remark that the exceptional del Pezzo surface with the smallest mld has the second smallest recorded anti-canonical volume; see Remark \ref{rem: vol 1/11 exceptional fano} below.

\medskip

\noindent\textbf{Upper bounds of $\alpha$-invariant and $n$-complements}. By adopting similar arguments as in \cite{Liu23}, Theorem \ref{thm: 1/13 exceptional Fano} will provide better explicit bounds on the boundedness of $n$-complements and Tian's $\alpha$-invariants of surfaces. These bounds are better than the bounds in \cite{Liu23} although they are still far from being optimal. 

\begin{cor}\label{cor: new bound alpha invariant}
Let $I_0:=2\cdot 11^{128\cdot 11^5}\approx 10^{10^{7.33}}$.
\begin{enumerate}
    \item For any exceptional del Pezzo surface $X$, $IK_X$ is Cartier for some $I<I_0$.
    \item Any $\Rr$-complementary surface $X$ is $n$-complementary for some $n<96I_0$.
    \item Tian's $\alpha$-invariant for any surface $X$ is $<3\sqrt{I_0}\approx 10^{10^{7.03}}$.
\end{enumerate}
\end{cor}
We remark that Corollary \ref{cor: new bound alpha invariant} also holds in the relative case with much stronger bounds. For example, for any morphism $\pi: X\rightarrow Z$ such that $\dim\pi(X)>0$ and $X/Z$ is $\Rr$-complementary, we have that $X/Z$ is $n$-complementary for some $n\in\{1,2,3,4,6\}$ if $\dim Z>0$ (\cite[3.1 Theorem]{Sho00}).

For all explicit geometry results we have proved above for surfaces, we refer to \cite{TW21,ETW22,ETW23,Tot23} for related results in higher dimensions, where many varieties with extremal properties (e.g. small mld, small volume, large $n$ for $n$-complement, large $\alpha$-invariant) are constructed. Unfortunately, although some values obtained by examples in these papers are conjecturally optimal, no optimal result is yet proved in dimension $\geq 3$.
\medskip

\noindent\textbf{A question on the gap of thresholds and the gap of mld's}. It is very interesting to notice that the first gap of $\Rr$-complementary thresholds in Theorem \ref{thm: 12/13 rct intro}, $\frac{1}{13}$, is also the first gap of global mld's with $1$ in dimension $3$ as proved in \cite[Theorem 1.4]{LX21a}. Moreover, it is easy to check that the only known $3$-fold example when $\mld=\frac{12}{13}$ is $\frac{1}{13}(3,4,5)$, which is exactly the orbifold cone $C(X,S)$ where $(X,S)$ is as in Example \ref{ex: 12/13}. These interesting behaviors indicate that there may be some deep relationship between $\Rr$-complementary thresholds and mld's. We first provide the notion of some gaps which naturally appear in birational geometry.

\begin{defn}[Gaps]\label{defn: gaps}
Let $d$ be a positive integer. We recall the definition of some folklore gaps. 
\begin{enumerate}
    \item $\delta_{\glct}(d)$ is the first gap of the global lc thresholds in dimension $d$, that is,
    \begin{align*}
    \delta_{\glct}(d):=1-\sup\Bigg\{t\Biggm|
    \begin{array}{r@{}l}
        (X,tS)\text{ is a klt pair of dimension } d,\\
        S\text{ is a non-zero effective Weil divisor, }K_X+tS\equiv 0 
    \end{array}\Bigg\},
    \end{align*}
    \item $\delta_{\tmld}(d+1)$ is the first gap of the global (total) mld's with $1$ in dimension $d+1$, that is,
    $$\delta_{\tmld}(d+1):=1-\sup\{\mld(X)\mid \mld(X)<1, \dim X=d+1\},$$
    \item $\delta_{\mld}(d+1)$ is the first gap of the local mld's with $1$ in dimension $d+1$, that is,
    $$\delta_{\mld}(d+1):=1-\sup\{\mld(X;x)\mid \mld(X;x)<1, \dim X=d+1\},$$
    \item $\delta_{\Rct}(d)$ is the first gap of $\Rr$-complementary thresholds in dimension $d$, that is,
    $$\delta_{\Rct}(d):=1-\sup\{t\mid t=\Rct(X,0;S)<1, \dim X=d, S\text{ is a non-zero effective Weil divisor}\},$$
    \item $\delta_{\Rct,\Ft}(d)$ is the first gap of the $\Rr$-complementary thresholds of Fano type varieties in dimension $d$, that is,
      \begin{align*}
    \delta_{\Rct,\Ft}(d):=1-\sup\Bigg\{t\Biggm|
    \begin{array}{r@{}l}
        X\text{ is Fano type of dimension }d,\\
       S\text{ is a non-zero effective Weil divisor, } t=\Rct(X,0;S)<1
    \end{array}\Bigg\},
    \end{align*}
        and
        \item $\delta_{\mld,\CY}(d)$ is the smallest mld of klt Calabi-Yau varieties in dimension $d$, that is,
        $$\delta_{\mld,\CY}(d):=\inf\{\mld(X)\mid \dim X=d, X\text{ is klt, }K_X\equiv 0\}.$$
\end{enumerate}
\end{defn}

\begin{ques}\label{ques: gap question}
Do we have
$$\delta_{\glct}(d)=\delta_{\tmld}(d+1)=\delta_{\mld}(d+1)=\delta_{\Rct}(d)=\delta_{\Rct,\Ft}(d)>0?$$
Even if the equations do not hold in general, can we get some inequalities between these values?
\end{ques}
It is obvious that $\delta_{\mld,\CY}(1)=1\not=\frac{1}{3}=\delta_{\glct}(1)$, so we do not consider $\delta_{\mld,\CY}(d)$ in Question \ref{ques: gap question}. Nevertheless, it is still interesting to consider the relationship between $\delta_{\mld,\CY}(d)$ and other invariants in higher dimensions.

\medskip

We hope Question \ref{ques: gap question} could help us to establish a global-to-local principle between $\Rr$-complementary thresholds in dimension $d$ and mld's in dimension $d+1$, and essentially help us to prove the ascending chain condition (ACC) conjecture for mld's \cite[Problem 5]{Sho88}. In fact, a positive answer for Question \ref{ques: gap question} will imply that $\delta_{\tmld}(d)$ and $\delta_{\mld}(d)$ are positive, and will deduce the $1$-gap conjecture (cf. \cite[Conjecture 1.2]{Jia21}) for mld's, a special case of the ACC conjecture for mld's. We summarize some partial results towards Question \ref{ques: gap question} in the following theorem.

\begin{thm}\label{thm: summary gap conjecture}
Let $d$ be a positive integer.
\begin{enumerate}
    \item Question \ref{ques: gap question} have a positive answer when $d=1$. We have $\delta_{\glct}(1)=\frac{1}{3}$.
    \item $\delta_{\glct}(2)=\delta_{\tmld}(3)=\delta_{\Rct}(2)=\delta_{\Rct,\Ft}(2)=\delta_{\mld,\CY}(2)=\frac{1}{13}\geq\delta_{\mld}(2)$.
    \item $\delta_{\mld,\CY}(d)\geq\delta_{\glct}(d)>0$, $\delta_{\tmld}(d+1)\geq\delta_{\tmld}(d+1)$, and $\delta_{\Rct,\Ft}(d)\geq\delta_{\Rct}(d)$.
    \item $\delta_{\glct}(d)\geq\delta_{\mld}(d+1)$ and $\delta_{\glct}(d)\geq\delta_{\Rct,\Ft}(d)$.
\end{enumerate}
\end{thm}

Finally, we remark that the set of global lc thresholds in dimension $d$ and the set of mld's in dimension $d+1$ may share more common properties than Question \ref{ques: gap question} expects. For example, we know that the largest accumulation point of global lc thresholds in dimension $2$ is $\frac{5}{6}$ (cf. \cite[Theorem 1.2]{Pro01} and \cite[Propositions 11.5, 11.7]{HMX14}), which is also equal to the second largest accumulation point of mld's in dimension $3$ (cf. \cite[Theorem 1.3]{LL22}). 

\medskip

\noindent\textit{Sketch of the proofs of the main theorems}. We first prove Theorem \ref{thm: global lct 12/13 intro}. By applying results in \cite{Sho00} and \cite{Abe97}, we may reduce the question to the case when $X$ is a rational del Pezzo surface of Picard number $1$, $S$ is a non-singular rational curve, $(X,S)$ is plt, and $(X,B=bS)$ is $\frac{1}{7}$-lt at closed points. By Bogomolov bounds (cf. \cite[9.2 Corollary]{KM99}) and \cite{Kol94}, we know that $b\leq\frac{41}{42}$ and the order of any singularity outside $S$ is $\leq 42$. This is done in Section \ref{sec: set-up}. 

With these properties settled and together with the extra restrictions on the indices of singularities on $X$, we can classify all possible singularities that may appear on $S$. Therefore, we can classify all possible singularities that may appear on $X$. This is done in Section \ref{sec: sings on S} (see Table \ref{Table: singularity on S}). At the mean time, we know that $X$ has at most $4$ singularities and $S$ contains at least $3$ singularities of $X$. Thus we get a big but finite list of singularities on $X$ for us to take care of. Now a key observation is that, by using the fact that $\rho(\tilde X)+K_{\tilde X}^2=10$ where $\tilde X$ is the minimal resolution of $X$, $b$ will satisfy a quadratic equation
$$b^2+sb+t=0$$
(cf. Lemma \ref{lem: gamma(X) equality}) where $s,t$ are rational numbers which only depend on the singularities of $X$ (and their behavior along $S$). Since $b$ is a rational number, $s^2-4t$ is a square of a rational number. To provide a detailed description of $s$ and $t$, we introduce two invariants of surface singularities, $\gamma(X\ni x)$ and $q(X\ni x,S)$. This is done in Section \ref{sec: control singularities}.

It turns out that for random choices of triples or quadruples of singularities on $X$, $s^2-4t$ is usually not a square of a rational number. We also have several other restrictions on the singularities of $X$ as well, which are also stated in Section \ref{sec: control singularities}. After running a computer program, we reduce the range of $(X,S)$ to a very small list of possible cases. Now we can treat all cases one by one, and prove Theorem  \ref{thm: global lct 12/13 intro}. Theorems \ref{thm: 12/13 rct intro} and \ref{thm: 1/13 exceptional} can be easily deduced from Theorem \ref{thm: global lct 12/13 intro}, and Corollary \ref{cor: calabi-yau surface mld} follows from Theorem \ref{thm: 1/13 exceptional}. Moreover, we can additionally show that all exceptional varieties with a non-$\frac{1}{11}$-lc complement are dominated by the terminalization of the log pair in Example \ref{ex: 12/13}. Theorem \ref{thm: 1/13 exceptional Fano} follows from a detailed classification of all such varieties (see Section \ref{sec: 1/11 exceptional fano}). Corollary \ref{cor: new bound alpha invariant} follows from Theorem \ref{thm: 1/13 exceptional Fano} and similar arguments as in \cite{Liu23}. 

Now we prove Theorem \ref{thm: vol 1/462}. It is not difficult to reduce to the case when $S$ is a prime divisor. In this case, $\vol(X,S)=(K_X+S)^2$ is controlled by the pseudo-effective threshold $\pet(X,0;S)$ (cf. Definition \ref{defn: pet}) and the intersection number $(K_X+S)\cdot S$. By adjunction, we know that $(K_X+S)\cdot S=\frac{1}{42}$ or $(K_X+S)\cdot S\geq\frac{1}{24}$. By Theorem \ref{thm: global lct 12/13 intro}, $\pet(X,0;S)$ is either $\frac{12}{13}$ or $\leq\frac{10}{11}$. It is then easy to see that we have the desired inequality unless $\pet(X,0;S)=\frac{12}{13}$ and $(K_X+S)\cdot S=\frac{1}{42}$. This tells us that $S$ is a non-singular rational curve which contains three singular points of $X$, one of which has order $7$. In particular, $\mld(X,S)=\frac{1}{7}$. Now we may run a $(K_X+(\frac{12}{13}-\epsilon)S)$-MMP for some $0<\epsilon\ll 1$. It is easy to check that this MMP will terminate with a del Pezzo surface $X'$ such that $\rho(X')=1$. In this case, $(X',\frac{12}{13}S')$ must be isomorphic to the log pair as in Example \ref{ex: 12/13} where $S'$ is the image of $S$ on $X'$. In particular, $\mld(X',S')$ is $\frac{1}{5}$-lc. This contradicts the fact that $X\rightarrow X'$ is $(K_X+S)$-non-negative. Theorem \ref{thm: vol 1/462} follows immediately.

Finally, we say a few words for Theorem \ref{thm: summary gap conjecture}. The key part of this theorem is the inequality $\delta_{\glct}(d)\geq\delta_{\mld}(d+1)$ in Theorem \ref{thm: summary gap conjecture}(4). We consider a log pair $(X,bS)$ of dimension $d$ such that $S$ is a non-zero effective Weil divisor, and $b=1-\delta_{\glct}(d)$. By induction on dimensions, we may reduce to the case when $\rho(X)=1$. Now the cone $C(X,S)$ of $X$ has $\mld=1-\delta_{\glct}(d)$ at the vertex of the cone. This construction implies the desired inequality.

\medskip

\noindent\textbf{Acknowledgement}. The first author would like to thank Chen Jiang for telling him the question in MSRI in Spring 2019. He would like to thank Nathan Chen, Christopher D. Hacon, Jingjun Han, Yuchen Liu, Yujie Luo, Fanjun Meng, Lu Qi, Lingyao Xie, Chenyang Xu, Qingyuan Xue, and Ziquan Zhuang for useful discussions. He would like to thank Louis Esser, Burt Totaro, and Chengxi Wang for useful communications on their paper \cite{ETW22}. He would like to thank Yujie Luo for helping him coding. We thank Justin Lacini for helping us identifying two surfaces in our classification with the surfaces enlisted in \cite{Lac20}.

\section{Preliminaries}\label{sec: preliminaries}

We adopt the standard notation and definitions in \cite{Sho92,KM98,BCHM10} and will freely use them. A \emph{surface} is a normal quasi-projective variety of dimension $2$.

\begin{agree}\label{defn sing}
A \emph{pair} $(X/Z\ni z,B)$ consists of a \emph{contraction} $\pi: X\rightarrow Z$ between normal quasi-projective varieties, a (not necessarily closed) point $z\in Z$, and an $\Rr$-divisor $B\geq 0$ on $X$. If $K_X+B$ is $\Rr$-Cartier over a neighborhood of $z$, then we say that $(X/Z\ni z,B)$ is a \emph{log pair}. If $(X/Z\ni z,B)$ is a (log) pair for any point $z\in Z$, then we say that $(X/Z,B)$ is a (log) pair. If $\pi$ is the identity morphism and $(X/Z\ni z,B)$ (resp. $(X/Z,B)$) is a (log) pair, then we denote it by $(X\ni z,B)$ (resp. $(X,B)$). If $(X\ni x,0)$ is a log pair and $x$ is a closed point, then we say that $X\ni x$ is a \emph{germ} of $X$ at $x$. The \emph{order} of a germ $X\ni x$ is defined as the order of its local algebraic fundamental group.
\end{agree}

\begin{agree}
Let $(X/Z\ni z,B)$ be a pair and $\pi: X\rightarrow Z$ its contraction. For any prime divisor $E$ over $X$, we say that $E$ is \emph{over} $z$ if $\bar z=\pi(\Center_XE)$.

Suppose that $(X/Z\ni z,B)$ is a log pair. Let $f: Y\rightarrow X$ be a log resolution of $(X,B)$ and $K_Y+B_Y:=f^*(K_X+B)$ over a neighborhood of $z$. For any prime divisor $E$ over $X/Z\ni z$ on $Y$, we recall that $a(E,X,B):=1-\mult_EB_Y$ denotes the \emph{log discrepancy} of $E$ with respect to $(X,B)$. We denote by
 $$\mld(X,B;z):=\inf\{a(E,X,B)\mid E\text{ is over }z\}$$
the \emph{minimal log discrepancy} (\emph{mld}) of $(X/Z\ni z,B)$. If $(X,B)$ is a log pair, then we denote by
 $$\mld(X,B):=\inf\{a(E,X,B)\mid E\text{ is an exceptional prime divisor of }X\}$$
the \emph{mld} of $(X,B)$, and
  $$\tmld(X,B):=\inf\{a(E,X,B)\mid E\text{ is an exceptional or non-exceptional prime divisor over }X\}$$
the \emph{total minimal log discrepancy} (\emph{tmld}) of $(X,B)$.
 
 Let $\epsilon$ be a non-negative real number. A pair $(X/Z\ni z,B)$ is called lc (resp. klt, $\epsilon$-lc, $\epsilon$-lt) if $(X,B;z)$ is a log pair and $\mld(X,B;z)\geq 0$ (resp. $>0,\geq\epsilon,>\epsilon$). If $\pi$ is the identity morphism, then we say that $(X\ni z,B)$ is lc (resp. klt, $\epsilon$-lc, $\epsilon$-lt) if $(X/Z\ni z,B)$ is lc. If $(X,B)$ is a log pair, then we say that $(X,B)$ is lc (resp. klt, $\epsilon$-lc, $\epsilon$-lt) if $\tmld(X,B)\geq 0$ (resp. $>0,\geq\epsilon,>\epsilon$).
 
 A pair $(X,B)$ is called \emph{plt} if $(X\ni x,B)$ is klt for any codimension $\geq 2$ point $x\in X$. We say that $(X\ni x,B)$ is \emph{plt} if $(X,B)$ is plt near $x$. 
 
 A pair $(X,B)$ is called \emph{dlt} if there exists a log resolution $f: Y\rightarrow X$ of $(X,\Supp B)$ such that $a(E,X,B)>0$ for any prime divisor $E$ on $Y$ that is exceptional over $X$. 
 
 For any lc pair $(X,B)$, a \emph{dlt modification} of $(X,B)$ is a projective birational morphism $f: Y\rightarrow X$ which only extracts divisors with log discrepancy $0$ with respect to $(X,B)$, and $(Y,B_Y)$ is $\Qq$-factorial dlt, where $K_Y+B_Y:=f^*(K_X+B)$. We recall that any lc pair has a dlt modification (cf. \cite[Proposition 3.3.1]{HMX14}), and for any dlt surface pair $(X,B)$, $X$ is $\Qq$-factorial.
 
 We remark that any lc (resp. klt, $\epsilon$-lc, $\epsilon$-lt, plt, dlt) pair is automatically a log pair.
 \end{agree} 
 
 \begin{defn}
A pair $(X,B)$ is called \emph{log Calabi-Yau} if $K_X+B\equiv 0$. A variety $X$ is called \emph{Calabi-Yau} if $(X,0)$ is log Calabi-Yau. A pair $(X,B)$ is called \emph{exceptional} if $|-(K_X+B)|_{\mathbb R}\not=\emptyset$, and for any $D\in |-(K_X+B)|_{\mathbb R}$, $(X,B+D)$ is klt. $(X,B)$ is called \emph{log Fano} if $-(K_X+B)$ is ample and $(X,B)$ is klt.  We remark that we will only consider projective log Calabi-Yau (resp. exceptional, log Fano) pairs in this paper.
 
 A variety $X$ is called \emph{of Fano type} if $(X,\Delta)$ is log Fano for some $\Delta$.  A variety $X$ is called \emph{Fano} if $-K_X$ is ample and $X$ is klt. A Fano surface is also called a \emph{del Pezzo} surface.
 
 We denote by $\rho(X)$ the Picard number of $X$ for any variety $X$.
 \end{defn}
 
 \begin{cons}\label{cons: terminalization}
 Let $(X,B)$ be a klt pair. Then there exists a birational morphism $f: Y\rightarrow X$ and a klt pair $(Y,B_Y)$, such that $K_Y+B_Y:=f^*(K_X+B)$ and $\mld(Y,B_Y)>1$. When $X$ is a surface, $f$ is unique, and we will call $f$ as \emph{the terminalization} of $(X,B)$.
 \end{cons}
 
 \begin{defn}
Let $(X/Z\ni z,B)$ and $(X/Z\ni z,B^+)$ be two pairs. We say that $(X/Z\ni z,B^+)$ is an \emph{$\Rr$-complement} of $(X/Z\ni z,B)$ if $B^+\geq B$, $(X,B^+)$ is lc over a neighborhood of $z$, and $K_X+B^+\sim_{\mathbb R}0$ over a neighborhood of $z$. We say that $(X/Z,B^+)$ is an \emph{$\Rr$-complement} of $(X/Z,B)$ if $(X/Z\ni z,B^+)$ is an $\Rr$-complement of $(X/Z\ni z,B)$ for any $z\in Z$. We say that $(X/Z\ni z,B)$ (resp. $(X/Z,B)$) is \emph{$\Rr$-complementary} if $(X/Z\ni z,B)$ (resp. $(X/Z,B)$) has an $\Rr$-complement. For any $\Rr$-divisor $D\geq 0$ on $X$, we define
 $$\Rct(X/Z\ni z,B;D):=\sup\{t\geq 0\mid (X/Z\ni z,B+tD)\}\text{ is }\Rr\text{-complementary}\}$$
  $$(\text{resp. }\Rct(X/Z,B;D):=\sup\{t\geq 0\mid (X/Z,B+tD)\}\text{ is }\Rr\text{-complementary}\})$$
  to be the $\Rr$-complementary threshold of $D$ with respect to $(X/Z\ni z,B)$ (resp. $(X/Z,B)$).
  
  If $Z$ is a point, then we may drop $Z\ni z$ or $Z$.
 \end{defn}
 
 \begin{defn}
 Let $(X/Z\ni z,B)$ and $(X/Z\ni z,B^+)$ be two pairs. We say that $(X/Z\ni z,B^+)$ is an \emph{$n$-complement} of $(X/Z\ni z,B)$ if $nB^+ \geq n\lfloor B\rfloor+\lfloor (n+1)\{B\}\rfloor$, $(X,B^+)$ is lc over a neighborhood of $z$, and $n(K_X+B^+)\sim 0$ over a neighborhood of $z$. We say that $(X/Z\ni z,B)$ is \emph{$n$-complementary} if $(X/Z\ni z,B)$ has an $n$-complement.
 \end{defn}
 
 \begin{defn}
 Let $(X/Z\ni z,B)$ (resp. $(X/Z,B)$) be a log pair and $D\geq 0$ an $\Rr$-Cartier $\Rr$-divisor on $X$. We define
 $$\lct(X/Z\ni z,B;D):=\sup\{t\mid (X/Z\ni z,B+tD)\text{ is lc}\}$$
$$(\text{resp. }\lct(X/Z,B;D):=\sup\{t\mid (X/Z,B+tD)\text{ is lc}\})$$
to be the \emph{lc threshold} of $D$ with respect to $(X/Z\ni z,B)$ (resp. $(X/Z,B)$). If $Z$ is a point, then we may drop $Z\ni z$ or $Z$.
 
 If $(X,B)$ is a projective pair such that $|-(K_X+B)|_{\mathbb R}\not=\emptyset$, then we define
 $$\alpha(X,B):=\inf\{\lct(X,B;D)\mid D\in |-(K_X+B)|_{\mathbb R}\}$$
to be \emph{Tian's $\alpha$-invariant} of $(X,B)$. If $X$ is $\Qq$-Gorenstein, then we define $\alpha(X):=\alpha(X,0)$ to be \emph{Tian's $\alpha$-invariant} of $X$. 
 \end{defn}
 
 \begin{defn}
 Let $n$ be a positive integer and $X$ a surface. A \emph{$(-n)$-curve} on $X$ is a smooth rational curve $C$ on $X$ such that $C^2=-n$.
 \end{defn}

\begin{defn}[Dual graph]\label{defn: dual graph}
Let $n$ be a non-negative integer, and $C=\cup_{i=1}^nC_i$ a collection of irreducible curves on a non-singular surface. We define the \emph{dual graph} $\Dd(C)$ of $C$ as follows.
\begin{enumerate}
    \item The vertices $v_i=v_i(C_i)$ of $\Dd(C)$ correspond to the curves $C_i$.
    \item For $i\neq j$, the vertices $v_i$ and $v_j$ are connected by $C_i\cdot C_j$ edges.
    \item Each vertex $v_i$ is labeled by the integer $e_i:=-C_i^2$. Nevertheless, if $e_i$ is not important or unknown, then we may not write it down in the dual graph.
\end{enumerate}
We sometimes write the notation of the curve $C_i$ near the vertex $v_i$. We may highlight some vertices by using filled circles. We say that $\mathcal{D}(C)$ contains a \emph{cycle} if there exists $C_{i_1},\dots,C_{i_l}$ for some $l\geq 3$ such that $C_{i_j}$ intersects $C_{i_k}$ when $|j-k|\leq 1$ or $\{j,k\}=\{1,l\}$. We say that $\mathcal{D}(C)$ is a \emph{tree} if
\begin{itemize}
    \item $\mathcal{D}(C)$ does not contain a cycle, and
    \item $C_i\cdot C_j\leq 1$ for any $i\not=j$.
\end{itemize}
The \emph{intersection matrix} of $\mathcal{D}(C)$ is defined as the matrix $(E_i\cdot E_j)_{1\leq i,j\leq n}$ if $C\not=\emptyset$. The \emph{determinant} of $\mathcal{D}(C)$ is defined as
$$\det(\mathcal{D}(C)):=\det(-(E_i\cdot E_j)_{1\leq i,j\leq n})$$
if $C\not=\emptyset$, and $\det(\mathcal{D}(C)):=1$ if $C=\emptyset$.

A \emph{fork} of $\mathcal{D}(C)$ is a curve $C_i$ such that $C_i\cdot C_j\geq 1$ for at least three different $j\not=i$. A \emph{tail} of $\mathcal{D}(C)$ is a curve $C_i$ such that $C_i\cdot C_j\geq 1$ for at most one $j\not=i$.

For any projective birational morphism $f: Y\rightarrow X$ between surfaces, let $E_1,\dots,E_n$ be the set of prime $f$-exceptional divisors, and let $E=\cup_{i=1}^nE_i$ be the reduced exceptional divisor with $n$ prime components. We define $\Dd(f):=\Dd(E)$. 
\end{defn}

\section{Set-up and reductions}\label{sec: set-up}

In this section, we make some preparations to prove Theorem \ref{thm: global lct 12/13 intro}. In fact, we prove a more general and precise version of Theorem \ref{thm: global lct 12/13 intro}.

\begin{thm}\label{thm: global lct 12/13}
Let $m$ be a positive integer, $b$ a real number, and $(X,B:=T+\sum_{i=1}^m b_iS_i)$ a complete lc surface pair, such that 
\begin{itemize}
    \item $T,S_1,\dots,S_m$ are effective Weil divisors and $S:=\sum_{i=1}^m S_i\not=0$,
    \item $0<b\leq b_1,\dots,b_m<1$, and
    \item $(X,B)$ is log Calabi-Yau.
\end{itemize}
Then:
\begin{enumerate}
    \item $b\leq\frac{12}{13}$.
    \item If $b>\frac{6}{7}$, then $T=0$, $m=1$, and $S$ is prime.
    \item If $b\geq\frac{10}{11}$, then $b_1=\frac{10}{11}$ or $\frac{12}{13}$. In addition, if $b_1=\frac{12}{13}$, then:
    \begin{enumerate}
        \item If $\rho(X)=1$, then
        $(X,B)$ is isomorphic to $$(X_0,B_0):=\left(\mathbb P(3,4,5),\frac{12}{13}(x^3y+y^2z+z^2x=0)\right).$$
        and the isomorphism $(X,B)\rightarrow (X_0,B_0)$ is unique. In particular, $(X,B)$ is $\frac{3}{13}$-lc at any closed point of $X$.
        \item $(X,B)$ is a crepant model of $(X_0,B_0)$. More precisely, there exists a projective birational morphism $X\rightarrow X_0$ which only extract divisors with log discrepancy $1$ with respect to $(X_0,B_0)$, and $K_X+B$ is the pullback of $K_{X_0}+B_0$.
    \end{enumerate}
\end{enumerate}
\end{thm}

The goal of this section is to reduce Theorem \ref{thm: global lct 12/13} to the case when $\rho(X)=1$, $T$ is $0$, $m=1$, $S$ is a prime divisor, $b_1=b$, $(X,B)$ is klt, and the singularities on $X$ are well-controlled. In particular, we will prove Theorem \ref{thm: global lct 12/13}(2) in this section (see Lemma \ref{lem: reduce to smooth rational curve}). We summarize the conditions of Theorem \ref{thm: global lct 12/13} in the following set-up:

\begin{setup}\label{setup: 12/13}
$m,b,X,B,T,S,b_1,\dots,b_m,S_1,\dots,S_m$ are as follows:
\begin{enumerate}
\item $m$ is a positive integer.
\item $0<b\leq b_1,\dots,b_m<1$.
    \item $(X,B=T+\sum_{i=1}^mb_iS_i)$ is a complete lc surface pair.
    \item $T,S_1,\dots,S_m$ are effective Weil divisors and $S:=\sum_{i=1}^m S_i\not=0$.
    \item $K_X+B\equiv 0$.
\end{enumerate}
\end{setup}

\subsection{Reduction to the klt case}

In this subsection, we reduce Theorem \ref{thm: global lct 12/13} to the case when $(X,B)$ is klt and $(X,S)$ is lc. In particular, $X$ is projective and $T=0$.

\begin{lem}\label{lem: nonklt 5/6 lemma}
Conditions as in Set-up \ref{setup: 12/13}. Suppose that there exists a component $C$ of $T$ such that $C$ intersects $S$. Then $b\leq\frac{5}{6}$.
\end{lem}
\begin{proof}
Let $f: Y\rightarrow X$ be a dlt modification of $(X,B)$. Let $K_Y+B_Y:=f^*(K_X+B)$, $S_{Y_i}:=f^{-1}_*S_i$ for any $i$, $S_Y:=\sum_{i=1}^mS_i$, and $T_Y$ the sum of the strict transform of $T$ on $Y$ and the reduced $f$-exceptional divisor. Then there exists a component $C_Y$ of $T_Y$ which intersects $S_Y$. Possibly replacing $X,B,T,S,S_1,\dots,S_m$ with $Y,B_Y,T_Y,S_Y,S_{Y_1},\dots,S_{Y_m}$ respectively, we may assume that $(Y,B_Y)$ is $\Qq$-factorial dlt. In particular, $C$ is normal.

Suppose that $b>\frac{5}{6}$. We consider 
$$K_C+B_C:=(K_X+B)|_C.$$
Since $(X,B)$ is dlt, we have
$$B_C=\sum_j\frac{n_j-1+\sum_i k_{i,j}b_i}{n_j}D_j+P,$$
where $D_j$ are distinct points, $k_{i,j}\in\mathbb N$, $n_j\in\mathbb N^+$, at least one $k_{i,j}$ is positive, $$0<n_j-1+\sum_ik_{i,j}b_i<1$$ for any $j$, and $P=(\lfloor B\rfloor-C)|_C$ is a sum of distinct points. Since $K_C+B_C\equiv 0$, we have $C=\mathbb P^1$ and $$2=\deg B_C=\sum_j \frac{n_j-1+\sum_i k_{i,j}b_i}{n_j}+\deg P.$$
Thus $\sum_j \frac{n_j-1+\sum_i k_{i,j}b_i}{n_j}\in\{1,2\}$. Since $\frac{5}{6}<b\leq b_i<1$ for all $i$, an elementary computation shows that it is not possible.
\end{proof}

\begin{lem}\label{lem: mfs 2/3}
Conditions as in Set-up \ref{setup: 12/13}. Suppose that there exists a contraction $f: X\rightarrow Z$ such that $Z$ is a curve and $S_i$ dominates $Z$ for some $i$. Then $b\leq\frac{2}{3}$.
\end{lem}
\begin{proof}
Let $F$ be a general fiber of $f$. We consider
$$K_F+B_F:=(K_X+B)|_F, T_F:=T|_F,\text{ and }S_{i,F}:=S_i|_F\text{ for each }i.$$
Then $S_{i,F}\not=0$ for some $i$ and $K_F+B_F\equiv 0$. Thus $F=\mathbb P^1$, hence
$$2=\deg B_F=\deg T_F+\sum b_i\deg S_{i,F}.$$
Since $1>b_i\geq b$ for all $i$, this equation cannot have a solution when $b>\frac{2}{3}$.
\end{proof}

\begin{prop}\label{prop: not klt case}
Conditions as in Set-up \ref{setup: 12/13}. If $(X,B)$ is not klt, then $b\leq\frac{5}{6}$. In particular, if $b>\frac{5}{6}$, then $X$ is klt and projective.
\end{prop}
\begin{proof}
Suppose that  $(X,B)$ is not klt. Let $f: Y\rightarrow X$ be a dlt modification of $(X,B)$, $K_Y+B_Y:=f^*(K_X+B)$, $S_Y$ the strict transform of $S$ on $Y$, $S_{i,Y}$ the strict transform of $S_i$ on $Y$ for each $i$, and $T_Y:=\lfloor B_Y\rfloor$. Then $T_Y, S_{i,Y}$ are effective Weil divisors, $S_Y=\sum_{i=1}^mS_{i,Y}$, and $(Y,B_Y=T_Y+\sum_{i=1}^ma_iS_{i,Y})$ is $\mathbb Q$-factorial dlt log Calabi-Yau. Possibly replacing $X,B,T,S,S_i$ with $Y,B_Y,T_Y,S_Y,S_{i,Y}$ respectively, we may assume that $(X,B)$ is $\mathbb Q$-factorial dlt. In particular, $T=\lfloor B\rfloor\not=0$.

Possibly reordering indices, we may assume that $S_1\not=0$. We may run a $(K_X+B-b_1S_1)$-MMP which terminates with a Mori fiber space. Then this MMP is $S_1$-positive. Hence $S_1$ is not contracted by this MMP. There are three cases.

\medskip

\noindent\textbf{Case 1}. A component of $T$ is contracted by a step $X'\rightarrow X''$ of this MMP. Possibly replacing $X$ with $X'$ and replacing $B,T,S_i,S$ with their images on $X'$ respectively, we may assume that $T$ intersects $S_1$, hence $T$ intersects $S$. The proposition follows from Lemma \ref{lem: nonklt 5/6 lemma}.

\medskip

\noindent\textbf{Case 2}. Any component of $T$ is not contracted by this MMP. In this case, we let $X\rightarrow X'$ be the contraction, where $X'$ is the resulting model of this MMP with $X'\rightarrow Z$ the Mori fiber space. Possibly replacing $X$ with $X'$ and replacing $B,T,S,S_i$ with their images on $X'$ respectively, we may assume that we have a $(K_X+B-b_1S_1)$-Mori fiber space $X\rightarrow Z$. In particular, $S_1$ dominates $Z$. By Lemma \ref{lem: mfs 2/3}, we may assume that $Z$ is a point. Thus $\rho(X)=1$ and $T$ intersects $S$. The proposition follows from Lemma \ref{lem: nonklt 5/6 lemma}.
\end{proof}

The following lemma follows from the precise classification of surface lc thresholds \cite{Kuw99}. A precise statement was provided in \cite[Corollary 3.3]{Pro02}.
\begin{lem}[{cf. \cite{Kuw99}, \cite[Corollary 3.3]{Pro02}}]\label{lem: surface lct 5/6}
Let $(X,T)$ be an lc pair such that $T$ is a Weil divisor, and $S$ is a non-zero effective $\Qq$-Cartier Weil divisor. Then $\lct(X,T;S)=1$ or $\lct(X,T;S)\leq\frac{5}{6}$.
\end{lem}

\begin{lem}\label{lem: (X,T+S) lc}
Conditions as in Set-up \ref{setup: 12/13}. If $b>\frac{5}{6}$, then $T=0$ and $(X,S)$ is lc.
\end{lem}
\begin{proof}
By Proposition \ref{prop: not klt case}, $T=0$. The lemma follows Lemma \ref{lem: surface lct 5/6}.
\end{proof}

\begin{lem}\label{lem: reduced to picard number 1}
Conditions as in Set-up \ref{setup: 12/13}. Suppose that $(X,B)$ is klt and $b>\frac{2}{3}$. Then there exists a contraction $f: X\rightarrow Y$ and a klt surface pair $(Y,B_Y=\sum_{i=1}^mb_iS_{i,Y})$, such that
\begin{enumerate}
\item $K_X+B=f^*(K_Y+B_Y)$,
    \item for each $i$, $S_{i,Y}=f_*S_{i}$ is an effective Weil divisor for each $i$, and $S_Y:=\sum b_iS_{i,Y}\not=0$,
    \item $\rho(Y)=1$ and $K_Y+B_Y\sim_{\mathbb R}0$, and
    \item if $b>\frac{5}{6}$, then $(Y,S_Y)$ is lc. 
\end{enumerate}
\end{lem}
\begin{proof}
By Set-up \ref{setup: 12/13}, we may assume that $S_1\not=0$. We run a $(K_X+B-b_1S_1)$-MMP, which terminates with a Mori fiber space $Y\rightarrow Z$.  We let $S_{i,Y}$ is the image of $S_i$ on $Y$ for each $i$. Then $S_{1,Y}$ dominates $Z$. In particular, $S_{1,Y}\not=0$. 

Since $K_X+B\equiv 0$, we get (1). Since $S_{1,Y}\not=0$, $S_Y\not=0$, and we get (2). By Lemma \ref{lem: mfs 2/3}, $Z$ is a point. Hence $\rho(Y)=1$ and we get (3). (4) follows from Lemma \ref{lem: (X,T+S) lc}.
\end{proof}

\subsection{Reduction to an irreducible boundary}

In this subsection, we reduce to the case when $S$ is a non-singular rational curve and $(X,S)$ is plt. In particular, $S$ is irreducible and $S\cap\Sing(X)$ only contains cyclic quotient singularities of $X$.

\begin{nota}
Let $(X,B)$ be a klt pair. We put
  \begin{align*}
    \delta(X,B):=\#\left\{E\biggm| \begin{array}{r@{}l} a(E,X,B)\leq\frac{1}{7}, E  
        \text{ is an exceptional}\\
       \text{or non-exceptional prime divisor over }X
    \end{array}\right\}.
    \end{align*}
\end{nota}

\begin{lem}[cf. {\cite[Theorem 4.15]{KM98}}]\label{lem: de type mld}
Let $(X\ni x,S)$ be an lc surface pair such that $X\ni x$ is a $D$-type or $E$-type singularity and $S$ is an effective Weil divisor. Then $X\ni x$ is $D$-type singularity, $(X\ni x,S)$ is not plt. Moreover, let $E$ be the  prime divisor corresponding to the unique fork of the minimal resolution of $X\ni x$, then $a(E,X,S)=0$.
\end{lem}

\begin{lem}\label{lem: reduce to smooth rational curve}
Conditions as in Set-up \ref{setup: 12/13}. Suppose that $b>\frac{6}{7}$. Then:
\begin{enumerate}
    \item $\delta(X,B)=1$. In particular, $T=0$, $m=1$, $S$ is prime, and $(X,B)$ is $\frac{1}{7}$-lt at any closed point $x\in X$.
    \item For any point $x\in S$ that is contained in the singular locus of $X$ but is not contained in the singular locus of $S$, $x$ is a cyclic quotient singularity, and $(X,S)$ is plt near $x$.
    \item One of the following holds:
    \begin{enumerate}
        \item $S$ is a non-singular rational curve.
        \item $S$ is a non-singular elliptic curve.
        \item $S$ is a singular rational curve with a unique singularity $x$, such that $x$ is a nodal singularity of $S$ and a non-singular point of $X$.
    \end{enumerate}
    \item If $b>\frac{10}{11}$, then $S$ is a non-singular rational curve. In particular, $(X,S)$ is plt.
\end{enumerate}
\end{lem}
\begin{proof}
By Proposition \ref{prop: not klt case}, $(X,B)$ is klt. Since $S\not=0$, $\delta(X,B)\geq 1$.

(1) By Lemma \ref{lem: reduced to picard number 1}, we may assume that $\rho(X)=1$. Now (1) follows from the classification result \cite[5.1.3]{Sho00} for the case when $\delta(X,B)=2$. To see this, note that $F=\{D\}$ as in \cite[5.1.1]{Sho00} must be $0$ as the coefficients of $B$ are strictly greater than $\frac{6}{7}$. Hence only \cite[Case (A$_2^6$)]{Sho00} could happen. In this case, the coefficients of $B$ are equal to $\frac{6}{7}$, which is also not possible.

(2) Let $x$ be a singular point of $X$ that is contained in $S$. By Lemma \ref{lem: (X,T+S) lc}, $(X,S)$ is lc near $x$. If $x$ is a cyclic quotient singularity of $X$, then by the classification of extended dual graphs (cf. \cite[Theorem 4.15]{KM98}), $(X,S)$ is plt near $x$. Thus we may assume that $x$ is not a cyclic quotient singularity of $X$. By Lemma \ref{lem: de type mld}, $x$ is a $D$-type singularity and $a(E,X,S)=0$, where $E$ is the prime divisor corresponding to the unique fork of the minimal resolution of $X$ near $x$. Thus $a(E,X,0)\leq 1$, and
$$a(E,X,B)=b_1a(E,X,S)+(1-b_1)a(E,X,0)\leq 1-b_1<\frac{1}{7},$$
which contradicts (1).

(3) The case when $\rho(X)=1$ follows from (1) and \cite[5.1.2]{Sho00}. The general case essentially follows from the proof of \cite[5.1.2]{Sho00} but we give a full proof here for the reader's convenience.  By Lemma \ref{lem: reduced to picard number 1}, there exists a contraction $f: X\rightarrow Y$ and a klt surface pair $(Y,B_Y=bS_Y)$, such that $K_X+B=f^*(K_Y+B_Y)$, $S_Y:=f_*S$, and $\rho(Y)=1$. By (1), $\delta(Y,B_Y)=1$. By Zariski's main theorem, if $S_Y$ is non-singular rational curve (resp. elliptic curve), then $S$ is non-singular (resp. elliptic curve). By (1) and \cite[5.1.2, 5.3 Corollary, 5.4 Proposition]{Sho00}, we may assume that $S_Y$ is a rational curve with exactly $1$ singularity $y$, such that $y$ is a nodal singularity of $S_Y$ and a non-singular point of $Y$. Since $S_Y$ is semi-normal near $y$, $S$ is semi-normal near $f^{-1}(y)$. Thus $S$ is either non-singular, or has at exactly $1$ nodal singularity $x$ near $f^{-1}(y)$. 

Suppose that $x$ is a singular point of $S$, then we only need to show that $x$ is a non-singular point of $X$. By Lemma \ref{lem: (X,T+S) lc}, $(X,S)$ is lc near $x$. Since $x$ is a nodal singularity of $S$, $x$ is a cyclic quotient singularity or a non-singular point. Let $E$ be a divisor contained in the  dual graph of the minimal resolution of $X\ni x$, then $a(E,X,0)\leq 1$ and $a(E,X,S)=0$, so if $x$ is a cyclic quotient singularity, then
$$a(E,X,B)=b_1a(E,X,S)+(1-b_1)a(E,X,0)\leq 1-b_1<\frac{1}{7},$$
which contradicts (1). Therefore, $x$ is a non-singular point.

(4) By Zariski's main theorem, the image of a singular curve under divisorial contraction cannot be non-singular. By Lemma \ref{lem: reduced to picard number 1}, we may assume that $\rho(X)=1$. By (2)(3) and the classification result \cite[Table]{Abe97} for the case when $S$ is a non-singular elliptic curve or a rational curve with at most $1$ nodal singularity, we know that $S$ is a non-singular rational curve. In fact, the case in \cite[Table]{Abe97} with maximal coefficient is \cite[Table, Case 10]{Abe97}, where the coefficient is exactly $\frac{10}{11}$. By (2), $(X,S)$ is plt.
\end{proof}

\subsection{Reduction to Picard number one} We introduce two new set-ups in this subsection. They will be used later in this paper.

Now we get a new set-up.
\begin{setup}\label{setup: klt set-up}
$X,B,b,S$ are as follows:
\begin{enumerate}
\item $(X,B=bS)$ is a projective klt surface pair that is $\frac{1}{7}$-lt at any closed point $x\in X$.
\item $S$ is a non-singular rational curve.
\item $K_X+B\sim_{\mathbb Q}0$ and $b\in (0,1)\cap\mathbb Q$.
\item $(X,S)$ is plt. In particular, $S$ only contains cyclic quotient singularities of $X$.
\end{enumerate}
\end{setup}

\begin{lem}\label{lem: S contain at least 3 sing}
Conditions as in Set-up \ref{setup: klt set-up}. Assume that $S^2>0$. Then:
\begin{enumerate}
    \item $S$ contains at least $3$ singularities of $X$. 
    \item If $S$ contains exactly $3$ singularities of $X$ with indices $p_1,p_2,p_3$, then $\sum_{i=1}^3\frac{1}{p_i}<1$.
    \item If $S$ contains exactly $4$ singularities of $X$ with indices $p_1,p_2,p_3,p_4$, then $\sum_{i=1}^4\frac{1}{p_i}<2$.
\end{enumerate}
\end{lem}
\begin{proof}
Let $x_1,\dots,x_n$ be the singular points of $X$ that are in $S$, $f: X'\rightarrow X$ a resolution of $X$, and $S'$ the strict transform of $S$ on $X'$. We let $p_i$ be the order of $x_i$ for each $i$. Since $x_i$ are cyclic quotient singularities and $S$ is a non-singular rational curve, we have
$$0=(K_X+bS)\cdot S<(K_X+S)\cdot S=(K_{X'}+S')\cdot S'+\sum \left(1-\frac{1}{p_i}\right)=-2+\sum \left(1-\frac{1}{p_i}\right).$$
The lemma immediately follows.
\end{proof}

\begin{thm}[{\cite[Theorem 1.2]{Bel08}, \cite[Theorem 1.1]{Bel09}, \cite[Theorem 1.2]{LX21b}}]\label{thm: picard 1 4 sing}
Let $(X,B)$ be a klt surface pair such that $B>0$ and $K_X+B\equiv 0$. Then the number of singular points on $X$ is $\le2\rho(X)+2$. In particular, if $\rho(X)=1$, then  the number of singular points on $X$ is $\leq 4$.
\end{thm}

\begin{lem}\label{lem: reduce to picard number 1}
Conditions as in Set-up \ref{setup: klt set-up}. Assume that $b>\frac{10}{11}$. Then there exists a contraction $f: X\rightarrow X'$ and a pair $(X',B'=bS')$, such that 
\begin{enumerate}
    \item $X',B',b,S'$ satisfy the conditions of Set-up \ref{setup: klt set-up},
    \item $B'=f_*B$ and $S'=f_*S$,
    \item $\rho(X')=1$, and
    \item $X'$ contains at most $4$ singular points, and $S'$ contains at least $3$ singular points of $X'$.
\end{enumerate}
\end{lem}
\begin{proof}
We let $f: X\rightarrow X'$ be a $K_X$-MMP which terminates with a Mori fiber space $X'\rightarrow Z$, $B':=f_*B$, and $S':=f_*S$.
It is clear that $B'=bS'$. (2) holds by our construction.

We prove (1). Since $K_X+B\sim_{\mathbb Q}0$ and $f$ is $B$-positive, $B$ is not contracted by the MMP and $K_X+B=f^*(K_{X'}+B')$. Thus $X',B',b,S'$ satisfy the conditions of Set-up \ref{setup: klt set-up}(1)(3). By Lemma \ref{lem: reduce to smooth rational curve}(2)(3), $X',B',b,S'$ also satisfy the conditions of Set-up \ref{setup: klt set-up}(2)(4).

By Lemma \ref{lem: mfs 2/3}, $Z$ is a point and $\rho(X')=1$. This implies (3). (4) follows from (3), Lemma \ref{lem: S contain at least 3 sing}, and Theorem \ref{thm: picard 1 4 sing}.
\end{proof}

\begin{setup}\label{setup: rho=1 set-up}
$X,B,b,S,x_1,x_2,x_3,x_4,p_1,p_2,p_3,p_4$ are as follows:
\begin{enumerate}
\item $X,B,b,S$ satisfy the conditions as in Set-up \ref{setup: klt set-up}.
\item $\rho(X)=1$.
\item $x_1,x_2,x_3$ are three cyclic quotient singularities on $X$ that are in $S$, and $p_1,p_2,p_3$ are the indices of $x_1,x_2,x_3$ respectively, such that $2\leq p_1\leq p_2\leq p_3$.
\item $x_4$ and $p_4$ satisfy the following:
\begin{enumerate}
    \item If $X$ has exactly $3$ singular points, then we do not define $x_4,p_4$.
    \item If $X$ has $4$ singular points, then $x_4\not=x_1,x_2,x_3$ is a singular point on $X$ and $p_4\geq 2$ is the order of $x_4$. Moreover, if $x_4$ is contained in $S$, then $p_3\leq p_4$.
\end{enumerate}
\end{enumerate}
\end{setup}

\begin{rem}
We can actually further assume that $x_4$ is not Du Val. Nevertheless, since we do not need this in the proof of the main theorems, we do not make such additional reduction.
\end{rem}

\begin{thm}\label{thm: euler bound}
Conditions as in Set-up \ref{setup: rho=1 set-up}. If $x_4$ exists, then
$\sum_{i=1}^4\frac{1}{p_i}\geq 1.$ 
\end{thm}
\begin{proof}
It follows from \cite[9.2 Corollary]{KM99}.
\end{proof}

\subsection{Two useful lemmas}

We need the following two useful lemmas in this paper. The first lemma is well-known (which should be known by the Italian schools in the 19th century). In some literature, it is called Noether's formula.

\begin{lem}\label{lem: rho+kx2=10}
Let $X$ be a non-singular rational surface. Then $\rho(X)+K_X^2=10$.
\end{lem}
\begin{proof}
We run a $K_X$-MMP 
$$X:=X_0\xrightarrow{f_1} X_1\xrightarrow{f_2}\dots\xrightarrow{f_n} X_n$$ which terminates with a Mori fiber space $X_n\rightarrow Z$. Then $X_n$ is either $\mathbb P^2$ or a Hirzebruch surface, so $\rho(X_n)+K_{X_n}^2=10$. For any $i$, $f_i$ is a blow-up of a point, so $K_{X_{i-1}}^2=K_{X_i}^2-1$ and $\rho(X_{i-1})=\rho(X_i)+1$ for any $1\leq i\leq n$. Thus $\rho(X_i)+K_{X_i}^2$ is a constant for each $i$, and the lemma follows.
\end{proof}

\begin{lem}\label{lem: intersection numbers}
Let $f: X\rightarrow Y$ be a contraction between projective surfaces such that the exceptional locus of $f$ is a prime divisor $E$. Let $C,D$ be two $\Rr$-Cartier $\Rr$-divisors on $X$ such that $f_*C,f_*D$ are $\Rr$-Cartier. Then
$$f_*C\cdot f_*D=C\cdot D-\frac{(C\cdot E)(D\cdot E)}{E^2}.$$
\end{lem}
\begin{proof}
We may assume that $f^*f_*C=C+cE$ and $f^*f_*D=D+dE$ for some real numbers $c,d$. Then $C\cdot E=-cE^2$ and $D\cdot E=-dE^2$, so
$$f_*C\cdot f_*D=f^*f_*C\cdot f^*f_*D=(C+cE)\cdot (D+dE)=C\cdot D-cdE^2=C\cdot D-\frac{(C\cdot E)(D\cdot E)}{E^2}.$$
\end{proof}

\section{Singularity control}\label{sec: control singularities}

\subsection{Coefficients control on the minimal surface}

\begin{lem}\label{lem: check s_y^2}
Let $(X,B=\sum_{i=1}^nb_iB_i)$ be a projective klt surface pair such that $X$ is non-singular, $B_i$ are the prime components of $\Supp B$, and $K_X+B\equiv 0$. Then there exists a birational contraction $f: X\rightarrow X'$ to a minimal surface (in the Italian sense) $X'$. For any such $f: X\rightarrow X'$, we let $B':=f_*B$ and $B_i':=f_*B_i$ be the images of $B$ and $B_i$ on $X'$ respectively. We let $X'\rightarrow Z$ be a $\mathbb P^1$-fibration if $\rho(X')=2$, and let $Z:=\{pt\}$ if $X'=\mathbb P^2$.
\begin{enumerate}
    \item $(X',B')$ is klt and $K_{X'}+B'\equiv 0$.
    \item For any $i$, if $B_i^2\geq 0$, then $B_i'\not=0$.
    \item If $Z=\{pt\}$, then $\sum_{i=1}^n (\deg B_i')b_i=3$, and for any $i,j$,
    \begin{enumerate}
        \item     $(\deg B_i')^2\geq B_i^2$, and
        \item  $(\deg B_i')(\deg B_j')\geq B_i\cdot B_j$ if $B_i'\not=0$ and $B_j'\not=0$.
    \end{enumerate}
    \item If $Z\not=\{pt\}$, then there exist non-negative integers $c_i$, such that
    $$\sum_{B_i'\not=0} c_ib_i=2,\text{ and }c_i>0\text{ if }B_i^2>0.$$
    \item For any $i$, if $B_i^2\geq 0$ and the equation
        $$\sum_{j=1}^n c_jb_j=2, c_i\in\mathbb N^+,c_j\in\mathbb N\text{ for any }j$$
        does not have a solution, then for some solution $(c_1,\dots,c_n)$ of the equation
        $$\sum_{j=1}^n c_jb_j=3, c_i\in\mathbb N^+,c_j\in\mathbb N\text{ for any }j,$$
        we have $c_jc_k\geq B_j\cdot B_k$ for any $j,k$ such that $c_j\not=0$ and $c_k\not=0$. In particular, $c_i^2\geq B_i^2$ for any $i$.
\end{enumerate}
\end{lem}
\begin{proof}
By contracting $(-1)$-curves we get $X\rightarrow X'$. Then either there exists a $\mathbb P^1$-fibration $X'\rightarrow Z$, or $X'=\mathbb P^2$.

Now we prove (1-5) for any such $f: X\rightarrow X'$ and $X'\rightarrow Z$.

Since $K_X+B\equiv 0$, $K_X+B=f^*(K_{X'}+B')$, and (1) follows.

(2) Since $X'$ is smooth, $f$ is a composition of ordinary blow-ups of $(-1)$-curves. By Lemma \ref{lem: intersection numbers}, $B_i'^2\geq B_i^2$, and $B_i^2=B_i'^2$ if and only if the MMP is an isomorphism near $B_i$. Thus if the MMP contracts $B_i$, then $0=B_i'^2>B_i^2\geq 0$, a contradiction.

(3) Since $K_{X'}+B'\equiv 0$, $\sum_{B_i'\not=0} (\deg B_i')b_i=3$. By Lemma \ref{lem: intersection numbers}, $(\deg B_i')^2=B_i'^2\geq B_i^2$ for any $i$, which is (3.a). By Lemma \ref{lem: intersection numbers} again, if $B_i,B_j$ are not contracted by $f$, then $(\deg B_i')(\deg B_j')=B_i'\cdot B_j'\geq B_i\cdot B_j$, which is (3.b).

(4) Let $F$ be a general fiber of $X'\rightarrow Z$, $B_F:=B'|_F$, and $B_{i,F}:=B_i'|_F$ for each $i$. Since $-K_{X'}$ is ample over $Z$, $-K_F$ is ample, so $F\cong\mathbb P^1$. Since $K_{X'}+B'\equiv 0$, $K_F+B_F\equiv 0$. Thus $\deg B_F=2$, and may let $c_i:=\deg B_{i,F}$ for each $i$.

(5) By (2), $B_i'\not=0$. By (4), $Z=\{pt\}$. By (3), $X'\cong\mathbb P^2$, and for some $(c_1,\dots,c_n)$ which satisfies the equation
$$\sum_{j=1}^n c_jb_j=3, c_i\in\mathbb N^+,c_j\in\mathbb N\text{ for any }j,$$
$\deg B_j'=c_j$ for any $j$. Now (5) follows from (3.a) and (3.b).
\end{proof}

\subsection{Invariants of singularities}

\begin{nota}\label{nota: gamma invariant}
Let $X\ni x$ (resp. $X$) be a klt surface germ (resp. a klt surface), $f: Y\rightarrow X$ the minimal resolution of $X\ni x$ (resp. $X$), and
$$K_Y+\sum_{i=1}^m b_iE_i=f^*K_X,$$
where $E_1,\dots,E_m$ are the prime $f$-exceptional divisors and $b_i=1-a(E_i,X,0)$ for each $i$. We put 
$$\gamma(X\ni x)\text{(resp. }\gamma(X)\text{)}:=m-\sum_{i=1}^mb_i(K_Y\cdot E_i).$$
It is clear that
$$\gamma(X)=\sum_{x\in X}\gamma(X\ni x).$$
\end{nota}
The following lemma provides an important characterization of $\gamma(X)$.
\begin{lem}\label{lem: global characterization gamma(X)}
Let $X$ be a projective klt surface  and $f: Y\rightarrow X$ the minimal resolution of $X$. Then 
$$\gamma(X)=\rho(Y/X)+K_Y^2-K_X^2.$$
\end{lem}
\begin{proof}
Indeed $K_Y+\sum_{i=1}^m b_iE_i=f^*K_X$ as in Notation \ref{nota: gamma invariant}. Hence
\begin{align*}
    \gamma(X)&=m-\sum_{i=1}^mb_i(K_Y\cdot E_i)=\rho(Y/X)-\left(K_Y+\sum_{i=1}^mb_iE_i\right)\cdot K_Y+K_Y^2\\
    &=\rho(Y/X)-f^*K_X\cdot K_Y+K_Y^2=\rho(Y/X)-f^*K_X\cdot \left(f^*K_X-\sum_{i=1}^mb_iE_i\right)+K_Y^2\\
    &=\rho(Y/X)-f^*K_X\cdot f^*K_X+K_Y^2=\rho(Y/X)+K_Y^2-K_X^2.
\end{align*}
\end{proof}

\begin{lem}\label{lem: gamma(X) equality}
Conditions as in Set-up \ref{setup: rho=1 set-up}. Then:
\begin{enumerate}
    \item
   $$\gamma(X)=9-\frac{b^2}{1-b}\left(-2+\sum_{x_i\in S}\left(1-\frac{1}{p_i}\right)\right).$$
\item $$\sqrt{(9-\gamma(X))^2+4(9-\gamma(X))\left(-2+\sum_{x_i\in S}\left(1-\frac{1}{p_i}\right)\right)}$$
is a rational number.
\end{enumerate}
\end{lem}
\begin{proof}
Let $f: Y\rightarrow X$ be the minimal resolution of $X$. Then $Y$ is a non-singular rational surface. By Lemma \ref{lem: rho+kx2=10}, $10=\rho(Y)+K_Y^2$. By Lemma \ref{lem: global characterization gamma(X)},
$$1=\rho(X)=\rho(Y)-\rho(Y/X)=10-K_Y^2-\rho(Y/X)=10-\gamma(X)-K_X^2.$$
Since 
$$(K_X+S)\cdot S=-2+\sum_{x_i\in S}\left(1-\frac{1}{p_i}\right)$$
and $K_X+bS\equiv 0$, we have
$$K_X^2=\frac{b^2}{1-b}\left(-2+\sum_{x_i\in S}\left(1-\frac{1}{p_i}\right)\right),$$
and (1) follows. (2) follows from (1) and the fact that $b$ is a rational number.
\end{proof}

\begin{notalem}\label{notalem: q(x,s)}
Let $(X,S)$ be a surface pair and $x$ a closed point in $X$, such that $(X,S)$ is plt near $x$. Then $x$ is a cyclic quotient singularity of $X$ and $S$ is a non-singular curve. We define $q(X\ni x,S)$ in the following way.

We let $f: Y\rightarrow X$ be the minimal resolution of $X$ and $E_1,\dots,E_k$ the prime $f$-exceptional divisors.  It is clear that the dual graph of $X\ni x$ is a chain, and by \cite[Theorem 4.15]{KM98}, we may suppose that exactly one tail of the chain intersects $S_Y:=f^{-1}_*S$. Possibly reordering indices, we may assume that the dual graph of $f$ and $S_Y$ is
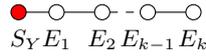
\begin{center}
\begin{tikzpicture}
         \draw[fill=red] (1.0,0) circle (0.1);   
          \node [below] at (1.1,-0.1) {\footnotesize$S_Y$};
     \draw (1.1,0)--(1.4,0);
         \draw (1.5,0) circle (0.1);
         \node [below] at (1.5,-0.1) {\footnotesize$E_1$};
         \draw (1.6,0)--(2,0);
         \draw (2.1,0) circle (0.1);
         \node [below] at (2.1,-0.1) {\footnotesize$E_2$};
         \draw [dashed] (2.2,0)--(2.6,0);
         \draw (2.7,0) circle (0.1);
         \node [below] at (2.7,-0.1) {\footnotesize$E_{k-1}$};
         \draw (2.8,0)--(3.2,0);
         \draw (3.3,0) circle (0.1);
         \node [below] at (3.3,-0.1) {\footnotesize$E_k$};



    \end{tikzpicture}
\end{center}
We define
$$q(X\ni x,S):=\det(\mathcal{D}(\cup_{i=2}^kE_i))$$
if $k\geq 2$, and $q(X\ni x,S):=1$ if $k=1$.

Notice that $q(X\ni x,S)=p\mult_{E_1}S$, where $p$ is the order of $X\ni x$. To see this, we let $g: Y'\rightarrow X$ be the extraction of $E_1$ and let $E_1'$ be the image of $E_1$ on $Y'$, then $E_1'^2=-\frac{p}{q(X\ni x,S)}$. Thus
$$0=S\cdot g_*E_1'=g^*S\cdot E_1'=(g^{-1}_*S+(\mult_{E_1}S)E_1')\cdot E_1'=1-\frac{p}{q(X\ni x,S)}\cdot\mult_{E_1}S,$$
and $q(X\ni x,S)=p\mult_{E_1}S$.
\end{notalem}

\begin{lem}\label{lem: control sy^2}
Conditions as in Set-up \ref{setup: rho=1 set-up}. Then
$$\frac{1}{1-b}\left(-2+\sum_{x_i\in S}\left(1-\frac{1}{p_i}\right)\right)-\sum_{x_i\in S}\frac{q(X\ni x_i,S)}{p_i}\in\mathbb Z.$$
\end{lem}
\begin{proof}
Let $f: Y\rightarrow X$ be the minimal resolution of $X$ and $S_Y:=f^{-1}_*S$. Then $S_Y^2$ is an integer, and
$$S_Y^2=S^2-\sum_{x_i\in S}\frac{q(X\ni x_i,S)}{p_i}.$$
Since
$$(K_X+S)\cdot S=-2+\sum_{x_i\in S}\left(1-\frac{1}{p_i}\right)$$
and $K_X+bS\equiv 0$, we have
$$S^2=\frac{1}{1-b}\left(-2+\sum_{x_i\in S}\left(1-\frac{1}{p_i}\right)\right),$$
hence
$$\frac{1}{1-b}\left(-2+\sum_{x_i\in S}\left(1-\frac{1}{p_i}\right)\right)-\sum_{x_i\in S}\frac{q(X\ni x_i,S)}{p_i}=S_Y^2\in\mathbb Z.$$
\end{proof}

\section{Classification of singularities on the boundary}\label{sec: sings on S}

In this section, we classify all possible singularities on $S$.

\begin{lem}\label{lem: classification 10/11}
Let $b\in (\frac{6}{7},1)$ be a real number. Let $(X\ni x,S)$ be a plt surface germ where $X\ni x$ is a cyclic quotient singularity and $S$ a prime divisor on $X$, such that $(X\ni x,B=bS)$ is $\frac{1}{7}$-lt for some $b$. Suppose that the dual graph of $X\ni x$ and $S$ is
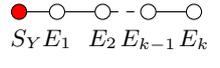
\begin{center}
\begin{tikzpicture}
        \draw[fill=red] (1.0,0) circle (0.1);   
          \node [below] at (1.1,-0.1) {\footnotesize$S_Y$};
     \draw (1.1,0)--(1.4,0);
         \draw (1.5,0) circle (0.1);
         \node [below] at (1.5,-0.1) {\footnotesize$E_1$};
         \draw (1.6,0)--(2,0);
         \draw (2.1,0) circle (0.1);
         \node [below] at (2.1,-0.1) {\footnotesize$E_2$};
         \draw [dashed] (2.2,0)--(2.6,0);
         \draw (2.7,0) circle (0.1);
         \node [below] at (2.7,-0.1) {\footnotesize$E_{k-1}$};
         \draw (2.8,0)--(3.2,0);
         \draw (3.3,0) circle (0.1);
         \node [below] at (3.3,-0.1) {\footnotesize$E_k$};



    \end{tikzpicture},
\end{center}
where $E_1,\dots,E_k$ are the exceptional prime divisors of the minimal resolution of $X\ni x$ and $S_Y$ is the strict transform of $S$ on the minimal resolution of $X\ni x$. Then the order of $X\ni x$ is $\leq\frac{7b}{7b-6}$. In particular, there are only finitely many types of $(X\ni x,S)$, i.e. the possibilities of $\mathcal{D}(S_Y\cup_{i=1}^kE_i)$ are finite.

 Moreover, if $b\geq\frac{10}{11}$, then the order of $X\ni x$ is $\leq 17$, and $(X\ni x,S)$ is of one of the types listed in Table \ref{Table: singularity on S} below. Here
    \begin{enumerate}
    \item ``Singularity" stands for the type of the cyclic quotient singularity $X\ni x$,
        \item ``Weights" stands for $(-E_1^2,\dots,-E_k^2)$,
        \item ``$X$ l.d." stands for $(a(E_1,X),\dots,a(E_k,X))$,
        \item ``Pair l.d." stands for $(a(E_1,X,B),\dots,a(E_k,X,B))$, 
        \item ``$\gamma$" is the value $\gamma(X\ni x)$", the constant as in Notation \ref{nota: gamma invariant}, 
        \item ``$q$" stands for $q(X\ni x,S)$, and
        \item ``$\sup b$" stands for the supremum of all possible values of $b$, which is a value in $[\frac{10}{11},1)$.
    \end{enumerate}
    In Table \ref{Table: singularity on S}, we denote $m^n$ the sequence $(m)_{i=1}^n$.
\begin{center}
    \begin{longtable}{|l|l|l|l|l|l|l|l|}
 			\caption{Singularities of $\frac{1}{7}$-lt $(X\ni x,bS)$ with $b\geq\frac{10}{11}$}\label{Table: singularity on S}\\
		\hline
		No.&Singularity& Weights & $X$ l.d. & Pair l.d. & $\gamma$ & $q$ & $\sup b$\\ 
		\hline
1& $\frac{1}{2}(1,1)$ & $(2)$ & $(1)$ & $(\frac{2-b}{2})$ & $1$ & $1$ & $1$ \\ 
\hline
2 & $\frac{1}{3}(1,1)$ & $(3)$ & $(\frac{2}{3})$ & $(\frac{2-b}{3})$ & $\frac{2}{3}$ & $1$ & $1$ \\
\hline
3 & $\frac{1}{4}(1,1)$ & $(4)$ & $(\frac{1}{2})$ & $(\frac{2-b}{4})$ &  $0$ & $1$ & $1$ \\
\hline
4 & $\frac{1}{5}(1,1)$ & $(5)$ & $(\frac{2}{5})$ & $(\frac{2-b}{5})$ & $-\frac{4}{5}$ & $1$ & $1$  \\
\hline
5 & $\frac{1}{6}(1,1)$ & $(6)$ & $(\frac{1}{3})$ & $(\frac{2-b}{6})$ & $-\frac{5}{3}$ & $1$ & $1$  \\
\hline
6 & $\frac{1}{7}(1,1)$ & $(7)$ & $(\frac{2}{7})$ & $(\frac{2-b}{7})$ & $-\frac{18}{7}$ & $1$ & $1$ \\
\hline
7 & $\frac{1}{7}(1,2)$ & $(4,2)$ & $(\frac{3}{7},\frac{5}{7})$ & $(\frac{3-2b}{7},\frac{5-b}{7})$ & $\frac{6}{7}$ & $2$ & $1$ \\
\hline
8 & $\frac{1}{8}(1,3)$ & $(3,3)$ & $(\frac{1}{2},\frac{1}{2})$ & $(\frac{4-3b}{8},\frac{4-b}{8})$ & $1$ & $3$ & $\frac{20}{21}$  \\
\hline
9 & $\frac{1}{5}(1,2)$ & $(3,2)$ & $(\frac{3}{5},\frac{4}{5})$ & $(\frac{3-2b}{5},\frac{4-b}{5})$ & $\frac{8}{5}$ & $2$ & $1$  \\
\hline
10 & $\frac{1}{7}(1,3)$ & $(3,2,2)$ & $(\frac{4}{7},\frac{5}{7},\frac{6}{7})$ & $(\frac{4-3b}{7},\frac{5-2b}{7},\frac{6-b}{7})$ & $\frac{18}{7}$ & $3$ & $1$ \\
\hline
11 & $\frac{1}{9}(1,4)$ & $(3,2,2,2)$ & $(\frac{5}{9},\frac{2}{3},\frac{7}{9},\frac{8}{9})$ & $(\frac{5-4b}{9},\frac{2-b}{3},\frac{7-2b}{9},\frac{8-b}{9})$ & $\frac{32}{9}$ & $4$ & $\frac{13}{14}$ \\
\hline
12 & $\frac{1}{3}(1,2)$ & $(2,2)$ & $(1,1)$ & $(\frac{3-2b}{3},\frac{3-b}{3})$ & $2$ & $2$ & $1$ \\
\hline
13 & $\frac{1}{5}(1,2)$ & $(2,3)$ & $(\frac{4}{5},\frac{3}{5})$ & $(\frac{4-3b}{5},\frac{3-b}{5})$ & $\frac{8}{5}$ & $3$ & $1$ \\
\hline
14 & $\frac{1}{7}(1,2)$ & $(2,4)$ & $(\frac{5}{7},\frac{3}{7})$ & $(\frac{5-4b}{7},\frac{3-b}{7})$ & $\frac{6}{7}$ & $4$ & $1$  \\
\hline
15 & $\frac{1}{9}(1,2)$ & $(2,5)$ & $(\frac{2}{3},\frac{1}{3})$ & $(\frac{6-5b}{9},\frac{3-b}{9})$ & $0$ & $5$ & $\frac{33}{35}$ \\
\hline
16 & $\frac{1}{4}(1,3)$ & $(2,2,2)$ & $(1,1,1)$ & $(\frac{4-3b}{4},\frac{4-2b}{4},\frac{4-b}{4})$& $3$ & $3$ & $1$ \\
\hline
17 & $\frac{1}{8}(1,5)$ & $(2,3,2)$ & $(\frac{3}{4},\frac{1}{2},\frac{3}{4})$ & $(\frac{6-5b}{8},\frac{2-b}{4},\frac{6-b}{8})$ & $\frac{5}{2}$ & $5$ & $\frac{34}{35}$ \\
\hline
18 & $\frac{1}{5}(1,4)$ & $(2,2,2,2)$ & $(1^4)$ & $(\frac{5-(5-i)b}{5})_{i=1}^4$ & $4$ & $4$ & $1$  \\
\hline
19 & $\frac{1}{7}(1,3)$ & $(2,2,3)$ & $(\frac{6}{7},\frac{5}{7},\frac{4}{7})$ & $(\frac{6-5b}{7},\frac{5-3b}{7},\frac{4-b}{7})$ & $\frac{18}{7}$ & $5$ & $1$ \\
\hline
20 & $\frac{1}{11}(1,7)$ & $(2,3,2,2)$ & $(\frac{8}{11},\frac{5}{11},\frac{7}{11},\frac{9}{11})$ & $(\frac{8-7b}{11},\frac{5-3b}{11},\frac{7-2b}{11},\frac{9-b}{11})$& $\frac{38}{11}$ & $7$ & $\frac{45}{49}$ \\
\hline
21 & $\frac{1}{6}(1,5)$ & $(2^5)$ & $(1^5)$ & $(\frac{6-(6-i)b}{6})_{i=1}^5$ & $5$ & $5$ & $1$  \\
\hline
22 & $\frac{1}{10}(1,3)$ & $(2,2,4)$ & $(\frac{4}{5},\frac{3}{5},\frac{2}{5})$ & $(\frac{8-7b}{10},\frac{3-2b}{5},\frac{4-b}{10})$& $\frac{9}{5}$ & $7$ & $\frac{46}{49}$  \\
\hline
23 & $\frac{1}{7}(1,6)$ & $(2^6)$ & $(1^6)$ & $(\frac{7-(7-i)b}{7})_{i=1}^6$ & $6$ & $6$ & $1$ \\
\hline
24 & $\frac{1}{9}(1,4)$ & $(2,2,2,3)$ & $(\frac{8}{9},\frac{7}{9},\frac{2}{3},\frac{5}{9})$ & $(\frac{8-7b}{9},\frac{7-5b}{9},\frac{2-b}{3},\frac{5-b}{9})$& $\frac{32}{9}$ & $7$ & $\frac{47}{49}$ \\
\hline
25 & $\frac{1}{11}(1,7)$ & $(2,2,3,2)$ & $(\frac{9}{11},\frac{7}{11},\frac{5}{11},\frac{8}{11})$ & $(\frac{9-8b}{11},\frac{7-5b}{11},\frac{5-2b}{11},\frac{8-b}{11})$ & $\frac{38}{11}$ & $8$ & $\frac{13}{14}$ \\
\hline
26 & $\frac{1}{8}(1,7)$ & $(2^7)$ & $(1^7)$ & $(\frac{8-(8-i)b}{8})_{i=1}^7$ & $7$ & $7$ & $\frac{48}{49}$ \\
\hline
27 & $\frac{1}{9}(1,8)$ & $(2^8)$ & $(1^8)$ & $(\frac{9-(9-i)b}{9})_{i=1}^8$ & $8$ & $8$ & $\frac{27}{28}$ \\
\hline
28 & $\frac{1}{11}(1,5)$ & $(2,2,2,2,3)$ & $(\frac{11-i}{11})_{i=1}^5$& $(\frac{11(1-b)-i(1-2b)}{11})_{i=1}^5$& $\frac{50}{11}$ & $9$ & $\frac{59}{63}$  \\
\hline
29 & $\frac{1}{13}(1,4)$ & $(2,2,2,4)$ & $(\frac{11}{13},\frac{9}{13},\frac{7}{13},\frac{5}{13})$ & $(\frac{11-10b}{13},\frac{9-7b}{13},\frac{7-4b}{13},\frac{5-b}{13})$ & $\frac{36}{13}$ & $10$ & $\frac{32}{35}$ \\
\hline
30 & $\frac{1}{10}(1,9)$ & $(2^9)$ & $(1^9)$ & $(\frac{10-(10-i)b}{10})_{i=1}^9$ & $9$ & $9$ & $\frac{20}{21}$ \\
\hline
31 & $\frac{1}{11}(1,10)$ & $(2^{10})$ & $(1^{10})$ & $(\frac{11-(11-i)b}{11})_{i=1}^{10}$ & $10$ & $10$ & $\frac{33}{35}$ \\
\hline
32 & $\frac{1}{13}(1,6)$ & $(2^5,3)$ & $(\frac{13-i}{13})_{i=1}^6$& $(\frac{13(1-b)-i(1-2b)}{13})_{i=1}^6$& $\frac{72}{13}$ & $11$ & $\frac{71}{77}$ \\
\hline
33 & $\frac{1}{12}(1,11)$ & $(2^{11})$ & $(1^{11})$ & $(\frac{12-(12-i)b}{12})_{i=1}^{11}$ & $11$ & $10$ & $\frac{72}{77}$ \\
\hline
34 & $\frac{1}{13}(1,12)$  & $(2^{12})$ & $(1^{12})$ & $(\frac{13-(13-i)b}{13})_{i=1}^{12}$ & $12$ & $11$ & $\frac{13}{14}$ \\
\hline
35 & $\frac{1}{15}(1,7)$ & $(2^6,3)$ & $(\frac{15-i}{15})_{i=1}^7$& $(\frac{15(1-b)-i(1-2b)}{15})_{i=1}^7$& $\frac{98}{15}$ & $13$ & $\frac{83}{91}$ \\
\hline
36 & $\frac{1}{14}(1,13)$ &  $(2^{13})$ & $(1^{13})$ & $(\frac{14-(14-i)b}{14})_{i=1}^{13}$ & $13$  & $13$ & $\frac{12}{13}$ \\
\hline
37 & $\frac{1}{15}(1,14)$ &  $(2^{14})$ & $(1^{14})$ & $(\frac{15-(15-i)b}{15})_{i=1}^{14}$ & $14$ & $14$ & $\frac{45}{49}$ \\
\hline
38 & $\frac{1}{16}(1,15)$ &  $(2^{15})$ & $(1^{15})$ & $(\frac{16-(16-i)b}{16})_{i=1}^{15}$ & $15$ & $15$ & $\frac{32}{35}$ \\
\hline
39 & $\frac{1}{17}(1,16)$ &  $(2^{16})$ & $(1^{16})$ & $(\frac{17-(17-i)b}{17})_{i=1}^{16}$ & $16$ & $16$ & $\frac{51}{56}$  \\
\hline
\end{longtable}
\end{center}
\end{lem}
\begin{proof}
We let 
\begin{itemize}
\item $w_1:=-E_1^2$,
\item $m:=\det(-(E_i\cdot E_j))_{i,j\geq 1}$,
\item $q:=q(X\ni x,S)=\det(-(E_i\cdot E_j))_{i,j\geq 2}$, and
\item $n:=\det(-(E_i\cdot E_j))_{i,j\geq 3}$ if $k\geq 3$, $n:=1$ if $k=2$, and $n:=0$ if $k=1$.
\end{itemize}
Then $m$ is the order of $X\ni x$, and $m=qw_1-n$. By \cite[3.1.8 Lemma, 3.1.10 Lemma]{Kol+92} (see also \cite[Page 13]{Ale93}, \cite[Lemma A.1]{CH21}), 
$$\frac{1}{7}<a(E_1,X,B)=\frac{(1-b)q+1}{m}=\frac{(1-b)q+1}{qw_1-n}.$$
This implies that
$$m<7(1-b)q+7\text{ and }q(w_1-7(1-b))<n+7.$$
Since $w_1\geq 2$ and $n\leq q-1$, we have
$$q(-5+7b)\leq q(w_1-7(1-b))<n+7\leq q+6,$$
so 
$$q<\frac{6}{7b-6}.$$
Thus
$$m<7(1-b)q+7<7(\frac{6-6b}{7b-6}+1)=\frac{7b}{7b-6}.$$
Thus there are only finitely many possibilities of $m$. Therefore, the possibilities of $\mathcal{D}(S_Y\cup_{i=1}^kE_i)$ are also finite.

Now we assume that $b\geq\frac{10}{11}$. Then
$$m<\frac{7\cdot\frac{10}{11}}{7\cdot\frac{10}{11}-6}=\frac{35}{2},$$
so $m\leq 17$. (In fact, we only need $b\geq\frac{101}{119}$ to make $m\leq 17$.)

We already can get Table \ref{Table: singularity on S} directly by computer enumeration, but for the reader's convenience, we provide a detailed enumeration here. We have
$$q\left(w_1-\frac{7}{11}\right)<n+7.$$
If $q=1$, then $n=0$, and we have $w_1<7+\frac{7}{11}$. Thus $w_1\in\{2,3,4,5,6,7\}$ and they correspond to No. 1--6 respectively. In the following, we assume that $q\geq 2$. Then since $n\leq q-1$, $w_1-\frac{7}{11}<1+\frac{6}{q}$. Thus $w_1<1+\frac{7}{11}+\frac{6}{q}$. In particular, $w_1\leq 4$. 

If $w_1=4$, then $q=2$, hence $n=0$, which corresponds to No. 7. If $w_1=3$, then $q\in\{2,3,4\}$ and $\frac{26}{11}q<n+7$. Thus $(q,n)\in \{(2,1),(3,1),(3,2),(4,3)\}$. They correspond to No. 8--11 respectively. In the following, we assume that $w_1=2$. Now we have
$$\frac{15}{11}q<n+7.$$
Since $q\geq n+1$, we have $n<\frac{31}{2}$, hence $n\leq 15$.

If $n=1$, then $2\leq q\leq 5$, so $q\in\{2,3,4,5\}$. They correspond to No. 12--15 respectively.

If $n=2$, then $3\leq q\leq 6$, so $q\in\{3,5\}$. They correspond to No. 16--17 respectively.

If $n=3$, then $4\leq q\leq 7$, so $q\in\{4,5,7\}$. They correspond to No. 18--20 respectively.

If $n=4$, then $5\leq q\leq 8$, so $q\in\{5,7\}$. They correspond to No. 21--22 respectively.

If $n=5$, then $6\leq q\leq 8$, so $q\in\{6,7,8\}$. They correspond to No. 23--25 respectively.

If $n=6$, then $7\leq q\leq 9$, so $q=7$, which corresponds to No. 26.

If $n=7$, then $8\leq q\leq 10$, so $q\in\{8,9,10\}$. They correspond to No. 27--29 respectively.

If $n=8$, then $9\leq q\leq 10$, so $q=9$, which corresponds to No. 30.

If $n=9$, then $10\leq q\leq 11$, so $q\in\{10,11\}$, which corresponds to No. 31--32 respectively.

If $n=10$, then $11\leq q\leq 12$, so $q=11$, which corresponds to No. 33.

If $n=11$, then $12\leq q\leq 13$, so $q\in\{12,13\}$, which corresponds to No. 34--35 respectively.

If $n=12$, then $q=13$, which corresponds to No. 36. 

If $n=13$, then $q=14$, which corresponds to No. 37. 

If $n=14$, then $q=15$, which corresponds to No. 38.

If $n=15$, then $q=16$, which corresponds to No. 39.
\end{proof}

\begin{rem}
Lemma \ref{lem: classification 10/11} implies that $\Rr$-complementary thresholds in dimension $2$ do not have any accumulation point that is strictly larger than $\frac{6}{7}$. We give a sketch of the proof here.

By standard arguments (cf. Proof of Theorem \ref{thm: 12/13 rct} below), we only need to show that for any projective klt log Calabi-Yau surface pair $(X,bS)$ where $S$ is a prime divisor and $b>\frac{6}{7}+\epsilon$. there are only finitely many possibilities of $b$. We may easily reduce to the case of Set-up \ref{setup: rho=1 set-up} with $b>\frac{6}{7}+\epsilon$. By Lemma \ref{lem: classification 10/11}, there are only finitely many possibilities of $x_1,x_2,x_3$, hence there are only finitely many possibilities of $p_1,p_2,p_3$. By Theorem \ref{thm: euler bound} there are only finitely many possibilities of $p_4$. By Lemma \ref{lem: gamma(X) equality}, there are only finitely many possibilities of $b$.
\end{rem}

\section{Proof of the gap theorem}\label{sec: proof of the gap theorem}

In this section, by using classification and with the help of computer programs, we prove Theorem \ref{thm: global lct 12/13}. First we introduce a new setting that will be used throughout the section.

\begin{setup}\label{setup: setup with resolution}
$X,B,b,S,x_1,x_2,x_3,x_4,p_1,p_2,p_3,p_4,f,Y,B_Y,S_Y,E_{i,j},b_{i,j}$ are as follows:
\begin{enumerate}
    \item $X,B,b,S,x_1,x_2,x_3,x_4,p_1,p_2,p_3,p_4$ satisfy the conditions as in Set-up \ref{setup: rho=1 set-up}.
    \item $f: Y\rightarrow X$ is the minimal resolution of $X$.
    \item $K_Y+B_Y=f^*(K_X+B)$.
    \item $S_Y$ is the strict transform of $S$ on $Y$.
    \item $E_{i,j}$ are the prime $f$-exceptional divisors, such that
    \begin{enumerate}
        \item $\Center_XE_{i,j}=x_i$ for any $i,j$,
        \item $E_{i,j}$ intersects $E_{i,k}$ if and only if $|j-k|\leq 1$, and
        \item $E_{i,1}$ intersects $S_Y$ for any $1\leq i\leq 3$.
    \end{enumerate}
    \item $b_{i,j}:=1-a(E_{i,j},X,B)$ for any $i,j$.
    \item If $x_4\not\in S$, then $x_4$ is a cyclic quotient singularity. Moreover, suppose that there are $m$ exceptional divisors on $Y$ over $x_4$ for some positive integer $m$, then either $b_{4,j}=b_{4,m+1-j}$ for any $1\leq j\leq m$, or $b_{4,k_0}>b_{4,m+1-k_0}$, where $k_0:=\min\{k\mid b_{4,j}\not=b_{4,m+1-j}\}$.
\end{enumerate}
By Notation-Lemma \ref{notalem: q(x,s)} and Conditions (5) and (7), $f,Y,S_Y,b_{i,j}$ are uniquely determined, and $E_{i,j}$ are uniquely determined when $1\leq i\leq 3$, and $E_{4,j}$ is uniquely determined possibly up to opposite ordering, i.e. $j\leftrightarrow m+1-j$.
\end{setup}

\subsection{Three singular points case}

\begin{lem}\label{lem: summary three singular points condition}
Conditions as in Set-up \ref{setup: rho=1 set-up}. Suppose that $b>\frac{10}{11}$ and $X$ has exactly $3$ singular points. Then \begin{enumerate}
\item $(X\ni x_i,S)$ is in Table \ref{Table: singularity on S} for each $i$, 
\item $\sum_{i=1}^3\frac{1}{p_i}<1$,
\item $\sum_{i=1}^3\gamma(X\ni x_i)=9-\frac{b^2}{1-b}(1-\sum_{i=1}^3\frac{1}{p_i})<9-\frac{100}{11}(1-\sum_{i=1}^3\frac{1}{p_i})$,
\item $\sqrt{(9-\sum_{i=1}^3\gamma(X\ni x_i))^2+4(9-\sum_{i=1}^3\gamma(X\ni x_i))(1-\sum_{i=1}^3\frac{1}{p_i})}\in\mathbb Q$,
\item $\frac{1}{1-b}(1-\sum_{i=1}^3\frac{1}{p_i})-\sum_{i=1}^3\frac{q(X\ni x_i,S)}{p_i}\in\mathbb Z$, and
\item $(X\ni x_i,bS)$ is $\frac{1}{7}$-lt for for any $i$.
\end{enumerate}
\end{lem}
\begin{proof}
(1)(6) follow immediately from our conditions. (2) follows from Lemma \ref{lem: S contain at least 3 sing}. (3) follows from  Lemma \ref{lem: gamma(X) equality}(1) and the assumption $b>\frac{10}{11}$. (4) follows from Lemma \ref{lem: gamma(X) equality}(2). (5) follows from Lemma \ref{lem: control sy^2}.
\end{proof}

We need the following lemma to show the uniqueness of $(X,B)$ when $b=\frac{12}{!3}$.

\begin{lem}\label{lem: uniquenss 12/13}
Conditions as in Set-up \ref{setup: setup with resolution}. Suppose that $b=\frac{12}{13}$, $X$ has exactly $3$ singular points, and
$$(p_1,p_2,p_3,q(X\ni x_1,S),q(X\ni x_2,S),q(X\ni x_3,S))=(3,4,5,2,3,2).$$
Then $(X,B=\frac{12}{13}S)$ is isomorphic to
$$\left(\mathbb P(3,4,5),\frac{12}{13}(x^3y+y^2z+z^2x=0)\right).$$
Moreover, the isomorphism $(X,B)\rightarrow \left(\mathbb P(3,4,5),\frac{12}{13}(x^3y+y^2z+z^2x=0)\right)$ is unique.
\end{lem}
\begin{proof}
We have
$$B_Y=bS_Y+\sum_{i,j}b_{i,j}E_{i,j}=\frac{12}{13}S_Y+\frac{8}{13}E_{1,1}+\frac{4}{13}E_{1,2}+\frac{9}{13}E_{2,1}+\frac{6}{13}E_{2,2}+\frac{3}{13}E_{2,3}+\frac{10}{13}E_{3,1}+\frac{5}{13}E_{3,2},$$
and $S_Y^2=1$. By Set-up \ref{setup: klt set-up}, $S_Y$ is a rational curve. Since $Y$ is a smooth rational surface, $S_Y$ is base-point-free. Let $g: Y\rightarrow\mathbb P^n$ be the morphism defined by $|S_Y|$ and $Z$ the image of $Y$ in $\mathbb P^n$.

We show that $Z=Y\cong\mathbb P^2$. Indeed, let $S_Z:=\mathcal{O}_{\mathbb P^n}(1)|_Z$. Since $S_Y=g^*\mathcal{O}_{\mathbb P^n}(1)$,
$$1=S_Y^2=(\deg(Y\rightarrow Z))\cdot S_Z^2.$$
Hence $S_Z^2=1$. Thus $Y\rightarrow Z$ is birational and $\deg Z=1$. Thus $Z=\mathbb P^2$. By Zariski's main theorem, $g$ is a birational contraction.

We let $F_i$ be the image of $E_{i,1}$ on $Z\cong\mathbb P^2$ for each $i$ and $B_Z:=g_*B_Y$. Since $g$ is defined by $|S_Y|$,
$$K_{\mathbb P^2}+B_Z=K_{\mathbb P^2}+\frac{12}{13}S_Z+\frac{8}{13}F_1+\frac{9}{13}F_2+\frac{10}{13}F_3\equiv 0,$$
and $\deg S_Z=\deg F_1=\deg F_2=\deg F_3=1$. Since $\frac{8}{13}+\frac{9}{13}+\frac{10}{13}>2$ and $(Z,B_Z)$ is klt, no three divisors of $S_Z,F_1,F_2,F_3$ have a common point. 
 Thus $$(Z,B_Z=\frac{12}{13}S_Z+\frac{8}{13}E_{1,1,Z}+\frac{9}{13}E_{2,1,Z}+\frac{10}{13}E_{3,1,Z})\cong(\mathbb P^2,\frac{12}{13}L_0+\frac{8}{13}L_1+\frac{9}{13}L_2+\frac{10}{13}L_3),$$
 where $L_0,L_1,L_2,L_3$ are four lines in general position. Since the coefficients of $L_0,L_1,L_2,L_3$ are different, the automorphism group of $(\mathbb P^2,\frac{12}{13}L_0+\frac{8}{13}L_1+\frac{9}{13}L_2+\frac{10}{13}L_3)$ is trivial. So the isomorphism $(Z,B_Z)\rightarrow (\mathbb P^2,\frac{12}{13}L_0+\frac{8}{13}L_1+\frac{9}{13}L_2+\frac{10}{13}L_3)$ is unique.

 Now by our construction and by computing the log discrepancies, $(Y,B_Y)$ and $(X,B)$ are uniquely determined by $(Z,B_Z)$. In particular, $(X,B)$ is unique up to isomorphism. Since $\left(\mathbb P(3,4,5),\frac{12}{13}(x^3y+y^2z+z^2x=0)\right)$ is an example which satisfies all our conditions, any $(X,B)$ satisfying our conditions is isomorphic to $\left(\mathbb P(3,4,5),\frac{12}{13}(x^3y+y^2z+z^2x=0)\right)$. Moreover, the isomorphism $(X,B)\rightarrow \left(\mathbb P(3,4,5),\frac{12}{13}(x^3y+y^2z+z^2x=0)\right)$ is unique.
\end{proof}

\begin{prop}\label{prop: three singular points}
Conditions as in Set-up \ref{setup: setup with resolution}. Suppose that $b>\frac{10}{11}$ and $X$ has exactly $3$ singular points. Then $b=\frac{12}{13}$. Moreover, when $b=\frac{12}{13}$, $(X,B)$ is isomorphic to $$\left(\mathbb P(3,4,5),\frac{12}{13}(x^3y+y^2z+z^2x=0)\right)$$
and the isomorphism is unique.
\end{prop}
\begin{proof}
We check all possibilities of $\{(p_i,q(X\ni x_i,S),\gamma(X\ni x_i))\}_{i=1}^3$ in Table \ref{Table: singularity on S} such that $b>\frac{10}{11}$ and satisfy the conditions of Lemma \ref{lem: summary three singular points condition}. By running a computer program, we know they are of one of the types in Table \ref{table: 3 singular points}:
\begin{center}
    \begin{longtable}{|l|l|l|l|l|l|}
 			\caption{Possible cases of $3$-singular points}\label{table: 3 singular points}\\
		\hline
		Case & $(p_i)_{i=1}^3$ & $(\gamma(X\ni x_i))_{i=1}^3$ & $((X\ni x_i,S))_{i=1}^3$ type No. & $(q(X\ni x_i,S))_{i=1}^3$ & $b$\\ 
		\hline
1 & $(3,4,5)$  & $(2,3,\frac{8}{5})$ & $(12,16,9)$ & $(2,3,2)$ & $\frac{12}{13}$ \\
\hline
2 & $(3,5,7)$  & $(2,-\frac{4}{5},-\frac{18}{7})$ & $(12,4,6)$ & $(2,1,1)$ & $\frac{33}{34}$ \\
\hline
3 & $(3,4,7)$  & $(\frac{2}{3},0,\frac{18}{7})$ & $(2,3,19)$ & $(1,1,5)$ & $\frac{22}{23}$\\
\hline
4 & $(4,5,8)$  & $(0,\frac{8}{5},1)$ & $(3,13,8)$ & $(1,3,3)$ & $\frac{16}{17}$ \\
\hline
5 & $(5,5,8)$  & $(-\frac{4}{5},-\frac{4}{5},\frac{5}{2})$ & $(4,4,17)$ & $(1,1,5)$& $\frac{18}{19}$\\
\hline
6 & $(6,7,10)$  & $(-\frac{5}{3},\frac{6}{7},\frac{9}{5})$ & $(5,7,22)$ & $(1,2,7)$ & $\frac{29}{31}$ \\
\hline
\end{longtable}
\end{center}

\noindent\textbf{Case 1}. This case is possible. By Table \ref{table: 3 singular points}, $$(p_1,p_2,p_3,q(X\ni x_1,S),q(X\ni x_2,S),q(X\ni x_3,S))=(3,4,5,2,3,2).$$
By Lemma \ref{lem: uniquenss 12/13}, $(X,\frac{12}{13}S)$ is isomorphic to $$\left(\mathbb P(3,4,5),\frac{12}{13}(x^3y+y^2z+z^2x=0)\right)$$
and the isomorphism is unique.

\medskip

\noindent\textbf{Case 2}. We have
$$K_Y+\frac{33}{34}S_Y+\frac{22}{34}E_{1,1}+\frac{11}{34}E_{1,2}+\frac{27}{34}E_{2,1}+\frac{29}{34}E_{3,1}=f^*(K_X+bS).$$
We consider all solutions of the equation
$$c_0b+\sum_{i,j}c_{i,j}b_{i,j}=2\text{ or }3, c_0\in\mathbb N^+,c_{i,j}\in\mathbb N.$$
Then $c_0=1$ and $c_0b+\sum_{i,j}c_{i,j}b_{i,j}=3$. By Lemma \ref{lem: check s_y^2}(5), $1=c_0^2\geq S_Y^2=10$, a contradiction.

\medskip

\noindent\textbf{Case 3}. We have
$$K_Y+\frac{22}{23}S_Y+\frac{15}{23}E_{1,1}+\frac{17}{23}E_{2,1}+\frac{19}{23}E_{3,1}+\frac{16}{23}E_{3,2}+\frac{13}{23}E_{3,3}=f^*(K_X+bS).$$
We consider all solutions of the equation
$$c_0b+\sum_{i,j}c_{i,j}b_{i,j}=2\text{ or }3, c_0\in\mathbb N^+,c_{i,j}\in\mathbb N.$$
Then $c_0=1$ and $c_0b+\sum_{i,j}c_{i,j}b_{i,j}=3$. By Lemma \ref{lem: check s_y^2}(5), $1=c_0^2\geq S_Y^2=5$, a contradiction.

\medskip

\noindent\textbf{Case 4}. We have
$$K_Y+\frac{16}{17}S_Y+\frac{25}{34}E_{1,1}+\frac{13}{17}E_{2,1}+\frac{10}{17}E_{2,2}+\frac{29}{34}E_{3,1}+\frac{21}{34}E_{3,2}=f^*(K_X+bS).$$
We consider all solutions of the equation
$$c_0b+\sum_{i,j}c_{i,j}b_{i,j}=2\text{ or }3, c_0\in\mathbb N^+,c_{i,j}\in\mathbb N.$$
Then $c_0=1$ and $c_0b+\sum_{i,j}c_{i,j}b_{i,j}=3$. By Lemma \ref{lem: check s_y^2}(5), $1=c_0^2\geq S_Y^2=6$, a contradiction.

\medskip

\noindent\textbf{Case 5}. We have
$$K_Y+\frac{18}{19}S_Y+\frac{15}{19}E_{1,1}+\frac{15}{19}E_{2,1}+\frac{16}{19}E_{3,1}+\frac{14}{19}E_{3,2}+\frac{7}{19}E_{3,3}=f^*(K_X+bS).$$
We consider all solutions of the equation
$$c_0b+\sum_{i,j}c_{i,j}b_{i,j}=2\text{ or }3, c_0\in\mathbb N^+,c_{i,j}\in\mathbb N.$$
Then $c_0=1$ or $2$ and $c_0b+\sum_{i,j}c_{i,j}b_{i,j}=3$. By Lemma \ref{lem: check s_y^2}(5), $4\geq c_0^2\geq S_Y^2=8$, a contradiction.

\medskip

\noindent\textbf{Case 6}. We have
$$K_Y+\frac{29}{31}S_Y+\frac{51}{62}E_{1,1}+\frac{26}{31}E_{2,1}+\frac{13}{31}E_{2,2}+\frac{53}{62}E_{3,1}+\frac{24}{31}E_{3,2}+\frac{43}{62}E_{3,3}=f^*(K_X+bS).$$
We consider all solutions of the equation
$$c_0b+\sum_{i,j}c_{i,j}b_{i,j}=2\text{ or }3, c_0\in\mathbb N^+,c_{i,j}\in\mathbb N.$$
Then $c_0=1$ and $c_0b+\sum_{i,j}c_{i,j}b_{i,j}=3$. By Lemma \ref{lem: check s_y^2}(5), $1=c_0^2\geq S_Y^2=8$, a contradiction.
\end{proof}

\subsection{Four singular points case, I}

\begin{lem}\label{lem: summary four singular points condition, I}
Conditions as in Set-up \ref{setup: rho=1 set-up}. Suppose that $b>\frac{10}{11}$, $X$ has exactly $4$ singular points, and $x_4\in S$. Then \begin{enumerate}
\item $(X\ni x_i,S)$ is in Table \ref{Table: singularity on S} for each $i$, 
\item $1\leq \sum_{i=1}^3\frac{1}{p_i}<2$,
\item $\sum_{i=1}^4\gamma(X\ni x_i)=9-\frac{b^2}{1-b}(2-\sum_{i=1}^4\frac{1}{p_i})<9-\frac{100}{11}(2-\sum_{i=1}^4\frac{1}{p_i})$,
\item $\sqrt{(9-\sum_{i=1}^4\gamma(X\ni x_i))^2+4(9-\sum_{i=1}^4\gamma(X\ni x_i)(2-\sum_{i=1}^4\frac{1}{p_i}))}\in\mathbb Q$,
\item $\frac{1}{1-b}(2-\sum_{i=1}^4\frac{1}{p_i})-\sum_{i=1}^4\frac{q(X\ni x_i,S)}{p_i}\in\mathbb Z$, and
\item $(X\ni x_i,bS)$ is $\frac{1}{7}$-lt for for any $i$.
\end{enumerate}
\end{lem}
\begin{proof}
(1)(6) follow immediately from our conditions. (2) follows from Lemma \ref{lem: S contain at least 3 sing} and Theorem \ref{thm: euler bound}. (3) follows from  Lemma \ref{lem: gamma(X) equality}(1) and the assumption that $b>\frac{10}{11}$. (4) follows from Lemma \ref{lem: gamma(X) equality}(2). (5) follows from Lemma \ref{lem: control sy^2}.
\end{proof}

\begin{prop}\label{prop: four singularities, i}
Conditions as in Set-up \ref{setup: rho=1 set-up}. Then the case when $b>\frac{10}{11}$, $X$ has exactly $4$ singular points, and $x_4\in S$ is not possible.
\end{prop}
\begin{proof}
We check all possibilities of $\{(p_i,q(X\ni x_i,S),\gamma(X\ni x_i))\}_{i=1}^4$ in Table \ref{Table: singularity on S} such that $b>\frac{10}{11}$ and satisfy the conditions of Lemma \ref{lem: summary four singular points condition, I}. By running a computer program, we find that there is no possible case.
\end{proof}

\subsection{Four singular points case, II}

In this subsection, we prove the following proposition:

\begin{prop}\label{prop: four singularities, ii}
Conditions as in Set-up \ref{setup: rho=1 set-up}. Then the case when $b>\frac{10}{11}$, $X$ has exactly $4$ singular points, and $x_4\not\in S$ is not possible.
\end{prop}

We need the following lemmas in the proof of Proposition \ref{prop: four singularities, ii}.

\begin{lem}\label{lem: special curve have negative self-intersection}
Let $(X,B)$ be a projective klt log Calabi-Yau surface pair such that $X$ is non-singular. Let $C$ be curve on $X$ such that $C$ is not a component of $\Supp B$ and $C$ intersects $\Supp B$.  Let $f: X\rightarrow Y$ be a projective birational morphism that contracts $C$. Then $C$ is a $(-1)$-curve.
\end{lem}
\begin{proof}
Since $(K_X+B)\cdot C=0$, $C\not\subset\Supp B$, and $C$ intersects $\Supp B$, then $B\cdot C>0$ and $K_X\cdot C<0$. Since $f$ is a birational contraction, $C^2<0$. Thus $K_X\cdot C=C^2=-1$ and $C$ is a $(-1)$-curve.
\end{proof}

\begin{lem}\label{lem: two ray game four sings}
Conditions as in Set-up \ref{setup: setup with resolution}. Suppose that $b>\frac{10}{11}$, $X$ has exactly $4$ singular points, $x_4$ is not a Du Val singularity, and $x_4\not\in S$. Suppose that there exists an index $j_0$ such that $a(E_{4,j_0},X,0)=\mld(X;x_4)$. Then either there exists a positive integer $n$ such that
$$2=b+nb_{4,j_0},$$
or there exists a non-singular rational curve $C$ on $Y$ satisfying the following.
\begin{enumerate}
    \item $C^2=-1$.
    \item $C\not=S_Y$ and $C\not=E_{i,j}$ for any $i,j$.
    \item $C$ intersects $S_Y\cup\cup_{1\leq i\leq 3}E_{i,j}$ and $C$ intersects $\cup_jE_{4,j}$.
    \item The intersection matrix of any connected component of $\mathcal{D}(C\cup\cup_{(i,j)\not=(4,j_0)}E_{i,j})$ is negative definite.
    \item Let $\Ii:=\{(i,j)\mid E_{i,j}\text{ intersects }C\}$. Then there exist positive integers $c_{i,j}$ and $c$, such that
    $$cb+\sum_{(i,j)\in\Ii}c_{i,j}b_{i,j}=1.$$
    Moreover, $c=1$ if $S_Y$ intersects $C$, and $c=0$ if $S_Y$  does not intersect $C$.
\end{enumerate}
\end{lem}
\begin{proof}
Since $a(E_{4,j_0},X,0)=\mld(X;x_4)$ and $x_4$ is not Du Val, $a(E_{4,j_0},X,0)\leq\frac{2}{3}$ (cf. \cite[Lemma 5.1]{Jia21}), so $b_{4,j_0}\geq\frac{1}{3}$. We let $h: Y'\rightarrow X$ be the extraction of $E_{4,j_0}$. We denote $S_{Y'}$ and $E_{4,j_0}'$ the images of $S_Y$ and $E_{4,j_0}$ on $Y'$, and let $g: Y\rightarrow Y'$ be the induced birational contraction. Then $\rho(Y')=2$. We let $R$ be the unique extremal ray $\overline{NE}(Y')$ that is not $[E'_{4,j_0}]$ and let $\pi: Y'\rightarrow Z$ be the divisorial contraction of $R$. There are two cases.

\medskip

\noindent\textbf{Case 1}. $\pi$ is a Mori fiber space. We let $F$ be a general fiber of $\pi$, $S_F:=S_{Y'}|_F$, and $E_F:=E_{4,j_0}'|_F$. Since $(-S_{Y'})\cdot E'_{4,j_0}=0$ and $-S_{Y'}$ is not pseudo-effective, $R\cdot (-S_{Y'})<0$. Thus $S_{Y'}\cdot R>0$ and $S_F\not=0$. Hence $F\cong\mathbb P^1$ and
$$2=b\deg S_F+b_{4,j_0}\deg E_F.$$
Since $b\in(\frac{10}{11},1)$, $\deg E_F>0$. Since $b\in(\frac{10}{11},1)$ and $1>b_{4,j_0}\geq\frac{1}{3}$, $\deg S_F=1$. We may let $n=\deg E_F$.

\medskip

\noindent\textbf{Case 2}. $\pi$ is not a Mori fiber space. Then we let $C'$ be the unique rational curve contracted by $\pi$ and let $C$ be the strict transform of $C$ on $Y$. We show that $C$ satisfies our requirement.

(1) It follows from Lemma \ref{lem: special curve have negative self-intersection} as $C$ is contracted by $\pi\circ g: Y\rightarrow Z$.

(2) It follows from our construction.

(3) We let $E_Z$ and $S_Z$ be the images of $E_{4,j_0}'$ and $S_{Y'}$ on $Z$ respectively. Since $\rho(Z)=1$, $E_Z$ intersects $S_Z$. Since $E_{4,j_0}'$ and $S_{Y'}$ do not intersect, $C'$ intersects $E_{4,j_0}'$ and $S_{Y'}$. Since $g$ only extracts $E_{i,j}$ and $S_Y\cup\cup_{1\leq i\leq 3}E_{i,j}$ and $\cup_jE_{4,j}$ do not intersect, we get (3).

(4) It follows from the fact that $\pi\circ g$ is a projective birational morphism between klt varieties which contracts $\mathcal{D}(C\cup\cup_{(i,j)\not=(4,j_0)}E_{i,j})$.

(5) By (1), $C^2=-1$ and $K_Y\cdot C=-1$. Since 
$$K_Y+bS_Y+\sum_{i,j} b_{i,j}E_{i,j}\equiv 0,$$
we have
$$1=b(S_Y\cdot C)+\sum_{i,j} b_{i,j}(E_{i,j}\cdot C)=b(S_Y\cdot C)+\sum_{(i,j)\in\Ii} b_{i,j}(E_{i,j}\cdot C).$$
We $c_{i,j}:=E_{i,j}\cdot C$ for any $(i,j)\in\Ii$ and let let $c:=S_Y\cdot C$ if $S_Y$ intersects $C$. Then (5) follows.
\end{proof}

\begin{lem}\label{lem: summary four singular points condition, II}
Conditions as in Set-up \ref{setup: rho=1 set-up}. Suppose that $b>\frac{10}{11}$, $X$ has exactly $4$ singular points, and $x_4\not\in S$. Then \begin{enumerate}
\item $(X\ni x_i,S)$ is in Table \ref{Table: singularity on S} for each $i$, 
\item $1\leq \sum_{i=1}^4\frac{1}{p_i}$, $\sum_{i=1}^3\frac{1}{p_i}<1$, and $p_4\leq 42$.
\item $\sum_{i=1}^4\gamma(X\ni x_i)=9-\frac{b^2}{1-b}(1-\sum_{i=1}^3\frac{1}{p_i})<9-\frac{100}{11}(1-\sum_{i=1}^3\frac{1}{p_i})$,
\item $\sqrt{(9-\sum_{i=1}^4\gamma(X\ni x_i))^2+4(9-\sum_{i=1}^4\gamma(X\ni x_i)(1-\sum_{i=1}^3\frac{1}{p_i}))}\in\mathbb Q$,
\item $\frac{1}{1-b}(1-\sum_{i=1}^3\frac{1}{p_i})-\sum_{i=1}^3\frac{q(X\ni x_i,S)}{p_i}\in\mathbb Z$, 
\item $(X\ni x_i,bS)$ is $\frac{1}{7}$-lt for for any $i$, and
\item $\gamma(X\ni x_4)<\frac{887}{77}$.
\end{enumerate}
\end{lem}
\begin{proof}
(1)(6) follow immediately from our conditions. By Lemma \ref{lem: S contain at least 3 sing} and Theorem \ref{thm: euler bound}, $1\leq \sum_{i=1}^4\frac{1}{p_i}$ and $\sum_{i=1}^3\frac{1}{p_i}<1$. Thus  $\sum_{i=1}^3\frac{1}{p_i}\leq\frac{41}{42}$, so $p_4\leq 42$, and we get (2). (3) follows from  Lemma \ref{lem: gamma(X) equality}(1) and the assumption that $b>\frac{10}{11}$. (4) follows from Lemma \ref{lem: gamma(X) equality}(2). (5) follows from Lemma \ref{lem: control sy^2}. We check all possibilities of $\{(p_i,\gamma(X\ni x_i))\}_{i=1}^3$ in Table \ref{Table: singularity on S} by running a computer program, and get $9-\sum_{i=1}^3\gamma(X\ni x_i)-\frac{100}{11}(1-\sum_{i=1}^3\frac{1}{p_i})\leq\frac{887}{77}$, so (7) follows from (3). 
\end{proof}

\begin{proof}[Proof of Proposition \ref{prop: four singularities, ii}]
We check all possibilities of $\{(p_i,q(X\ni x_i,S),\gamma(X\ni x_i))\}_{i=1}^3$ in Table \ref{Table: singularity on S} and $(p_4,\gamma(X\ni x_4))$ such that $b>\frac{10}{11}$ and satisfy the conditions of Lemma \ref{lem: summary four singular points condition, II}. Note there are only finitely many possibilities of $x_4$ as $p_4\leq 42$, so we can enumerate all of them (see Appendix \ref{sec: classification of x4}
 for details). By running a computer program, we know they are only of types in Table \ref{table: 4 singular points, II}. In particular, $x_4$ is a cyclic quotient singularity, so we may let $f,Y,B_Y,S_Y,E_{i,j},b_{i,j}$ be objects satisfying the conditions as in Set-up \ref{setup: setup with resolution}.
 
\begin{center}
    \begin{longtable}{|l|l|l|l|l|l|l|}
 			\caption{Possible cases of $4$-singular points}\label{table: 4 singular points, II}\\
		\hline
		Case & $(p_i)_{i=1}^4$ & $x_4$ &$(q(X\ni x_i,S))_{i=1}^3$ & $((X\ni x_i,S))_{i=1}^3$ type No. & $(\gamma(X\ni x_i))_{i=1}^4$ & $b$\\ 
		\hline
1 & $(2,3,7,19)$  & $\frac{1}{19}(1,6)$ & $(1,2,2)$ & $(1,12,7)$ & $(1,2,\frac{6}{7},\frac{90}{19})$ & $\frac{18}{19}$\\
\hline
2 & $(2,3,7,31)$  & $\frac{1}{31}(1,22)$ & $(1,2,4)$ & $(1,12,14)$ & $(1,2,\frac{6}{7},\frac{138}{31})$ & $\frac{30}{31}$\\
\hline
3 & $(2,3,7,29)$ & $\frac{1}{29}(1,16)$ & $(1,1,6)$ & $(1,2,23)$ & $(1,\frac{2}{3},6,\frac{20}{29})$ & $\frac{28}{29}$ \\
\hline
4 & $(2,5,7,3)$ & $\frac{1}{3}(1,1)$ & $(1,2,2)$ & $(1,9,7)$ & $(1,\frac{8}{5},\frac{6}{7},\frac{2}{3})$ & $\frac{32}{33}$\\
\hline
5 & $(3,4,5,4)$ & $\frac{1}{4}(1,1)$& $(2,3,2)$ & $(12,16,9)$ &$(2,3,\frac{8}{5},0)$ & $\frac{12}{13}$ \\
\hline
6 & $(3,5,7,2)$ & $\frac{1}{2}(1,1)$& $(1,3,4)$ & $(2,13,14)$ &$(\frac{2}{3},\frac{8}{5},\frac{6}{7},1)$ & $\frac{16}{17}$ \\
\hline
7 & $(4,5,6,2)$ & $\frac{1}{2}(1,1)$& $(1,2,1)$ & $(3,9,5)$ &$(0,\frac{8}{5},-\frac{5}{3},1)$ & $\frac{22}{23}$ \\
\hline
\end{longtable}
\end{center}

\noindent\textbf{Case 1}. We have
$$K_Y+\frac{18}{19}S_Y+\frac{9}{19}E_{1,1}+\frac{12}{19}E_{2,1}+\frac{6}{19}E_{2,2}+\frac{16}{19}E_{3,1}+\frac{8}{19}E_{3,2}+\sum_{j=1}^6\frac{14-2j}{19}E_{4,j}=f^*(K_X+bS).$$
This contradicts Lemma \ref{lem: two ray game four sings}. More precisely, since $b=\frac{18}{19}$ and $\mld(X; x_4)=\frac{7}{19}$, no positive integer $n$ satisfies
$$2=b+(1-\mld(X;x_4))n.$$
By Lemma \ref{lem: two ray game four sings}(1)(2)(5), there exists a non-singular rational curve $C$ on $Y$ such that $C^2=-1$, $C\not=S_Y$, $C\not=E_{i,j}$ for any $i,j$, and there exist non-negative integers $c$ and $c_{i,j}$ such that
$$cb+\sum c_{i,j}b_{i,j}=1$$
and $C$ intersects $E_{i,j}$ if and only if $c_{i,j}>0$. Thus $c=0$, $c_{1,1}=1$, and $c_{4,j}\geq 1$ for some integer $j\geq 2$. This contradicts Lemma \ref{lem: two ray game four sings}(4) as the intersection matrix of the dual graph of $C\cup\cup_{j=2}^5 E_{4,j}\cup E_{1,1}$ is not negative definite in any case.

\medskip

\noindent\textbf{Case 2}. We have
\begin{align*}
    K_Y+&\frac{30}{31}S_Y+\frac{15}{31}E_{1,1}+\frac{20}{31}E_{2,1}+\frac{10}{31}E_{2,2}+\frac{26}{31}E_{3,1}\\
    +&\frac{22}{31}E_{3,2}+\frac{8}{31}E_{4,1}+\frac{16}{31}E_{4,2}+\sum_{i=1}^4\frac{30-6i}{31}E_{4,2+i}=f^*(K_X+bS).
\end{align*}
This contradicts Lemma \ref{lem: two ray game four sings}. More precisely, since $b=\frac{30}{31}$ and $\mld(X;x_4)=\frac{7}{31}$, no positive integer $n$ satisfies
$$2=b+(1-\mld(X;x_4))n.$$
By Lemma \ref{lem: two ray game four sings}(1)(2)(5), there exists a non-singular rational curve $C$ on $Y$ such that $C^2=-1$, $C\not=S_Y$, $C\not=E_{i,j}$ for any $i,j$, and there exist non-negative integers $c$ and $c_{i,j}$ such that
$$cb+\sum c_{i,j}b_{i,j}=1$$
and $C$ intersects $E_{i,j}$ if and only if $c_{i,j}>0$. Thus $c=0$, and one of the following holds:
\begin{itemize}
    \item $c_{1,1}=1$, $c_{4,1}=2$, and all other $c_{i,j}=0$. In this case, the dual graph of $C\cup E_{4,1}\cup E_{1,1}$ is not negative definite. 
    \item $c_{1,1}=1$, $c_{4,2}=1$, and all other $c_{i,j}=0$.  In this case, the dual graph of $C\cup E_{4,2}\cup E_{1,1}$ is not negative definite. 
\end{itemize}
Thus either case contradicts Lemma \ref{lem: two ray game four sings}(4).

\medskip

\noindent\textbf{Case 3}. We have
$$K_Y+\frac{28}{29}S_Y+\frac{14}{29}E_{1,1}+\frac{19}{29}E_{2,1}+\sum_{i=1}^6\frac{28-4i}{29}E_{3,i}+\frac{12}{29}E_{4,1}+\sum_{i=1}^3\frac{32-8i}{29}E_{4,i+1}=f^*(K_X+bS).$$
This contradicts Lemma \ref{lem: two ray game four sings}. More precisely, since $b=\frac{28}{29}$ and $\mld(X;x_4)=\frac{5}{29}$, no positive integer $n$ satisfies
$$2=b+(1-\mld(X;x_4))n.$$
By Lemma \ref{lem: two ray game four sings}(2)(5), there exists a non-singular rational curve $C$ on $Y$ such that $C\not=S_Y$, $C\not=E_{i,j}$ for any $i,j$, and there exist non-negative integers $c$ and $c_{i,j}$ such that
$$cb+\sum c_{i,j}b_{i,j}=1$$
and $C$ intersects $E_{i,j}$ if and only if $c_{i,j}>0$. Thus $c=0,c_{4,j}=0$ for any $j$, which contradicts Lemma \ref{lem: two ray game four sings}(3).

\medskip

\noindent\textbf{Case 4}. We have 
$$K_Y+\frac{32}{33}S_Y+\frac{16}{33}E_{1,1}+\frac{26}{33}E_{2,1}+\frac{13}{33}E_{2,2}+\frac{28}{33}E_{3,1}+\frac{14}{33}E_{3,2}+\frac{1}{3}E_{4,1}=f^*(K_X+bS).$$
This contradicts Lemma \ref{lem: two ray game four sings}. More precisely, since $b=\frac{32}{33}$ and $\mld(X;x_4)=\frac{2}{3}$, no positive integer $n$ satisfies
$$2=b+(1-\mld(X;x_4))n.$$
By Lemma \ref{lem: two ray game four sings}(2)(5), there exists a non-singular rational curve $C$ on $Y$ such that $C\not=S_Y$, $C\not=E_{i,j}$ for any $i,j$, and there exist non-negative integers $c$ and $c_{i,j}$ such that
$$cb+\sum c_{i,j}b_{i,j}=1$$
and $C$ intersects $E_{i,j}$ if and only if $c_{i,j}>0$. Thus $c=0$, $c_{4,1}=3$, and all other $c_{i,j}=0$. This contradicts Lemma \ref{lem: two ray game four sings}(3).

\medskip

\noindent\textbf{Case 5}. We have
$$K_Y+\frac{12}{13}S_Y+\frac{8}{13}E_{1,1}+\frac{4}{13}E_{1,2}+\frac{9}{13}E_{2,1}+\frac{6}{13}E_{2,2}+\frac{3}{13}E_{2,3}+\frac{10}{13}E_{3,1}+\frac{5}{13}E_{3,2}+\frac{1}{2}E_{4,1}=f^*(K_X+bS)$$
which contradicts Lemma \ref{lem: two ray game four sings}. More precisely, since $b=\frac{12}{13}$ and $\mld(X;x_4)=\frac{1}{2}$, no positive integer $n$ satisfies 
$$2=b+(1-\mld(X;x_4))n.$$
By Lemma \ref{lem: two ray game four sings}(2)(5), there exists a non-singular rational curve $C$ on $Y$ such that $C\not=S_Y$, $C\not=E_{i,j}$ for any $i,j$, and there exist non-negative integers $c$ and $c_{i,j}$ such that
$$cb+\sum c_{i,j}b_{i,j}=1$$
and $C$ intersects $E_{i,j}$ if and only if $c_{i,j}>0$. Thus $c=0$, and either 
\begin{itemize}
    \item $c_{4,1}=2$ and all other $c_{i,j}=0$, or
    \item $c_{4,1}=0$.
\end{itemize}
But both cases contradict Lemma \ref{lem: two ray game four sings}(3).

\medskip

\noindent\textbf{Case 6}. We have
$$K_Y+\frac{16}{17}S_Y+\frac{11}{17}E_{1,1}+\frac{12}{17}E_{2,1}+\frac{10}{17}E_{2,2}+\frac{14}{17}E_{3,1}+\frac{12}{17}E_{3,2}=f^*(K_X+bS)$$
and $S_Y^2=4$. We consider all solutions of the equation
$$c_0b+\sum_{i,j}c_{i,j}b_{i,j}=2\text{ or }3, c_0\in\mathbb N^+,c_{i,j}\in\mathbb N.$$
Then $c_0=1$ and $c_0b+\sum_{i,j}c_{i,j}b_{i,j}=3$. By Lemma \ref{lem: check s_y^2}(5), $1=c_0^2\geq S_Y^2=4$, a contradiction.

\medskip

\noindent\textbf{Case 7}. We have
$$K_Y+\frac{22}{23}S_Y+\frac{17}{23}E_{1,1}+\frac{18}{23}E_{2,1}+\frac{9}{23}E_{2,2}+\frac{19}{23}E_{3,1}=f^*(K_X+bS)$$
and $S_Y^2=8$. We consider all solutions of the equation
$$c_0b+\sum_{i,j}c_{i,j}b_{i,j}=2\text{ or }3, c_0\in\mathbb N^+,c_{i,j}\in\mathbb N.$$
Then $c_0=1$ and $c_0b+\sum_{i,j}c_{i,j}b_{i,j}=3$. By Lemma \ref{lem: check s_y^2}(5), $1=c_0^2\geq S_Y^2=8$, a contradiction.
\end{proof}

\section{Proof of the main theorems}\label{sec: proof of the main theorems}

In this section, we prove all our main theorems except Theorem \ref{thm: 1/13 exceptional Fano} and Corollary \ref{cor: new bound alpha invariant}.

\subsection{Proof of Theorem \ref{thm: global lct 12/13 intro}}
\begin{proof}[Proof of Theorem \ref{thm: global lct 12/13}] (2) follows from Lemma \ref{lem: reduce to smooth rational curve}(1), and (3.b) follows from (3.a) and Lemma \ref{lem: reduced to picard number 1}. Thus we only left to prove (1), the main part of (3), and (3.a). 

By Lemmas \ref{lem: reduce to smooth rational curve} and possibly replacing $b$ with $b_1$, we may assume that $X,B,b$, and $S$ satisfy Set-up \ref{setup: klt set-up}. By Lemma \ref{lem: reduce to picard number 1}, we may assume that $X,B,b$, and $S$ satisfy Set-up \ref{setup: rho=1 set-up}. Then one of the following holds.
\begin{enumerate}
    \item $X$ has exactly three singular points, and all of them are in $S$. By Notation-Lemma \ref{notalem: q(x,s)}, $X,B,b$, and $S$ satisfy Set-up \ref{setup: setup with resolution}. In this case, the theorem follows from  Proposition \ref{prop: three singular points}.
    \item $X$ has exactly four singular points, and all of them are in $S$.  By Notation-Lemma \ref{notalem: q(x,s)}, $X,B,b$, and $S$ satisfy Set-up \ref{setup: setup with resolution}. In this case, the theorem follows from  Proposition \ref{prop: four singularities, i}.
    \item $X$ has exactly four singular points, and one singular point of $X$ is not contained in $S$. In this case, the theorem follows from Proposition \ref{prop: four singularities, ii}.
\end{enumerate}
\end{proof}
\begin{proof}[Proof of Theorem \ref{thm: global lct 12/13 intro}] It is a special case of Theorem \ref{thm: global lct 12/13}.
\end{proof}

The following corollary can be viewed as the relative version of Theorem \ref{thm: global lct 12/13 intro}.
\begin{cor}\label{cor: relative version 12/13}
Let $(X,T+bS)$ be a projective lc surface pair and $X\rightarrow Z$ a proper surjective morphism, such that \begin{enumerate}
    \item $T$ is an effective Weil divisor,
    \item $S$ is a non-zero effective Weil divisor that is horizontal over $Z$, and
    \item $K_X+T+bS\equiv_Z0$.
\end{enumerate} 
Then $b=1$, or $b=\frac{12}{13}$, or $b\leq\frac{10}{11}$. Moreover, if $b\in (\frac{2}{3},1)$, then $Z$ is a point.
\end{cor}
\begin{proof}
Since $S$ is horizontal over $Z$, $\dim Z=0$ or $1$. If $\dim Z=0$, then the corollary follows from Theorem \ref{thm: global lct 12/13 intro}. If $\dim Z=1$, then we let $F$ be a general fiber of $X\rightarrow Z$, $T_F:=T|_F$, and $S_F:=S|_F$. Then $K_F+T_F+bS_F\equiv 0$. Thus $b\deg(S_F)=1$ or $2$, so $b=1$ or $b\leq\frac{2}{3}$.
\end{proof}

\subsection{Proof of Theorem \ref{thm: 12/13 rct intro} }

First, we state a more general version of Theorem \ref{thm: 12/13 rct intro}:

\begin{thm}[First and second gaps of $\Rr$-complementary thresholds, general version]\label{thm: 12/13 rct}
Let $(X,T)$ be a surface pair and $X\rightarrow Z$ a proper surjective morphism such that $(X/Z,T)$ is $\Rr$-complementary and $T$ is an effective Weil divisor. Let $S$ be a non-zero effective $\Qq$-Cartier Weil divisor on $X$. Let
$$b:=\Rct(X/Z,T;S):=\sup\{t\geq 0\mid (X,T+tS)\text{ is }\Rr\text{-complementary}\}.$$
Then either $b=1$, or $b=\frac{12}{13}$, or $b\leq\frac{10}{11}$. Moreover, if $b\in [\frac{6}{7},1)$, then $Z$ is a point and $T=0$.
\end{thm}
\begin{proof}
The proof is similar to \cite[Proof of Theorem 5.20]{HLS19} and \cite[Lemma 14]{Sho20} so we only provide a sketch of the proof here. By \cite[2.3 Inductive Theorem, 3.1 Theorem]{Sho00}, if $Z\not=\{pt\}$ or $T\not=0$, then $(X/Z,T+\frac{6}{7}S)$ has an $n$-complement for some $n\in\{1,2,3,4,6\}$. Then $b=1$ or $b<\frac{6}{7}$ and we are done. Thus we may assume that $Z=\{pt\}$ and $b\geq\frac{6}{7}$.

By the ACC for $\Rr$-complementary thresholds \cite[Theorem 8.20]{HLS19}, \cite[Theorem 21]{Sho20}, $b$ is achieved, i.e. the supremum is a maximum. Note that although $X$ may not be of Fano type, the ACC for $\Rr$-complementary thresholds still holds in our case because the existence of good minimal models holds in dimension $2$ (cf. \cite[Theorem 5.20]{HLS19}). By \cite[Lemma 14]{Sho20}, $\Rct(X,T;S)$ is either equal to the lc threshold $\lct(X,T;S)$, or equal to the anti-canonical threshold $\act(X,T;S)$. By Lemma \ref{lem: surface lct 5/6}, we may assume that  $\Rct(X,T;S)$ is equal to the anti-canonical threshold $\act(X,T;S)$ and $(X,T+bS)$ is lc. Let that $X\rightarrow Z'$ be a contraction defined by $-(K_X+T+bS)$. Then $K_X+T+bS\equiv 0$. Since $b=\act(X,T;S)$, $S$ dominates $Z'$. The theorem follows from Corollary \ref{cor: relative version 12/13}. 
\end{proof}

\begin{proof}[Proof of Theorem \ref{thm: 12/13 rct intro}]
It is a special case of Theorem \ref{thm: 12/13 rct}.
\end{proof}

\subsection{Proof of Theorem \ref{thm: 1/13 exceptional} and Corollary \ref{cor: calabi-yau surface mld}}

\begin{thm}\label{thm: 1/13 exceptional strong}
Let $X$ be an exceptional surface. Then for any $B\in |-K_X|_{\mathbb Q}$, $(X,B)$ is $\frac{1}{13}$-lc.
\end{thm}
\begin{proof}
Suppose that there exists $B\in |-K_X|_{\mathbb Q}$ such that $(X,B)$ is not $\frac{1}{13}$-lc. Since $X$ is exceptional, $(X,B)$ is klt, so there exists a divisor $E$ over $X$ such that $0<a:=a(E,X,B)<\frac{1}{13}$. We let $f: Y\rightarrow X$ be the identity morphism if $E$ is on $X$, and let $f$ be an extraction of $E$ if $E$ is exceptional over $X$. Let $K_Y+B_Y:=f^*(K_X+B)$. Then $(Y,B_Y)$ is klt Calabi-Yau. Since $B_Y\geq (1-a)E$, $(Y,(1-a)E)$ is $\Rr$-complementary. Thus $\Rct(Y,0;E)\geq 1-a>\frac{12}{13}$. By Theorem \ref{thm: 12/13 rct}, $\Rct(Y,0;E)=1$. Thus $(Y,E)$ has an $\Rr$-complement $(Y,E+G_Y)$ for some $G_Y\geq 0$. Let $G:=f_*(E+G_Y)$, then $(X,G)$ is lc log Calabi-Yau and $a(E,X,G)=0$. This contradicts the assumption that $X$ is exceptional.
\end{proof}

\begin{proof}[Proof of Theorem \ref{thm: 1/13 exceptional}]
It follows from Theorem \ref{thm: 1/13 exceptional strong}.
\end{proof}

\begin{proof}[Proof of Corollary \ref{cor: calabi-yau surface mld}]
By abundance (cf. \cite[V.4.6 Theorem]{Nak04}), any  klt Calabi-Yau variety is exceptional. So Corollary \ref{cor: calabi-yau surface mld} is a special case of Theorem \ref{thm: 1/13 exceptional}.
\end{proof}

We need the following result in order to prove Theorem \ref{thm: 1/13 exceptional Fano}.

\begin{lem}\label{lem: exceptional fano surface 1/13 klt}
Exceptional del Pezzo surfaces are $\frac{1}{13}$-lt.
\end{lem}
\begin{proof}
Let $X$ be an exceptional del Pezzo surface and $E$ a prime divisor over $X$ such that $a(E,X,0)=\mld(X)$. We pick $B\in |-K_X|_{\mathbb Q}$ such that $\Center_XE$ is contained in $\Supp B$. By Theorem \ref{thm: 1/13 exceptional strong}, $(X,B)$ is $\frac{1}{13}$-lc, so $\mld(X)=a(E,X,0)>a(E,X,B)\geq\frac{1}{13}$.
\end{proof}

\subsection{Proof of Theorem \ref{thm: vol 1/462}}

\begin{lem}\label{lem: plt germ mld}
Let $p$ be a positive integer and $(X\ni x,S)$ a plt surface pair such that the order of $X\ni x$ is $p$. Then $\mld(X,S;x)=\frac{1}{p}$.
\end{lem}
\begin{proof}
Let $K_S+B_S:=(K_X+S)|_S$, then by \cite[3.9 Proposition]{Sho92} and inversion of adjunction,
$$\frac{1}{p}=\mld(S,B_S;x)=\mld(X,S;x).$$
\end{proof}

\begin{defn}\label{defn: pet}
Let $(X,B)$ be a projective lc pair and $D\geq 0$ an $\Rr$-Cartier $\Rr$-divisor on $X$. We define
$$\pet(X,B;D):=\inf\{t\geq 0\mid (X,B+tD)\text{ is lc},  K_X+B+tD\text{ is pseudo-effective}\}.$$
\end{defn}

\begin{thm}\label{thm: 12/13 pet}
Let $(X,T)$ be a $\Qq$-factorial projective dlt surface pair such that $T$ and $S$ are effective Weil divisors and $\pet(X,T;S)\in [\frac{10}{11},1)$. Then:
\begin{enumerate}
    \item $\pet(X,T;S)=\frac{12}{13}$ or $\frac{10}{11}$.
    \item If $\pet(X,T;S)=\frac{12}{13}$, then $S$ is a non-singular rational curve.
    \item If $K_X+T+S$ is nef, then: \begin{enumerate}
        \item $T=0$.
        \item If $\pet(X,T;S)=\frac{12}{13}$, then $S$ only contains cyclic quotient singularities of $X$ of indices $\leq 5$.
    \end{enumerate}
\end{enumerate}
\end{thm}
\begin{proof}
First we prove (1). We let $b:=\pet(X,T;S)$. We may assume that $b>\frac{10}{11}$. Pick $\epsilon\in (0,b-\frac{10}{11})$ such that $\frac{12}{13}\not\in(b-\epsilon,b)$. We run a $(K_X+T+(b-\epsilon)S)$-MMP, which terminates with a Mori fiber space $f: X'\rightarrow Z$. Let $T'$ and $S'$ be the images of $T$ and $S$ on $X'$ respectively. Since $(X',T'+(b-\epsilon)S')$ is lc, $(X,T'+\frac{10}{11}S')$ is lc. By Lemma \ref{lem: surface lct 5/6}, $(X',T'+S')$ is lc. Since $K_X+T+bS$ is pseudo-effective, $K_{X'}+T'+bS'$ is pseudo-effective, so $K_{X'}+T'+bS'$ is nef$/Z$. Since $K_{X'}+T'+(b-\epsilon)S'$ is anti-ample$/Z$, $S'$ is ample$/Z$, and there exists $b'\in (b-\epsilon,b]$ such that $K_{X'}+T'+b'S'\equiv_Z0$. In particular, $S'$ is horizontal$/Z$. By Corollary \ref{cor: relative version 12/13}, $Z$ is a point, $\rho(X')=1$, and $b'=b=\frac{12}{13}$. This implies (1).

Next we prove (2). Since $\rho(X')=1$, by Theorem \ref{thm: global lct 12/13}, $T'=0$. By Lemma \ref{lem: reduce to smooth rational curve}, $S'$ is a non-singular rational curve, $(X',S')$ is plt, and $S'$ only contains cyclic quotient singularities of $X$. This implies (2). Moreover, by Theorem \ref{thm: global lct 12/13}, $S'$ contains three singularities of $X$ with indices $\leq 5$.

Finally, we prove (3). Suppose that $K_X+T+S$ is nef. By the negativity lemma, $a(E,X,T+S)\geq a(E,X',S')$ for all exceptional and non-exceptional prime divisors $E$ over $X$. Since $(X',S')$ is plt, $T=0$ and $(X,S)$ is plt, and we get (3.a). Moreover, if $b=\frac{12}{13}$, then by Theorem \ref{thm: global lct 12/13}, $S'$ contains three singularities with indices $\leq 5$. By Lemma \ref{lem: plt germ mld}, $a(E,X',S')\geq\frac{1}{5}$ for any exceptional or non-exceptional prime divisor $E$ over $X'$. Thus $a(E,X,S)\geq\frac{1}{5}$ for any exceptional or non-exceptional prime divisor $E$ over $X$. Since $(X,S)$ is plt, by Lemma \ref{lem: de type mld}, $S$ only contains cyclic quotient singularities of $X$. By Lemma \ref{lem: plt germ mld}, $S$ only contains cyclic quotient singularities of $X$ of indices $\leq 5$, and we get (3.b).
\end{proof}

\begin{lem}\label{lem: pet small case}
Let $(X,B=T+S)$ be a $\Qq$-factorial projective dlt surface pair such that $T$ and $S$ are effective Weil divisors. Suppose that $K_X+B$ is nef, $\pet(X,T;S)\not=\frac{12}{13}$, and $(K_X+B)\cdot S>0$. Then:
\begin{enumerate}
    \item $(K_X+B)^2\geq\frac{1}{462}$.
    \item If $(K_X+B)^2=\frac{1}{462}$, then $T=0$ and $S$ is a non-singular rational curve.
\end{enumerate}
\end{lem}
\begin{proof}
Let $K_S+B_S:=(K_X+B)|_S$, where $B_S$ is the different of $T$ on $S$. We have
$$0<(K_X+B)\cdot S=\deg(K_X+B)|_S=\deg(K_S+B_S).$$
Since $B_S$ only has standard coefficients, an elementary computation shows that $\deg(K_S+B_S)\geq\frac{1}{42}$, and if $\deg(K_S+B_S)=\frac{1}{42}$, then $S$ is a non-singular rational curve.

Thus $(K_X+B)\cdot S\geq\frac{1}{42}$. By Theorem \ref{thm: 12/13 pet}, $K_X+T+\frac{10}{11}S$ is pseudo-effective. Thus
$$(K_X+B)^2=(K_X+B)\cdot \left(K_X+T+\frac{10}{11}S+\frac{1}{11}S\right)\geq\frac{1}{11}(K_X+B)\cdot S\geq\frac{1}{462}.$$
This implies (1). Moreover, if $(K_X+B)^2=\frac{1}{462}$, then $\deg(K_S+B_S)=\frac{1}{42}$ and $\pet(X,T;S)=\frac{10}{11}$. Thus $S$ is a non-singular rational curve. By Theorem \ref{thm: 12/13 pet}(3.a), $T=0$. This implies (2).
\end{proof}

\begin{lem}\label{lem: 237}
Let $(X,S)$ be a $\Qq$-factorial projective plt surface pair such that $S$ is prime divisor, $K_X+S$ is nef, $\pet(X,0;S)=\frac{12}{13}$, and $(K_X+S)\cdot S>0$. Then $$(K_X+S)^2>\frac{1}{462}.$$
\end{lem}
\begin{proof}
Suppose that $(K_X+S)^2\leq\frac{1}{462}$. Since $(X,S)$ is plt, $S$ is normal. We let $x_1,\dots,x_k$ be the singular points of $X$ along $S$ such that $x_i$ has order $p_i$ for each $i$ and let $K_S+B_S:=(K_X+S)|_S$. Then
\begin{align*}
    \frac{1}{462}&\geq(K_X+S)^2=(K_X+S)\cdot \left(K_X+\frac{12}{13}S+\frac{1}{13}S\right)\\
    &\geq\frac{1}{13}(K_X+S)\cdot S=\frac{1}{13}\deg(K_X+S)|_S=\deg(K_S+B_S)\\
    &=\frac{1}{13}\left(2p_a(S)-2+\sum_{i=1}^k\left(1-\frac{1}{p_i}\right)\right)>0.
\end{align*}
Thus $p_a(S)=0,k=3,$ and $\{p_1,p_2,p_3\}=\{2,3,7\}$. But this contradicts Theorem \ref{thm: 12/13 pet}(3.b).
\end{proof}

\begin{proof}[Proof of Theorem \ref{thm: vol 1/462}]
Let $S$ be a component of $\Supp B$, $f: Y\rightarrow X$ a $\Qq$-factorial dlt modification of $(X,B)$, $K_Y+B_Y:=f^*(K_X+B)$, and $S_Y:=f^{-1}_*S$. We only need to show that $\vol(Y,B_Y):=(K_Y+B_Y)^2\geq\frac{1}{462}$. We have
$$(K_Y+B_Y)\cdot S_Y=f^*(K_X+B)\cdot S_Y=(K_X+B)\cdot S>0.$$
By Lemma \ref{lem: pet small case}, we may assume that $\pet(Y,B_Y-S_Y;S_Y)=\frac{12}{13}$. By Lemma \ref{thm: 12/13 pet}(3.a), $B_Y=S_Y$. The theorem now follows from Lemma \ref{lem: 237}.
\end{proof}

\begin{rem}\label{rem: connection with stable degeneration}
We briefly recall the connection between the lower bound of the volume $(K_X+S)^2$ in Theorem \ref{thm: vol 1/462} and the degeneration of stable surfaces of general type.

Let $f: \mathcal{X}\rightarrow C$ be a stable degeneration of surfaces of general type, $X_g$ a general fiber of $f$, and $X_0$ a special fiber of $f$. One natural thing in the study of the structure of the special fiber $X_0$ is the number of irreducible components of $X_0$. We let $\tilde X_0\rightarrow X_0$ be the normalization of $X_0$, and let $X_1,\dots,X_m$ be the irreducible components of $\tilde X_0$. Let $S_i\subset X_i$ be the preimage of the double locus of $X_i$. Then each $S_i$ is a Weil divisor, and we have
$$K_{X_g}^2=K_{X_0}^2=\sum_{i=1}^m(K_{X_i}+S_i)^2.$$

If $S_i=0$ for some $i$, then the image of $X_i$ in $X_0$ is normal. By the connectedness of $X_0$, $X_0$ is normal. Thus $m=1$. Therefore, we may assume that $S_i\not=0$ for any $i$. In this case, by Theorem \ref{thm: vol 1/462}, $(K_{X_i}+S_i)^2\geq\frac{1}{462}$ for every $i$. In particular, $m\leq 462K_{X_g}^2$. In particular, this tells us that, for any stable degeneration of surfaces of general type, if the volume of a special fiber is $<\frac{1}{231}$, then any special fiber has at most $1$ component.

We can also improve our estimation by considering the following:
\begin{enumerate}
    \item If $S_i$ is not a prime divisor, then by the same arguments as in the proof of Theorem \ref{thm: vol 1/462}, it is easy to show that $(K_{X_i}+S_i)^2\geq2\cdot\frac{1}{42}\cdot\frac{1}{13}=\frac{1}{273}$.
    \item If $S_i$ is a prime divisor, the either there exists an involution $X_i\rightarrow X_i$, or there exists a unique index $j\not=i$ such that $S_i$ intersects $S_j$. Now by the connectedness of $X_0$, either $m=2$, or $S_j$ intersects $S_k$ for some $k\not=i,j$, or there exists an involution $X_j\rightarrow X_j$. 
\end{enumerate}
With these observations, it will be interesting to find the optimal lower bound of the minimal volume of $K_{X_i}+S_i$ which admits an involution, and the minimal volume of $K_{X_i}+S_i$ such that $S_i$ is not a prime divisor. A better estimation of these volumes will provide us a better bound of $m$.
\end{rem}

\subsection{Proof of Theorem \ref{thm: summary gap conjecture}}

\begin{proof}[Proof of Theorem \ref{thm: summary gap conjecture}]
(1) Let $(X,bS)$ be a pair such that $\dim X=1$, $b\in (0,1)$, and $S$ a non-zero effective Weil divisor on $X$. If $(X,bS)$ is $\Rr$-complementary or $K_X+bS\equiv 0$, then $X=\mathbb P^1$. Thus if $b=\Rct(X,0;S)$, then $K_X+bS\equiv 0$, and we have $b=\frac{2}{\deg S}\leq\frac{2}{3}$. Thus $\delta_{\glct}(1)=\delta_{\Rct}(1)=\delta_{\Rct,\Ft}(1)=\frac{1}{3}$. By \cite[Lemma 5.1]{Jia21}, $\delta_{\tmld}(2)=\delta_{\mld}(2)=\frac{1}{3}$. 

(2) By Theorem \ref{thm: global lct 12/13 intro}, Theorem \ref{thm: 12/13 rct intro}, Corollary \ref{cor: calabi-yau surface mld} (=\cite[Proposition 6.1]{ETW22}), and \cite[Theorem 1.4]{LX21a}, $\delta_{\glct}(2)=\delta_{\tmld}(3)=\delta_{\Rct}(2)=\delta_{\mld,\CY}(2)=\frac{1}{13}$. By definition, $\delta_{\Rct,\Ft}(2)\geq\delta_{\Rct}(2)=\frac{1}{13}$. By considering Example \ref{ex: 12/13}, we know $\delta_{\Rct,\Ft}(2)\leq\frac{1}{13}$, so $\delta_{\Rct,\Ft}(2)=\frac{1}{13}$. By definition, $\frac{1}{13}=\delta_{\tmld}(3)\geq\delta_{\mld}(3)$.

(3) $\delta_{\tmld}(d+1)\geq\delta_{\mld}(d+1)$ and $\delta_{\Rct,\Ft}(d)\geq\delta_{\Rct}(d)$ follow from the definitions. By the global ACC \cite[Theorem 1.5]{HMX14}, $\delta_{\mld,\CY}(d)>0$ and $\delta_{\glct}(d)>0$. 

By (1), we may assume that $d\geq 2$. Then $\delta_{\mld,\CY}(d)<1$, and there exists a klt Calabi-Yau variety $X$ such that $\mld(X)=\delta_{\mld,\CY}(d)<1$. We let $f: Y\rightarrow X$ be an extraction which extracts a prime divisor $E$ such that $a(E,X,0)=\mld(X)$. Then $(Y,(1-a)E)$ is klt of dimension $d$ and $K_Y+(1-\mld(X))E\equiv 0$. Thus $0<\delta_{\glct}(d)\leq \mld(X)=\delta_{\mld,\CY}(d)$.

(4) We prove (4) by induction on dimension. When $d=1$ it follows from (1). Now we may suppose that $d\geq 2$.

We let $(X,bS)$ be a klt pair of dimension $d$ such that $S$ is a non-zero effective Weil divisor, $K_X+bS\equiv 0$, and $b=1-\delta_{\glct}(d)$. Possibly replacing $(X,bS)$ with a small $\Qq$-factorialization, we may assume that $X$ is $\Qq$-factorial. We run a $K_X$-MMP which terminates with a Mori fiber space $X'\rightarrow Z$. Let $S'$ be the image of $S$ on $X'$, $F$ a general fiber of $X'\rightarrow Z$, and $S_F:=S'|_F$. Then $S_F\not=0$, $K_F+bS_F\equiv 0$, and $(F,bS_F)$ is klt.

If $\dim F<\dim X$, then by induction on dimension, we have $\delta_{\mld}(\dim F+1)\leq\delta_{\glct}(\dim F)\leq 1-b$ and $\delta_{\Rct,\Ft}(\dim F)\leq\delta_{\glct}(\dim F)\leq 1-b$. Therefore:
\begin{itemize}
    \item There exists a germ $W\ni w$ of dimension $\dim F+1$ such that $b\leq\mld(W;w)<1$. Let $W':=W\times\mathbb C^{d-\dim F}$ and $\eta$ the generic point of $\bar w\times\mathbb{C}^{d-\dim F}$ in $W'$. Now
$$b\leq\mld(W';\eta)=\mld(W;w)<1,$$
so $1-\mld(W'\ni\eta)\leq\delta_{\mld}(d+1)$, and so $\delta_{\mld}(d+1)\leq 1-b=\delta_{\glct}(d)$.
\item By the ACC for $\Rr$-complementary thresholds for Fano type varieties (\cite[Theorem 8.20]{HLS19}, \cite[Theorem 21]{Sho20}), there exists a Fano type variety $V$ and a non-zero effective $\Qq$-Cartier Weil divisor $S_V$ on $V$ such that $\Rct(V,0;S_V)=1-\delta_{\Rct,\Ft}(\dim F)\geq b$. Now $V':=V\times\mathbb P^{d-\dim F}$ is a Fano type variety and $\Rct(V',0;S_{V'}:=S_V\times\mathbb P^{d-\dim F})\geq b$, so $$1-\delta_{\Rct,\Ft}(d)\geq\Rct(V',0;S_{V'})\geq b=1-\delta_{\glct}(d).$$
Thus $\delta_{\Rct,\Ft}(d)\leq\delta_{\glct}(d)$.
\end{itemize}

If $\dim F=\dim X$, then $\rho(X')=1$. In particular, $X'$ is Fano and 
$$1-\delta_{\Rct,\Ft}(d)\geq\Rct(X',0;S')=b=1-\delta_{\glct}(d),$$
so $\delta_{\Rct,\Ft}(d)\leq\delta_{\glct}(d)$. Moreover, $S'$ is an ample Weil divisor. We consider the cone $V:=C(X',S')$ (note that although $S'$ is not necessarily Cartier, such cone construction is still valid, cf. \cite[2.3]{Ber21}). Let $S_V$ be the cone of $S'$ in $V$. We let $v$ be the vertex of $V$ and $f: W\rightarrow V$ be the blow-up of $v$ with exceptional divisor $E$. Then $E\cong X'$, and we have $(K_W+E)|_E=K_E\cong K_{X'}\sim_{\mathbb Q}-bS'$
and $-E|_E\sim_{\mathbb Q}S'$. Thus $(K_W+(1-b)E)|_E\sim_{\mathbb Q}0$, so $$\mld(V;v)=a(E,V,0)=b<1,$$
hence $$\delta_{\mld}(d+1)\leq 1-\mld(V;v)=1-b=\delta_{\glct}(d)$$
and we are done.
\end{proof}

\section{Exceptional del Pezzo surfaces with small mld's}\label{sec: 1/11 exceptional fano}

In this section, we classify all exceptional del Pezzo surfaces that are not $\frac{1}{11}$-lt and prove Theorem \ref{thm: 1/13 exceptional Fano} and Corollary \ref{cor: new bound alpha invariant}.

\begin{setup}\label{setup: terminalization of example}
$X_0,B_0,S_0,h,W,B_W,S_W,E_{i,j}$ are as follows:
\begin{enumerate}
\item $(X_0,B_0=\frac{12}{13}S_0):=(\mathbb P(3,4,5),\frac{12}{13}(x^3y+y^2z+z^2x=0))$.
    \item $h: W\rightarrow X_0$ is the terminalization of $(X_0,B_0)$ (cf. Construction \ref{cons: terminalization}).
    \item $K_W+B_W:=h^*(K_{X_0}+B_0)$.
    \item $S_W$ is the strict transform of $S_0$ on $W$.
    \item $E_{1,1},\dots, E_{1,12},E_{2,1},\dots,E_{2,17},E_{3,1},\dots,E_{3,21}$ are the prime exceptional divisors of $h$, such that for any $i\in\{1,2,3\}$ and any $j,k$,
    \begin{itemize}
        \item $E_{i,1}$ intersects $S_W$, and
        \item $E_{i,j}$ intersects $E_{i,k}$ if and only if $|j-k|\leq 1$.
    \end{itemize}
\end{enumerate}
\end{setup}
\begin{lem}\label{lem: non-1/11-lt exceptional fano dominated by termialization}
Conditions as in Set-up \ref{setup: terminalization of example}. Let $X$ be an exceptional surface that is not $\frac{1}{11}$-lt, $S$ an exceptional$/X$ prime divisor such $a(S,X,0)=\mld(X)$, and $\phi: Y\rightarrow X$ the extraction of $S$. Then $h$ factors through $Y$, i.e. we have the following diagram of projective birational morphisms:
\begin{center}$\xymatrix{
W\ar@{->}[r]^h\ar@{->}[d] & X_0\\
Y\ar@{->}[r]^\phi\ar@{->}[ur]^{\psi} & X.
}$
\end{center}
Moreover, $S_0$ is the image of $S$ on $X_0$.
\end{lem}
\begin{proof}
Since $X$ is of Fano type and $a(S,X,0)=\mld(X)\leq\frac{1}{11}$, there exists a klt Calabi-Yau pair $(X,G)$ such that $a:=a(S,X,G)<\frac{1}{11}$. We let $G_Y:=\phi^{-1}_*G$. Then $(Y,(1-a)S+G_Y)$ is klt log Calabi-Yau, so $(Y,(1-a)S)$ is $\Rr$-complementary. If $(Y,S)$ is $\Rr$-complementary, then we let $(Y,S^+)$ be an $\Rr$-complement of $(Y,S)$ and let $\Delta:=\phi_*S^+$. Then $(X,\Delta)$ is log Calabi-Yau but not klt, which contradicts the assumption that $X$ is exceptional. Thus $\Rct(Y,0;S)\in [1-a,1)\subset (\frac{10}{11},1)$. By Theorem \ref{thm: 12/13 rct intro}, $\Rct(Y,0;S)=\frac{12}{13}$.

Let $P$ and $N$ be the positive and negative part of the Zariski decomposition of $-(K_Y+\frac{12}{13}S)$. Let $f: Y\rightarrow Y'$ be the contraction of $N$, and let $S':=f_*S$. Then $(Y',\frac{12}{13}S')$ is a maximal model (cf. \cite[Construction 2]{Sho20}) of $(Y,\frac{12}{13}S)$. By \cite[Proposition 8, Addendum 8]{Sho20}, there exists a one-to-one correspondence of $\mathbb R$-complements of $(Y,\frac{12}{13}S)$ and $(Y',\frac{12}{13}S')$. Since $(Y,\frac{12}{13}S)$ is $\Rr$-complementary, $(Y',\frac{12}{13}S')$ is $\Rr$-complementary. In particular, $(Y',\frac{12}{13}S')$ is lc. By Lemma \ref{lem: surface lct 5/6}, $(Y',S')$ is lc. By \cite[Lemma 14]{Sho20}, $K_{Y'}+\frac{12}{13}S'\equiv 0$. By Theorem \ref{thm: global lct 12/13}(3.b), the lemma immediately follows.
\end{proof}

\begin{setup}\label{setup: classify non-1/11-lt}
$X,X_0,B_0,S_0,h,W,B_W,S_W,E_{i,j},f,g,\bar X,\pi,B,\bar B,\bar S,\bar E_{i,j,X},E_{i,j,X}$ are as follows.
\begin{enumerate}
\item $X$ is an exceptional Fano type surface that is not $\frac{1}{11}$-lt.
    \item $X_0,B_0,S_0,h,W,B_W,S_W,E_{i,j}$ are as in Set-up \ref{setup: terminalization of example}.
    \item $f: W\rightarrow X$ is a projective birational morphism, whose existence is guaranteed by Lemma \ref{lem: non-1/11-lt exceptional fano dominated by termialization}.
    \item $g: \bar X\rightarrow X$ is the minimal resolution of $X$.
    \item $\bar f: W\rightarrow \bar X$ is the induced birational map, which is a projective birational morphism as $W$ is non-singular.
    \item $B$ is the image of $B_W$ on $X$, $\bar B$ is the image of $B_W$ on $\bar X$, and $\bar S$ is the image of $S_W$ on $\bar X$.
    \item $E_{i,j,\bar X}$ is the image of $E_{i,j}$ on $\bar X$ and $E_{i,j,X}$ is the image of $E_{i,j}$ on $X$ for any $i,j$.
\end{enumerate}
\end{setup}

\begin{lem}\label{lem: 1/11-lt contract divisors}
Conditions as in Set-up \ref{setup: classify non-1/11-lt}.
\begin{enumerate}
    \item $\bar S\not=0$ and $\bar S$ is contracted by $g$. In particular, $S_W$ is not contracted by $\bar f$ and is contracted by $f$.
    \item Any divisor contracted by $f$ is either $S_W$ or $E_{i,j}$.
    \item If $X$ is a del Pezzo surface, then for any $i,j$ such that $E_{i,j,X}\not=0$, 
    \begin{enumerate}
    \item $K_{\bar X}\cdot E_{i,j,\bar X}\leq K_X\cdot E_{i,j,X}<0$, and
    \item  $E_{i,j,\bar X}^2=-1$.
    \end{enumerate}
\end{enumerate}
\end{lem}
\begin{proof}
(1) Let $Y\rightarrow X$ be the extraction of a prime divisor $S$ such that $a(S,X,0)=\mld(X)$. Since $\mld(X)\leq\frac{1}{11}<1$, the center of $S$ on $\bar X$ is a divisor. In particular, $g$ factors through $Y$. By Lemma \ref{lem: non-1/11-lt exceptional fano dominated by termialization}, $h$ factors through $Y$ and $S_0$ is the image of $S$ on $X_0$.  Thus $\bar S$ is the strict transform of $S$ on $\bar X$, and (1) follows.

(2) Let $F$ be a prime divisor contracted by $f$ such that $F\not=S_W$. Then $F$ is contracted by the morphism $W\rightarrow Y$. Hence $F$ is contracted by the morphism $W\rightarrow X_0$. Since $E_{i,j}$ are the only divisors contracted by $W\rightarrow X_0$, we get (2).

(3) Since $X$ is a del Pezzo surface, we get (3.a). Since $g$ only contracts divisors with self-intersection $\leq -2$, by (3.a), we have 
$0>K_X\cdot E_{i,j,X}\geq K_{\bar X}\cdot E_{i,j,\bar X}$. Since $\bar X$ is non-singular, $K_{\bar X}\cdot E_{i,j,\bar X}\leq -1$. Since $E_{i,j,\bar X}$ is contracted by $\bar X\rightarrow X_0$, $E_{i,j,\bar X}$ is an exceptional$/X_0$ non-singular rational curve and $E_{i,j,\bar X}^2\leq -1$. This implies (3.b).
\end{proof}

By Lemma \ref{lem: 1/11-lt contract divisors}, there are only finitely many possibilities of $X$ and all of them can be precisely described. At this point, we can run a computer program and check the structure of these varieties by enumeration to get our desired classification. Nevertheless, we still want to provide a proof which can be computed by hand. We have the following lemma.

\begin{lem}\label{lem: 1/11-lt exceptional classification step 1}
Conditions as in Set-up \ref{setup: classify non-1/11-lt}. Assume that $X$ is a del Pezzo surface. Then for each $i\in\{1,2,3\}$, we have the following lists of possibilities of the set of $j$ such that $E_{i,j,X}\not=0$. In each of the following tables for each fixed $i\in\{1,2,3\}$, we have the following.
\begin{enumerate}
\item ``Divisors on $X$" contains all indices $j$ such that $E_{i,j,X}\not=0$.  If no such index exists, we fill in ``No Divisor". In addition, we let $j_i:=\min\{j,22\mid E_{i,j,X}\not=0\}$.
\item ``Change near $S_W$" contains $(u;\bm{v})$ satisfying the following. 
    \begin{enumerate}
        \item We let $W\rightarrow X_i$ be the contraction of all $E_{i,j}$ where $j$ appears in ``$\overline{f}$-exceptional divisors", and $S_i$ the image of $S_W$ on $X_i$. We define $u:=S_i^2-S_W^2$. 
        \item There exists a unique non-negative integer $k_i$ such that       $$\cup_{j<j_i,E_{i,j,\bar X}\not=0}E_{i,j,\bar X}=E_{i,l_1,\bar X}\cup E_{i,l_2,\bar X}\cup\dots\cup E_{i,l_{k_i},\bar X},$$
        where $l_1\leq l_2\leq\dots\leq l_{k_i}$. We let $\bm{v}:=(-E_{i,l_1,\bar X}^2,\dots,-E_{i,l_{k_i},\bar X}^2)$ if $k_i\geq 1$ and $\bm{v}:=0$ if $k_i=0$.
    \end{enumerate}
    \item We let $E_{i,j,X_i}$ be the image of $E_{i,j}$ on $X_i$ for each $i$. ``Left intersection" contains $K_{X_i}\cdot E_{i,j_i,X_i}$ if $j_i\leq 21$, and contains ``No Intersection" if $j_i=22$.
\end{enumerate}
\begin{center}
    \begin{longtable}{|l|l|l|l|l|}
 			\caption{Possible contractions of $E_{1,j}$}\label{table: contraction E1}\\
		\hline
		Case & Divisors on $X$ & Change near $S_W$ & Left intersection\\
				\hline
1 & No Divisor &  $(8;2,2)$ & No Intersection \\
		\hline
2 & $1$ & $(0;0)$ & $-\frac{12}{13}$ \\
\hline
3 & $2$  & $(1;0)$ & $-\frac{21}{23}$ \\
\hline
4 & $2,8$  & $(1;0)$ & $-\frac{10}{11}$ \\
\hline
5 & $3$  &  $(2;0)$ & $-\frac{9}{10}$ \\
\hline
6 & $3,8$  & $(2;0)$ & $-\frac{8}{9}$ \\
\hline
7 & $4$  &  $(3;0)$ & $-\frac{15}{17}$ \\
\hline
8 & $4,8$  &  $(3;0)$ & $-\frac{6}{7}$ \\
\hline
9 & $5$  &  $(4;0)$ & $-\frac{6}{7}$ \\
\hline
10 & $5,8$  &  $(4;0)$ & $-\frac{4}{5}$ \\
\hline
11 & $6$  & $(5;0)$ & $-\frac{9}{11}$ \\
\hline
12 & $6,8$  & $(5;0)$ & $-\frac{2}{3}$ \\
\hline
13 & $6,10$  & $(5;0)$ & $-\frac{4}{5}$ \\
\hline
14 & $7$  & $(6;0)$ & $-\frac{3}{4}$ \\
\hline
15 & $7,10$  & $(6;0)$ & $-\frac{2}{3}$ \\
\hline
16 & $8$  &  $(6;2)$ & $-\frac{6}{13}$ \\
\hline
16 & $9$  &  $(7;0)$ & $-\frac{3}{5}$ \\
\hline
17 & $10$  & $(7;2)$ & $-\frac{3}{7}$ \\
\hline
\end{longtable}
\end{center}
\begin{center}
    \begin{longtable}{|l|l|l|l|l|}
 			\caption{Possible contractions of $E_{2,j}$}\label{table: contraction E2}\\
		\hline
		Case & Divisors on $X$ & Change near $S_W$ & Left intersection\\
				\hline
1 & No Divisor  & $(9;2,2,2)$ & No Intersection \\
		\hline
2 & $1$ & $(0;0)$ & $-\frac{12}{13}$ \\
\hline
3 & $2$  & $(1;0)$ & $-\frac{32}{35}$ \\
\hline
4 & $2,8$  & $(1;0)$ & $-\frac{10}{11}$ \\
\hline
5 & $2,13$  & $(1;0)$ & $-\frac{21}{23}$ \\
\hline
6 & $3$  & $(2;0)$ & $-\frac{28}{31}$ \\
\hline
7 & $3,8$  &  $(2;0)$ & $-\frac{8}{9}$ \\
\hline
8 & $3,13$  & $(2;0)$ & $-\frac{9}{10}$ \\
\hline
9 & $4$  & $(3;0)$ & $-\frac{8}{9}$ \\
\hline
10 & $4,8$  & $(3;0)$ & $-\frac{6}{7}$ \\
\hline
11 & $4,13$  & $(3;0)$ & $-\frac{15}{17}$ \\
\hline
12 & $5$  & $(4;0)$ & $-\frac{20}{23}$ \\
\hline
13 & $5,8$  &  $(4;0)$ & $-\frac{4}{5}$ \\
\hline
14 & $5,10$  & $(4;0)$ & $-\frac{6}{7}$ \\
\hline
15 & $5,13$  &  $(4;0)$ & $-\frac{6}{7}$ \\
\hline
16 & $6$  &  $(5;0)$ & $-\frac{16}{19}$ \\
\hline
17 & $6,8$  & $(5;0)$ & $-\frac{2}{3}$ \\
\hline
18 & $6,10$  & $(5;0)$ & $-\frac{4}{5}$ \\
\hline
19 & $6,10,13$  & $(5;0)$ & $-\frac{4}{5}$ \\
\hline
20 & $6,13$  & $(5;0)$ & $-\frac{9}{11}$ \\
\hline
21 & $7$  & $(6;0)$ & $-\frac{4}{5}$ \\
\hline
22 & $7,10$  & $(6;0)$ & $-\frac{2}{3}$ \\
\hline
23 & $7,10,13$  & $(6;0)$ & $-\frac{2}{3}$ \\
\hline
24 & $7,13$  & $(6;0)$ & $-\frac{3}{4}$ \\
\hline
25 & $8$  & $(6;2)$ & $-\frac{6}{13}$ \\
\hline
26 & $9$  & $(7;0)$ & $-\frac{8}{11}$ \\
\hline
27 & $9,12$  & $(7;0)$ & $-\frac{2}{3}$ \\
\hline
28 & $9,13$  & $(7;0)$ & $-\frac{3}{5}$ \\
\hline
29 & $9,15$  & $(7;0)$ & $-\frac{5}{7}$ \\
\hline
30 & $10$  & $(7;2)$ & $-\frac{4}{9}$ \\
\hline
31 & $10,13$  &  $(7;2)$ & $-\frac{3}{7}$ \\
\hline
32 & $11$  & $(8;0)$ & $-\frac{4}{7}$ \\
\hline
33 & $11,15$  & $(8;0)$ & $-\frac{1}{2}$ \\
\hline
34 & $12$  & $(8;2)$ & $-\frac{2}{5}$ \\
\hline
35 & $13$  & $(8;2,2)$ & $-\frac{4}{13}$ \\
\hline
36 & $15$  & $(9;3)$ & $-\frac{3}{5}$ \\
\hline
\end{longtable}
\end{center}
\begin{center}
\begin{longtable}{|l|l|l|l|l|}
 			\caption{Possible contractions of $E_{3,j}$}\label{table: contraction E3}\\
		\hline
		Case & Divisors on $X$ & Change near $S_W$ & Left intersection\\
				\hline
1 & No Divisor & $(10;3,2)$ & No Intersection \\
\hline
2 & $1$ & $(0;0)$ & $-\frac{12}{13}$ \\
\hline
3 & $2$ & $(1;0)$ & $-\frac{43}{47}$ \\
\hline
4		& $2,8$ & $(1;0)$ & $-\frac{10}{11}$ \\
\hline
5		& $2,13$ & $(1;0)$ & $-\frac{21}{23}$ \\
\hline
6 & $2,18$ & $(1;0)$ & $-\frac{32}{35}$ \\
\hline
7 & $3$ & $(2;0)$ & $-\frac{19}{21}$ \\
\hline
8 & $3,8$ & $(2;0)$ & $-\frac{8}{9}$ \\
\hline
9 & $3,13$ & $(2;0)$ & $-\frac{9}{10}$ \\
\hline
10& $3,18$ & $(2;0)$ & $-\frac{28}{31}$ \\
\hline
11& $4$ & $(3;0)$ & $-\frac{33}{37}$ \\
\hline
12& $4,8$ & $(3;0)$ & $-\frac{6}{7}$ \\
\hline
13& $4,10$ & $(3;0)$ & $-\frac{8}{9}$ \\
\hline
14& $4,13$ & $(3;0)$ & $-\frac{15}{17}$ \\
\hline
15& $4,18$ & $(3;0)$ & $-\frac{8}{9}$ \\
\hline
16& $5$ & $(4;0)$ & $-\frac{7}{8}$ \\
\hline
17& $5,8$ & $(4;0)$ & $-\frac{4}{5}$ \\
\hline
18& $5,10$ & $(4;0)$ & $-\frac{6}{7}$ \\
\hline
19& $5,10,18$ & $(4;0)$ & $-\frac{6}{7}$ \\
\hline
20& $5,13$ & $(4;0)$ & $-\frac{6}{7}$ \\
\hline
21& $5,18$ & $(4;0)$ & $-\frac{20}{23}$ \\
\hline
22& $6$ & $(5;0)$ & $-\frac{23}{27}$ \\
\hline
23& $6,8$ & $(5;0)$ & $-\frac{2}{3}$ \\
\hline
24& $6,10$ & $(5;0)$ & $-\frac{4}{5}$ \\
\hline
25& $6,10,13$ & $(5;0)$ & $-\frac{4}{5}$ \\
\hline
26& $6,10,18$ & $(5;0)$ & $-\frac{4}{5}$ \\
\hline
27& $6,13$ & $(5;0)$ & $-\frac{9}{11}$ \\
\hline
28& $6,15$ & $(5;0)$ & $-\frac{11}{13}$ \\
\hline
29& $6,18$ & $(5;0)$ & $-\frac{16}{19}$ \\
\hline
30& $7$ & $(6;0)$ & $-\frac{9}{11}$ \\
\hline
31& $7,10$ & $(6;0)$ & $-\frac{2}{3}$ \\
\hline
32& $7,10,13$ & $(6;0)$ & $-\frac{2}{3}$ \\
\hline
33& $7,10,18$ & $(6;0)$ & $-\frac{2}{3}$ \\
\hline
34& $7,12$ & $(6;0)$ & $-\frac{4}{5}$ \\
\hline
35& $7,13$ & $(6;0)$ & $-\frac{3}{4}$ \\
\hline
36& $7,15$ & $(6;0)$ & $-\frac{4}{5}$ \\
\hline
37& $7,18$ & $(6;0)$ & $-\frac{4}{5}$ \\
\hline
38& $8$ & $(6;2)$ & $-\frac{6}{13}$ \\
\hline
39& $9$ & $(7;0)$ & $-\frac{13}{17}$ \\
\hline
40& $9,12$ & $(7;0)$ & $-\frac{2}{3}$ \\
\hline
41& $9,12,15$ & $(7;0)$ & $-\frac{2}{3}$ \\
\hline
42& $9,12,18$ & $(7;0)$ & $-\frac{2}{3}$ \\
\hline
43& $9,13$ & $(7;0)$ & $-\frac{3}{5}$ \\
\hline
44& $9,15$ & $(7;0)$ & $-\frac{5}{7}$ \\
\hline
45& $9,15,18$ &  $(7;0)$ & $-\frac{5}{7}$ \\
\hline
46& $9,17$ & $(7;0)$ & $-\frac{3}{4}$ \\
\hline
47& $9,18$ & $(7;0)$ & $-\frac{8}{11}$ \\
\hline
48& $10$ & $(7;2)$ & $-\frac{13}{29}$ \\
\hline
49& $10,13$ & $(7;2)$ & $-\frac{3}{7}$ \\
\hline
50& $10,18$ & $(7;2)$ & $-\frac{4}{9}$ \\
\hline
51& $11$ & $(8;0)$ & $-\frac{2}{3}$ \\
\hline
52& $11,15$ & $(8;0)$ & $-\frac{1}{2}$ \\
\hline
53& $11,15,18$ &  $(8;0)$ & $-\frac{1}{2}$ \\
\hline
54& $11,17$ & $(8;0)$ & $-\frac{3}{5}$ \\
\hline
55& $11,18$ & $(8;0)$ & $-\frac{4}{7}$ \\
\hline
56& $12$ & $(8;2)$ & $-\frac{8}{19}$ \\
\hline
57& $12,15$ & $(8;2)$ & $-\frac{2}{5}$ \\
\hline
58& $12,18$ & $(8;2)$ & $-\frac{2}{5}$ \\
\hline
59& $13$ & $(8;2,2)$ & $-\frac{4}{13}$ \\
\hline
60& $14$ & $(9,0)$ & $-\frac{3}{7}$ \\
\hline
61& $14,20$ &  $(9;0)$ & $-\frac{2}{5}$ \\
\hline
62& $15$ & $(9;3)$ & $-\frac{5}{8}$ \\
\hline
63& $15,18$ &  $(9;3)$ & $-\frac{3}{5}$ \\
\hline
64& $16$ & $(9;2)$ & $-\frac{1}{3}$ \\
\hline
65& $17$ & $(9;2,2)$ & $-\frac{3}{11}$ \\
\hline
66& $18$ & $(9;2,2,2)$ & $-\frac{3}{13}$ \\
\hline
67 & $20$ & $(10;4)$ & $-\frac{2}{3}$ \\
\hline
\end{longtable}
\end{center}
\end{lem}
\begin{proof}
By Lemma \ref{lem: 1/11-lt contract divisors}(2), only $E_{i,j}$ and $S_W$ can be contracted by $f$, and the rest of the proof is enumeration and applying Lemma \ref{lem: 1/11-lt contract divisors}(3.a) to $E_{i,j,X}$ such that $j\not=j_i$.

An important observation is that Lemma \ref{lem: 1/11-lt contract divisors}(3.b) guarantees that, except the curve $E_{i,j_i}$, only very few curves can be not contracted by $f$. The only possible indices of $j$ such that $j\not=j_i$ and $E_{i,j,X}\not=0$ are $8,10$ (for $i\in\{1,2,3\}$), $13,15$ (for $i\in\{2,3\}$), and $17,18,20$ (for $i=3$). This greatly simplifies the enumeration, and makes the enumeration computable by hand.
\end{proof}

\begin{lem}\label{lem: nonexc 12346}
Let $(X,B)$ be an $\Rr$-complementary surface pair such that the coefficients of $B$ belong to $\{1-\frac{1}{m}\mid m\in\mathbb N^+\}\cup [\frac{6}{7},1]$. If $(X,B)$ is not exceptional, then $(X,B)$ has a monotonic $n$-complement for some $n\in\{1,2,3,4,6\}$.
\end{lem}
\begin{proof}
Let $(X,B+G)$ be an $\Rr$-complement of $(X,B+G)$ such that $(X,B+G)$ is not klt. Let $f: Y\rightarrow X$ be a dlt modification of $(X,B+G)$, $E$ the reduced $f$-exceptional divisor, and $B_Y:=f^{-1}_*B+E$. Then $(Y,B_Y)$ is $\Rr$-complementary and the coefficients of $B_Y$ belong to $\{1-\frac{1}{m}\mid m\in\mathbb N^+\}\cup [\frac{6}{7},1]$. 

Let $P$ and $N$ be the positive and negative part of the Zariski decomposition of $-(K_Y+B_Y)$. Let $h: Y\rightarrow Y'$ be the contraction of $N$, and let $B_{Y'}:=h_*B_Y$. Then $(Y',B_{Y'})$ is a maximal model (cf. \cite[Construction 2]{Sho20}) of $(Y,B_Y)$ and the coefficients of $B_{Y'}$ belong to  $\{1-\frac{1}{m}\mid m\in\mathbb N^+\}\cup [\frac{6}{7},1]$. By \cite[2.3 Inductive Theorem]{Sho00} (applying conditions (NK), (NEF), (M), and (EEC) of \cite{Sho00}),  $(Y',B_{Y'})$ has an $n$-complement $(Y',B_{Y'}^+)$ for some $n\in\{1,2,3,4,6\}$. Since the coefficients of $B_{Y'}$ belong to $\{1-\frac{1}{m}\mid m\in\mathbb N^+\}\cup [\frac{6}{7},1]$, by \cite[2.7. Monotonicity Lemma]{Sho00}, $(Y',B_{Y'}^+)$ is a monotonic $n$-complement of $(Y',B_{Y'})$. Let $K_{Y}+B_Y^+:=h^*(K_{Y'}+B_{Y'}^+)$. By \cite[Proposition 8, Addendum 8]{Sho20}, $(Y,B_Y^+)$ is a monotonic $n$-complement of $(Y,B_Y)$. Thus $(X,B^+:=f_*B_Y^+)$ is a monotonic $n$-complement of $(X,B)$. 
\end{proof}

\begin{exlem}\label{exlem: special exc fano}
Conditions as in Set-up \ref{setup: classify non-1/11-lt}. Suppose that there exists a unique pair of integers $(i,j)$ such that $E_{i,j,X}\not=0$, and $j=2$. Then $X$ is one of three exceptional del Pezzo surfaces of Picard number $1$, precisely described in Table \ref{table: picard number 1 fano}. Here
\begin{itemize}
    \item ``Divisors on $X$" consists of all $(i,j)$ such that $E_{i,j,X}\not=0$, and
    \item ``Worst Sing" consists of the type of the singularity of the closed point $\Center_XS_W$.
\end{itemize}
\begin{center}
    \begin{longtable}{|l|l|l|l|l|l|}
 			\caption{Three exceptional del Pezzo surfaces of Picard number 1}\label{table: picard number 1 fano}\\
		\hline
		Case & Divisors on $X$ & Worst Sing & $\mld(X)$ & $K_X^2$ \\
				\hline
1 & $(1,2)$  &  $\frac{1}{97}(1,74)$ & $\frac{9}{97}$  & $\frac{1}{2231}$ \\
\hline
2 & $(2,2)$ &  $\frac{1}{89}(1,61)$ & $\frac{8}{89}$ & $\frac{1}{3115}$ \\
\hline
3 & $(3,2)$ &  $\frac{1}{79}(1,54)$ & $\frac{7}{79}$ & $\frac{1}{3713}$ \\
\hline
\end{longtable}
\end{center}
All these three surfaces are $13$-complementary: $(X,h_*B_W)$ is a $13$-complement of $X$.
\end{exlem}
\begin{proof}
By our construction, $X$ is a Fano type surface of Picard number $1$. So $X$ is a del Pezzo surface. We are only left to show that $X$ is exceptional. We let $x$ be the center of $S_W$ on $X$, $\phi: Y\rightarrow X$ the extraction of $S_W$, and $E_Y$ the image of $E_{i,2}$ on $Y$. Let $S_Y$ be the image of $S_W$ on $Y$. By Lemma \ref{lem: non-1/11-lt exceptional fano dominated by termialization}, there exists a contraction $\psi: Y\rightarrow X_0$ of $S_Y$. We let $S:=\psi_*S_Y$.

All algebraic invariants in Table \ref{table: picard number 1 fano} can be computed directly by hand. In particular, when $i=1$ (resp. $2$,$3$), we have $a(S,X,0)=\mld(X;x)=\mld(X)=\frac{9}{97}$ (resp. $\frac{8}{89}$, $\frac{7}{79}$).

Suppose that $X$ is not exceptional. By Lemma \ref{lem: nonexc 12346}, $X$ has an $n$-complement $(X,G)$ for some $n\in\{1,2,3,4,6\}$. Since $a(S,X,0)<\frac{1}{6}$, $a(S,X,G)=0$. This is not possible. Therefore, $X$ is exceptional.
\end{proof}

\begin{thm}[Classification of non-$\frac{1}{11}$-lt exceptional del Pezzo surfaces]\label{thm: 1/11-lt exceptional classification step 2}
Conditions as in Set-up \ref{setup: classify non-1/11-lt}. Assume that $X$ is a del Pezzo surface. Then $X$ is one of surfaces in Table \ref{table: 1/11-lt fano}. Here
\begin{itemize}
    \item ``Divisors on $X$“ consists of all $(i,j)$ such that $E_{i,j,X}\not=0$, and
    \item ``Worst Sing" consists of the type of the singularity of the closed point $\Center_XS_W$.
\end{itemize}
Table \ref{table: 1/11-lt fano} gives the complete list of non-$\frac{1}{11}$-lt exceptional del Pezzo surfaces. In other words,
\begin{enumerate}
\item any non-$\frac{1}{11}$-lt exceptional del Pezzo surface belongs to Table \ref{table: 1/11-lt fano}, and
\item for any $X$ in Table \ref{table: 1/11-lt fano}, $X$ is a non-$\frac{1}{11}$-lt exceptional del Pezzo surface.
\end{enumerate}
Moreover, 
\begin{enumerate}
    \item[(3)] any $X$ in Table \ref{table: 1/11-lt fano} has a $13$-complement.
\end{enumerate}

\begin{center}
    \begin{longtable}{|l|l|l|l|l|l|l|}
 			\caption{List of non-$\frac{1}{11}$-lt exceptional del Pezzo surfaces}\label{table: 1/11-lt fano}\\
		\hline
		Case & Divisors on $X$ & Worst Sing & $\mld(X)$ & $K_X^2$ & $\rho(X)$\\
				\hline
1 & $(2,2),(3,2),(3,18)$  &  $\frac{1}{46}(1,3)$ & $\frac{2}{23}$  & $\frac{2}{10465}$ & $3$ \\
\hline
2 & $(2,2),(3,2)$ &  $\frac{1}{46}(1,3)$ & $\frac{2}{23}$ & $\frac{9}{37835}$ & $2$ \\
		\hline
3 & $(1,2),(2,2),(2,13)$ &  $\frac{1}{68}(1,5)$ & $\frac{3}{34}$ & $\frac{1}{5083}$ & $3$\\
\hline
4 & $(1,2),(2,2)$ &  $\frac{1}{68}(1,5)$ & $\frac{3}{34}$ & $\frac{4}{13685}$  & $2$\\
		\hline	
5 & $(1,2),(3,2),(3,13)$ &  $\frac{1}{57}(1,4)$ & $\frac{5}{57}$ & $\frac{2}{17043}$ & $3$\\
\hline
6 & $(1,2),(3,2),(3,18)$ &  $\frac{1}{57}(1,4)$ & $\frac{5}{57}$ & $\frac{127}{596505}$ & $3$\\
\hline
7 & $(1,2),(3,2)$ &  $\frac{1}{57}(1,4)$ & $\frac{5}{57}$ & $\frac{16}{61617}$  & $2$\\
\hline
8 & $(1,2),(2,2),(3,18)$ &  $\frac{1}{57}(1,4)$ & $\frac{5}{57}$ & $\frac{127}{596505}$ & $3$\\
\hline
9 & $(1,2),(2,2),(2,13),(3,18)$ &  $\frac{1}{57}(1,4)$ & $\frac{5}{57}$ & $\frac{2}{17043}$ & $3$\\
\hline
10 & $(1,8),(2,2)$ &  $\frac{1}{81}(1,43)$ & $\frac{7}{81}$ & $\frac{2}{36855}$  & $2$\\
\hline
11 & $(1,8),(3,2)$ &  $\frac{1}{70}(1,37)$ & $\frac{3}{35}$ & $\frac{1}{21385}$  & $2$\\
\hline
12 & $(2,13),(3,2)$ &  $\frac{1}{69}(1,47)$ & $\frac{2}{23}$ & $\frac{2}{14053}$  & $2$\\
\hline
13 & $(2,13),(3,2),(3,18)$ &  $\frac{1}{69}(1,47)$ & $\frac{2}{23}$ & $\frac{1}{10465}$ & $3$\\
\hline
14 & $(2,2),(3,13)$ &  $\frac{1}{69}(1,47)$ & $\frac{2}{23}$ & $\frac{1}{10465}$  & $2$\\
\hline
15 & $(2,2)$ &  $\frac{1}{89}(1,61)$ & $\frac{8}{89}$ & $\frac{1}{3115}$  & $1$ \\
\hline
16 & $(2,2),(2,13)$ &  $\frac{1}{89}(1,61)$ & $\frac{8}{89}$ & $\frac{6}{26611}$  & $2$\\
\hline
17 & $(1,2),(2,13)$ &  $\frac{1}{89}(1,61)$ & $\frac{8}{89}$ & $\frac{6}{26611}$  & $2$\\
\hline
18 & $(3,2)$ &  $\frac{1}{79}(1,54)$ & $\frac{7}{79}$ & $\frac{1}{3713}$  & $1$ \\
\hline
19 & $(3,2),(3,13)$ &  $\frac{1}{79}(1,54)$ & $\frac{7}{79}$ & $\frac{3}{23621}$  & $2$ \\
\hline
20 & $(3,2),(3,18)$ &  $\frac{1}{79}(1,54)$ & $\frac{7}{79}$ & $\frac{8}{35945}$  & $2$\\
\hline
21 & $(2,2),(3,18)$ &  $\frac{1}{79}(1,54)$ & $\frac{7}{79}$ & $\frac{8}{35945}$  & $2$\\
\hline
22 & $(2,2),(2,13),(3,18)$ &  $\frac{1}{79}(1,54)$ & $\frac{7}{79}$ & $\frac{3}{23621}$ & $3$\\
\hline
23 & $(1,2),(3,13)$ &  $\frac{1}{79}(1,54)$ & $\frac{7}{79}$  & $\frac{3}{23621}$  & $2$ \\
\hline
24 & $(1,2),(2,13),(3,18)$ &  $\frac{1}{79}(1,54)$ & $\frac{7}{79}$  & $\frac{3}{23621}$ & $3$ \\
\hline
25 & $(1,2),(3,18)$ &  $\frac{1}{88}(1,67)$ & $\frac{1}{11}$  & $\frac{1}{3289}$ & $2$\\
\hline
\end{longtable}
\end{center}
\end{thm}
\begin{proof}
Let $x_0:=\Center_XS_W$. By the classification as in Lemma \ref{lem: 1/11-lt exceptional classification step 1}, for any closed point $x\in X$ such that $x\not=x_0$,
$\mld(X;x)>\frac{1}{11}$. Therefore,
\begin{equation}\label{equ: mld 1/11}
    \mld(X;x_0)\leq\frac{1}{11}.
\end{equation}
Since $X$ is a del Pezzo surface,
\begin{equation}\label{equ: x fano}
    K_X\cdot E_{i,j,X}<0 \text{ if } E_{i,j,X}\not=0.
\end{equation}
Table \ref{table: 1/11-lt fano} is provided as a table of all possibilities of $X$ such that (\ref{equ: mld 1/11}) and (\ref{equ: x fano}) hold. The detailed computation is made by applying Lemma \ref{lem: 1/11-lt exceptional classification step 1} and enumeration. (We remark that the enumeration is actually not very complicated and can be computed by hand. For example, if the dual graph of $X\ni x_0$ contains a single point, then $-26\leq\bar S^2\leq -22$ as $\mld(X;x_0)\leq\frac{1}{11}$. This implies that the ``Change near $S_W$" in Tables \ref{table: contraction E1}, \ref{table: contraction E2}, and \ref{table: contraction E3} should be of the form $(\alpha_1,0),(\alpha_2,0),(\alpha_3,0)$ with $\alpha_1+\alpha_2+\alpha_3\leq 4$. This makes the possible cases very restrictive.)

By Lemmas \ref{lem: non-1/11-lt exceptional fano dominated by termialization} and \ref{lem: 1/11-lt contract divisors}, any non-$\frac{1}{11}$-lt exceptional del Pezzo surface belongs to Table \ref{table: 1/11-lt fano}. We left to show that $X$ is an exceptional del Pezzo surface for any $X$ in Table \ref{table: 1/11-lt fano}. 

First we show that $X$ is exceptional. Note that for any $X$ in Table \ref{table: 1/11-lt fano}, there exists a variety $T$ in Table \ref{table: picard number 1 fano} and a projective birational morphism $h_X: X\rightarrow T$. By Example-Lemma \ref{exlem: special exc fano}, $T$ is exceptional. Suppose that $X$ is not exceptional, then $X$ has an lc but not klt $\Rr$-complement $(X,G)$. Now $(T,(h_X)_*G)$ is an lc but not klt $\Rr$-complement of $T$, so $T$ is not exceptional, a contradiction. Thus $X$ is exceptional for any $X$ in Table \ref{table: 1/11-lt fano}.

Now we show that $X$ is a del Pezzo surface for any $X$ in Table \ref{table: 1/11-lt fano}. Table \ref{table: 1/11-lt fano} shows that $K_X^2>0$. To clarify our notation, in the following, we let $X_n$ be the variety $X$ as in Table \ref{table: 1/11-lt fano}, Case $n$ for any $1\leq n\leq 25$. There are three cases.

\medskip

\noindent\textbf{Case 1}. $n\in\{15,18\}$. Then $\rho(X_n)=1$ so $X_n$ is a del Pezzo surface.

\medskip

\noindent\textbf{Case 2}. $1\leq n\leq 9$. In this case, for any $i,j$ such that $E_{i,j,X_n}\not=0$, $E_{i,j,X_n}^2<0$. Thus $E_{i,j,X_n}$ are the $K_{X_n}$-negative extremal rays. Thus
$$\overline{NE}(X_n)=\sum_{i,j} \mathbb R_+[E_{i,j,X_n}],$$
so $K_{X_n}\cdot C<0$ for any curve $C\subset X_n$. By Kleiman's Criterion, $X_n$ is a del Pezzo surface.

\medskip

\noindent\textbf{Case 3}. $10\leq n\leq 25$ and $n\not\in\{15,18\}$. For $n\in\{10,14,16,21\}$, we have $K_{X_n}\sim_{\mathbb Q}-\frac{1}{13}E_{2,2,X_{n}}$ and $E_{2,2,X_n}^2>0$. For $n\in\{11,12,13,19,20\}$, we have $K_{X_n}\sim_{\mathbb Q}-\frac{1}{13}E_{3,2,X_{n}}$ and $E_{3,2,X_n}^2>0$. For $n\in\{17,23,24,25\}$,  we have $K_{X_n}\sim_{\mathbb Q}-\frac{1}{13}E_{1,2,X_{n}}$ and $E_{1,2,X_n}^2>0$. Thus $-K_{X_n}$ is nef and big. 

Suppose that $X_n$ is not del Pezzo. Then we let $\phi_n: X_n\rightarrow X_n'$ be the ample model of $-K_{X_n}$. Since $K_{X_n}\cdot E_{i,j,X_n}<0$ for any $i,j$ such that $E_{i,j,X_n}\not=0$, $\phi_n$ does not contract any $E_{i,j,X_n}$ such that $E_{i,j,X_n}\not=0$. By construction, $X_n'$ is a del Pezzo surface of Picard number $1$ and $\mld(X_n')=\mld(X_n)\leq\frac{1}{11}$. Since $\phi_n^*K_{X_n'}=K_{X_n}$ and $X_n$ is exceptional, $X_n'$ is exceptional. By Lemma \ref{lem: non-1/11-lt exceptional fano dominated by termialization}(2), the exceptional divisors of $\phi_n$ are $E_{i,j,X_n}$ for some $i,j$, a contradiction. Thus $X_n$ is del Pezzo.

Finally, for any $X$ in Table \ref{table: 1/11-lt fano},  $(X,B:=h_*B_W)$ is a $13$-complement of $X$.
\end{proof}

\begin{proof}[Proof of Theorem \ref{thm: 1/13 exceptional Fano}]
It follows immediately from Theorem \ref{thm: 1/11-lt exceptional classification step 2}.
\end{proof}

\begin{rem}\label{rem: some non-1/11-lt exc del pezzos}
It is easyy to see that our classification, as demonstrated in the proof of Theorem \ref{thm: 1/11-lt exceptional classification step 2}, can be employed to construct additional exceptional Fano surfaces whose mlds range between $\frac{1}{11}$ and $\leq\frac{1}{7}$. For instance, Case 1 of Table \ref{table: picard number 1 fano} is an exceptional Fano surface with an mld equal to $\frac{9}{97}$.

We should note that certain exceptional surfaces constructed in \cite{CPS10} also have mlds between $\frac{1}{11}$ and $\frac{1}{7}$. These surfaces can be constructed in the same manner as presented in the proof of Theorem \ref{thm: 1/11-lt exceptional classification step 2}, but with a different set-up comparing to Set-up \ref{setup: classify non-1/11-lt}. The primary reason for this is that, for such a surface $X$ and the extraction $f: Y\rightarrow X$ of the divisor $E$ on $X$ computing $\mld(X)$, it is possible for $\Rct(Y,0;E)\leq\frac{10}{11}$. Hence, we cannot apply Lemma \ref{lem: non-1/11-lt exceptional fano dominated by termialization}. For example, the exceptional del Pezzo surface with the smallest mld in \cite[The Long Table]{CPS10} is a degree $256$ hypersurface in $\mathbb P(11,49,69,128)$, with an mld equal to $\frac{5}{49}$. In this case, it can be computed that $\Rct(Y,0;E)=\frac{10}{11}$.

On the other hand, once we classify all klt log Calabi-Yau surface pairs $(X,B=bS)$ such that $S$ is a non-zero Weil divisor and $b>\frac{6}{7}$ (which we anticipate can be accomplished using similar arguments as in Sections \ref{sec: sings on S} and \ref{sec: proof of the gap theorem} although more computational algorithms are needed), we can expect to employ analogous reasoning as presented in Lemmas \ref{lem: non-1/11-lt exceptional fano dominated by termialization} and \ref{lem: 1/11-lt exceptional classification step 1}, and Theorem \ref{thm: 1/11-lt exceptional classification step 2}. Consequently, we can anticipate a classification of exceptional del Pezzo surfaces with mld's $\leq\frac{1}{7}$ or even $\leq\frac{1}{6}$, with a different Set-up \ref{setup: classify non-1/11-lt}, using similar methods of the paper. In particular, the degree $256$ hypersurface in $\mathbb P(11,49,69,128)$ in \cite{CPS10} is expected to be obtained by using similar methods as in our paper.
\end{rem}

\begin{rem}
We can also classify all exceptional Fano type surfaces that are not $\frac{1}{11}$-lt. Indeed, by Lemma \ref{lem: non-1/11-lt exceptional fano dominated by termialization}, all such exceptional Fano type surfaces are obtained by extracting divisors of the form $E_{i,j}$ over some $X$ listed in Table \ref{table: 1/11-lt fano}. In particular, there are only finitely many exceptional Fano type surfaces that are not $\frac{1}{11}$-lt.

On the other hand, any surface obtained by extracting divisors of the form $E_{i,j}$ over some $X$ listed in Table \ref{table: 1/11-lt fano} is exceptional and of Fano type, but not all of them are not $\frac{1}{11}$-lt. For simplicity, we will not make a detailed classification here. Nevertheless, as an immediate application of Lemma \ref{lem: non-1/11-lt exceptional fano dominated by termialization}, we have Corollary \ref{cor: classification non-1/11-lt exceptional fano type} below.
\end{rem}

\begin{cor}\label{cor: classification non-1/11-lt exceptional fano type}
Conditions as in Set-up \ref{setup: terminalization of example}. 
\begin{enumerate}
    \item For any exceptional Fano type surface $\tilde X$ that is not $\frac{1}{11}$-lt, there exists a surface $X$ in Table \ref{table: 1/11-lt fano} and a projective birational morphism $\tilde X\rightarrow X$ which only contracts the images of $E_{i,j}$ for some $i,j$.
    \item For any $X$ in Table \ref{table: 1/11-lt fano} and projective birational morphism $\tilde X\rightarrow X$ which only contracts the images of $E_{i,j}$ for some $i,j$, $\tilde X$ is an exceptional surface.
    \item Exceptional Fano type surfaces that are not $\frac{1}{11}$-lt have isolated moduli, i.e. there are only finitely many such varieties.
    \item The smallest mld of exceptional Fano type surfaces is $\frac{3}{35}$.
     \item The smallest mld of exceptional Fano type surfaces of Picard number $1$ is $\frac{7}{79}$.
    \item Any exceptional Fano type surface that is not $\frac{1}{11}$-lt has a $13$-complement.
\end{enumerate}
\end{cor}
\begin{proof}
(1) By Lemma \ref{lem: non-1/11-lt exceptional fano dominated by termialization} and Theorem \ref{thm: 1/11-lt exceptional classification step 2}, we may let $\tilde X\rightarrow X$ be the ample model of $-K_{\tilde X}$. 

(2) By Theorem \ref{thm: 1/11-lt exceptional classification step 2}, $X$ is exceptional, so $\tilde X$ is exceptional. 

(3) It immediately follows from (1).

(4) Let $X$ be an exceptional Fano type surface. We may run a $(-K_X)$-MMP $X\rightarrow X'$ which terminates with a model $X'$. Then $X'$ is an exceptional surface such that $-K_X$ is big and nef and semi-ample. Now we let $X''$ be the ample model of $X'$, then $X''$ is an exceptional del Pezzo surface. By Theorem \ref{thm: 1/13 exceptional Fano}(2), we have
$$\mld(X)\geq\mld(X')\geq\mld(X'')\geq\frac{3}{35}.$$

(5) Exceptional Fano type surfaces of Picard number $1$ are del Pezzo surfaces, so it follows from Theorem \ref{thm: 1/13 exceptional Fano}(3).

(6) By (1), we may let $\tilde h: W\rightarrow\tilde X$ be the induced projective birational morphism, then $(\tilde X,\tilde h_*B_W)$ is a $13$-complement of $\tilde X$.
\end{proof}

\begin{proof}[Proof of Corollary \ref{cor: new bound alpha invariant}] 
(1) By Theorem \ref{thm: 1/13 exceptional Fano}, we may assume that $X$ is $\frac{1}{11}$-lt. (1) follows from Theorem \ref{thm: 1/13 exceptional Fano} and \cite[Lemma 2.7]{Liu23}. 

(2) By \cite[Theorem 4.6(1)(2)]{Liu23}, we may assume that $X$ is an exceptional Fano type surface. Possibly replacing $X$ with the ample model of $-K_X$, we may assume that $X$ is exceptional del Pezzo.  By Theorem \ref{thm: 1/13 exceptional Fano}, we may assume that $X$ is $\frac{1}{11}$-lt. (2) follows from (1) and the effective base-point-freeness theorem (\cite[Theorem 1.1, Remark 1.2]{Fuj09} \cite[1.1 Theorem]{Kol93}).

(3) By \cite[Theorem 4.6]{Liu23}, we may assume that $X$ is an exceptional Fano type surface. Possibly replacing $X$ with the ample model of $-K_X$, we may assume that $X$ is an exceptional del Pezzo surface.  By Theorem \ref{thm: 1/13 exceptional Fano}, we may assume that $X$ is $\frac{1}{11}$-lt. By (1), $K_X^2>\frac{1}{I_0}$. By  \cite[Theorems A,D]{BJ20}), $\alpha(X)<3\sqrt{I_0}$.
\end{proof}

\begin{rem}\label{rem: vol 1/11 exceptional fano}
The anti-canonical volume of the exceptional Fano surface with $\mld=\frac{3}{35}$ (Case 11 of Table \ref{table: 1/11-lt fano}), $\frac{1}{21385}$, is the second smallest known volume of exceptional del Pezzo surfaces. It is smaller than the volume of any exceptional del Pezzo surface in \cite{CPS10} but is larger than Totaro's latest example of a degree $354$ non-quasi-smooth hypersurface in $\mathbb P(177,118,49,11)$ which has anti-canonical volume $\frac{1}{31801}$ \cite{Tot23}.
\end{rem}

\begin{rem}
\cite{Lac20} has classified all exceptional Fano surfaces of Picard number $1$ . In Table \ref{table: 1/11-lt fano}, Case 15 and Case 18 are exceptional Fano surfaces of Picard number $1$, so they should be enlisted in \cite{Lac20}. J. Lacini informed us that Case 15 of Table \ref{table: 1/11-lt fano} corresponds to \cite[LDP16, $(s,r)=(3,1)$, (3), $k=8$]{Lac20}. We can also check that Case 18 of Table \ref{table: 1/11-lt fano} corresponds to \cite[LDP16, $(s,r)=(3,1)$, (1), $k=9$]{Lac20}.
\end{rem}

\appendix

\section{Classification of some singularities related to Proposition  \ref{prop: four singularities, ii}}\label{sec: classification of x4}

In this appendix, we make a rough classification of all possible $x_4$ that satisfy the conditions of Lemma \ref{lem: summary four singular points condition, II}. This will significantly reduce the possibilities of singularities in Table \ref{table: 4 singular points, II} and satisfied the proof of Proposition \ref{prop: four singularities, ii}.

\begin{lem}\label{lem: x4 not a type}
Conditions as in Set-up \ref{setup: rho=1 set-up}. If $x_4$ is not a cyclic quotient singularity, then $x_4$ is of one of the following types listed in the following table:
\begin{center}
    \begin{longtable}{|l|l|l|l|}
 			\caption{Not $A$-type}\label{d type table i}\\
		\hline
	No.&Dual graph & $\gamma(X\ni x)$ & Order\\ 
		\hline
		1 &
    \begin{tikzpicture}
         \draw (3.3,0) circle (0.1);
         \node [below] at (3.3,-0.1) {\footnotesize$2$};
         \draw (3.4,0)--(3.8,0);

         \draw (3.9,0) circle (0.1);
         \draw (3.9,0.6) circle (0.1);
         \draw (4.5,0) circle (0.1);
         \node [below] at (4.5,-0.1) {\footnotesize$2$};

         \draw (3.9,0.1)--(3.9,0.5);
         \draw (4.0,0)--(4.4,0);

         \node [below] at (3.9,-0.1) {\footnotesize$4$};
         \node [right] at (4,0.6) {\footnotesize$2$};
    \end{tikzpicture}
    & $\frac{12}{5}$ & $40$ \\
\hline
2 &    \begin{tikzpicture}
         \draw (3.3,0) circle (0.1);
         \node [below] at (3.3,-0.1) {\footnotesize$2$};
         \draw (3.4,0)--(3.8,0);

         \draw (3.9,0) circle (0.1);
         \draw (3.9,0.6) circle (0.1);
         \draw (4.5,0) circle (0.1);
         \node [below] at (4.5,-0.1) {\footnotesize$2$};

         \draw (3.9,0.1)--(3.9,0.5);
         \draw (4.0,0)--(4.4,0);

         \node [below] at (3.9,-0.1) {\footnotesize$3$};
         \node [right] at (4,0.6) {\footnotesize$2$};
    \end{tikzpicture} & $\frac{10}{3}$ & $24$ \\
    	\hline
3 & \begin{tikzpicture}
         \draw (3.3,0) circle (0.1);
         \node [below] at (3.3,-0.1) {\footnotesize$3$};
         \draw (3.4,0)--(3.8,0);

         \draw (3.9,0) circle (0.1);
         \draw (3.9,0.6) circle (0.1);
         \draw (4.5,0) circle (0.1);
         \node [below] at (4.5,-0.1) {\footnotesize$2$};

         \draw (3.9,0.1)--(3.9,0.5);
         \draw (4.0,0)--(4.4,0);

         \node [below] at (3.9,-0.1) {\footnotesize$2$};
         \node [right] at (4,0.6) {\footnotesize$2$};
    \end{tikzpicture}  & $\frac{7}{2}$  & $24$\\
\hline
4 & \begin{tikzpicture}
         \draw (2.7,0) circle (0.1);
         \node [below] at (2.7,-0.1) {\footnotesize$3$};
         \draw (2.8,0)--(3.2,0);
         \draw (3.3,0) circle (0.1);
         \node [below] at (3.3,-0.1) {\footnotesize$2$};
         \draw (3.4,0)--(3.8,0);

         \draw (3.9,0) circle (0.1);
         \draw (3.9,0.6) circle (0.1);
         \draw (4.5,0) circle (0.1);
         \node [below] at (4.5,-0.1) {\footnotesize$2$};

         \draw (3.9,0.1)--(3.9,0.5);
         \draw (4.0,0)--(4.4,0);

         \node [below] at (3.9,-0.1) {\footnotesize$2$};
         \node [right] at (4,0.6) {\footnotesize$2$};
    \end{tikzpicture}  & $\frac{9}{2}$ & $40$\\
    \hline
		\end{longtable}
		\end{center}
\end{lem}
\begin{proof}
By \cite[Satz 2.11]{Bri68}, the order of any other non-cyclic quotient singularity is $>42$. The lemma now follows from Lemma \ref{lem: summary four singular points condition, II}(2).
\end{proof}

\begin{lem}
Conditions as in Set-up \ref{setup: rho=1 set-up} and assume that $x_4\not\in S$. Suppose that
\begin{enumerate}
    \item the order of $X\ni x_4$ is $\leq 24$,
    \item $\gamma(X\ni x)<\frac{887}{77}$, and
    \item $X\ni x_4$ is not a Du Val singularity.
\end{enumerate}
Then $(X\ni x_4)$ is of one of the types listed in Table \ref{table order<=24 A type}.

\begin{center}
    \begin{longtable}{|l|l|l|l|l|}
 			\caption{A type, I}\label{table order<=24 A type}\\
		\hline
		No.&Singularity& Weights & L.d. & $\gamma(X\ni x_4)$\\ 
		\hline
1 & $\frac{1}{3}(1,1)$ & $(3)$ & $(\frac{2}{3})$ & $\frac{2}{3}$\\
\hline
2 & $\frac{1}{4}(1,1)$ & $(4)$ & $(\frac{1}{2})$ & $0$\\
\hline
3 & $\frac{1}{5}(1,1)$ & $(5)$ & $(\frac{2}{5})$ & $-\frac{4}{5}$\\
\hline
3 & $\frac{1}{5}(1,2)$ & $(3,2)$ & $(\frac{3}{5},\frac{4}{5})$ & $\frac{8}{5}$\\
\hline
4 & $\frac{1}{6}(1,1)$ & $(6)$ & $(\frac{1}{3})$ & $-\frac{5}{3}$ \\
\hline
5 & $\frac{1}{7}(1,1)$ & $(7)$ & $(\frac{2}{7})$ & $-\frac{18}{7}$ \\
\hline
6 & $\frac{1}{7}(1,2)$ & $(4,2)$ & $(\frac{3}{7},\frac{5}{7})$ & $\frac{6}{7}$ \\
\hline
7 & $\frac{1}{7}(1,3)$ & $(3,2,2)$ & $(\frac{4}{7},\frac{5}{7},\frac{6}{7})$ & $\frac{18}{7}$ \\
\hline
8 & $\frac{1}{8}(1,1)$ & $(8)$ & $(\frac{1}{4})$ & $-\frac{7}{2}$ \\
\hline
9 & $\frac{1}{8}(1,3)$ & $(3,3)$ & $(\frac{1}{2},\frac{1}{2})$ & $1$ \\
\hline
10 & $\frac{1}{8}(1,5)$ & $(2,3,2)$ & $(\frac{3}{4},\frac{1}{2},\frac{3}{4})$ & $\frac{5}{2}$ \\
\hline
11 & $\frac{1}{9}(1,1)$ & $(9)$ & $(\frac{2}{9})$ & $-\frac{40}{9}$ \\
\hline
12 & $\frac{1}{9}(1,2)$ & $(5,2)$ & $(\frac{1}{3},\frac{2}{3})$ & $0$ \\
\hline
13 & $\frac{1}{9}(1,4)$ & $(3,2,2,2)$ & $(\frac{5}{9},\frac{2}{3},\frac{7}{9},\frac{8}{9})$ & $\frac{32}{9}$ \\
\hline
14 & $\frac{1}{10}(1,1)$ & $(10)$ & $(\frac{1}{5})$ & $-\frac{27}{5}$ \\
\hline
15 & $\frac{1}{10}(1,3)$ & $(4,2,2)$ & $(\frac{2}{5},\frac{3}{5},\frac{4}{5})$ & $\frac{9}{5}$ \\
\hline
16 & $\frac{1}{11}(1,1)$ & $(11)$ & $(\frac{2}{11})$ & $-\frac{70}{11}$ \\
\hline
17 & $\frac{1}{11}(1,2)$ & $(6,2)$ & $(\frac{3}{11},\frac{7}{11})$ & $-\frac{10}{11}$ \\
\hline
18 & $\frac{1}{11}(1,3)$ & $(4,3)$ & $(\frac{4}{11},\frac{5}{11})$ & $\frac{2}{11}$ \\
\hline
19 & $\frac{1}{11}(1,5)$ & $(3,2,2,2,2)$ & $(\frac{5+i}{11})_{i=1}^5$ & $\frac{50}{11}$ \\
\hline
20 & $\frac{1}{11}(1,7)$ & $(2,3,2,2)$ & $(\frac{8}{11},\frac{5}{11},\frac{7}{11},\frac{9}{11})$ & $\frac{38}{11}$ \\
\hline
21 & $\frac{1}{12}(1,1)$ & $(12)$ & $(\frac{1}{6})$ & $-\frac{22}{3}$ \\
\hline
22 & $\frac{1}{12}(1,5)$ & $(3,2,3)$ & $(\frac{1}{2},\frac{1}{2},\frac{1}{2})$ & $2$ \\
\hline
23 & $\frac{1}{12}(1,7)$ & $(2,4,2)$ & $(\frac{2}{3},\frac{1}{3},\frac{2}{3})$ & $\frac{5}{3}$ \\
\hline
24 & $\frac{1}{13}(1,1)$ & $(13)$ & $(\frac{2}{13})$ & $-\frac{108}{13}$ \\
\hline
25 & $\frac{1}{13}(1,2)$ & $(7,2)$ & $(\frac{3}{13},\frac{8}{13})$ & $-\frac{24}{13}$ \\
\hline
26 & $\frac{1}{13}(1,3)$ & $(5,2,2)$ & $(\frac{4}{13},\frac{7}{13},\frac{10}{13})$ & $\frac{12}{13}$ \\
\hline
27 & $\frac{1}{13}(1,4)$ & $(4,2,2,2)$ & $(\frac{5}{13},\frac{7}{13},\frac{9}{13},\frac{11}{13})$ & $\frac{36}{13}$ \\
\hline
28 & $\frac{1}{13}(1,5)$ & $(3,3,2)$ & $(\frac{6}{13},\frac{5}{13},\frac{9}{13})$ & $\frac{24}{13}$ \\
\hline
29 & $\frac{1}{13}(1,6)$ & $(3,2^5)$ & $(\frac{6+i}{13})_{i=1}^6$ & $\frac{72}{13}$ \\
\hline
30 & $\frac{1}{14}(1,3)$ & $(5,3)$ & $(\frac{2}{7},\frac{3}{7})$ & $-\frac{5}{7}$ \\
\hline
31 & $\frac{1}{14}(1,9)$ & $(2,3,2,2,2)$ & $(\frac{5}{7},\frac{3}{7},\frac{4}{7},\frac{5}{7},\frac{6}{7})$ & $\frac{31}{7}$ \\
\hline
32 & $\frac{1}{15}(1,2)$ & $(8,2)$ & $(\frac{1}{5},\frac{3}{5})$ & $-\frac{14}{5}$ \\
\hline
33 & $\frac{1}{15}(1,4)$ & $(4,4)$ & $(\frac{1}{3},\frac{1}{3})$ & $-\frac{2}{3}$ \\
\hline
34 & $\frac{1}{15}(1,7)$ & $(3,2^6)$ & $(\frac{7+i}{15})_{i=1}^7$ & $\frac{98}{15}$ \\
\hline
35 & $\frac{1}{15}(1,11)$ & $(2,2,3,2,2)$ & $(\frac{4}{5},\frac{3}{5},\frac{2}{5},\frac{3}{5},\frac{4}{5})$ & $\frac{22}{5}$ \\
\hline
36 & $\frac{1}{16}(1,3)$ & $(6,2,2)$ & $(\frac{1}{4},\frac{1}{2},\frac{3}{4})$ & $0$ \\
\hline
37 & $\frac{1}{16}(1,5)$ & $(4,2,2,2,2)$ & $(\frac{2+i}{8})_{i=1}^5$ & $\frac{15}{4}$ \\
\hline
38 & $\frac{1}{16}(1,7)$ & $(3,2,2,3)$ & $(\frac{1}{2},\frac{1}{2},\frac{1}{2},\frac{1}{2})$ & $3$ \\
\hline
39 & $\frac{1}{16}(1,9)$ & $(2,5,2)$ & $(\frac{5}{8},\frac{1}{4},\frac{5}{8})$ & $\frac{3}{4}$ \\
\hline
40 & $\frac{1}{17}(1,2)$ & $(9,2)$ & $(\frac{3}{17},\frac{10}{17})$ & $-\frac{64}{17}$ \\
\hline
41 & $\frac{1}{17}(1,3)$ & $(6,3)$ & $(\frac{4}{17},\frac{7}{17})$ & $-\frac{28}{17}$ \\
\hline
42 & $\frac{1}{17}(1,4)$ & $(5,2,2,2)$ & $(\frac{5}{17},\frac{8}{17},\frac{11}{17},\frac{14}{17})$ & $\frac{22}{17}$ \\
\hline
43 & $\frac{1}{17}(1,5)$ & $(4,2,3)$ & $(\frac{6}{17},\frac{7}{17},\frac{8}{17})$ & $\frac{3}{17}$ \\
\hline
44 & $\frac{1}{17}(1,8)$ & $(3,2^7)$ & $(\frac{8+i}{17})_{i=1}^8$ & $\frac{128}{17}$ \\
\hline
45 & $\frac{1}{17}(1,10)$ & $(2,4,2,2)$ & $(\frac{11}{17},\frac{5}{17},\frac{9}{17},\frac{13}{17})$ & $\frac{44}{17}$ \\
\hline
46 & $\frac{1}{17}(1,11)$ & $(2,3,2^4)$ & $(\frac{12}{17},(\frac{5+2i}{17})_{i=1}^5)$ & $\frac{92}{17}$ \\
\hline
47 & $\frac{1}{18}(1,5)$ & $(4,3,2)$ & $(\frac{1}{3},\frac{1}{3},\frac{2}{3})$ & $1$ \\
\hline
48 & $\frac{1}{18}(1,7)$ & $(3,3,2,2)$ & $(\frac{4}{9},\frac{1}{3},\frac{5}{9},\frac{7}{9})$ & $\frac{25}{9}$ \\
\hline
49 & $\frac{1}{19}(1,2)$ & $(10,2)$ & $(\frac{3}{19},\frac{11}{19})$ & $-\frac{90}{19}$ \\
\hline
50 & $\frac{1}{19}(1,3)$ & $(7,2,2)$ & $(\frac{4}{19},\frac{9}{19},\frac{14}{19})$ & $-\frac{18}{19}$ \\
\hline
51 & $\frac{1}{19}(1,4)$ & $(5,4)$ & $(\frac{5}{19},\frac{6}{19})$ & $-\frac{30}{19}$ \\
\hline
52 & $\frac{1}{19}(1,6)$ & $(4,2^5)$ & $(\frac{5+2i}{19})_{i=1}^6$ & $\frac{90}{19}$ \\
\hline
53 & $\frac{1}{19}(1,7)$ & $(3,4,2)$ & $(\frac{8}{19},\frac{5}{19},\frac{12}{19})$ & $\frac{18}{19}$ \\
\hline
54 & $\frac{1}{19}(1,8)$ & $(3,2,3,2)$ & $(\frac{9}{19},\frac{8}{19},\frac{7}{19},\frac{13}{19})$ & $\frac{54}{19}$ \\
\hline
55 & $\frac{1}{19}(1,9)$ & $(3,2^8)$ & $(\frac{9+i}{19})_{i=1}^9$ & $\frac{162}{19}$ \\
\hline
56 & $\frac{1}{19}(1,14)$ & $(2,2,3,2,2,2)$ & $(\frac{15}{19},\frac{11}{19},(\frac{4+3i}{19})_{i=1}^4)$ & $\frac{102}{19}$ \\
\hline
57 & $\frac{1}{20}(1,3)$ & $(7,3)$ & $(\frac{1}{5},\frac{2}{5})$ & $-\frac{13}{5}$ \\
\hline
58 & $\frac{1}{20}(1,9)$ & $(3,2,2,2,3)$ & $(\frac{1}{2})_{i=1}^5$ & $4$ \\
\hline
59 & $\frac{1}{20}(1,11)$ & $(2,6,2)$ & $(\frac{3}{5},\frac{1}{5},\frac{3}{5})$ & $-\frac{1}{5}$ \\
\hline
60 & $\frac{1}{20}(1,13)$ & $(2,3,2^5)$ & $(\frac{7}{10},(\frac{3+i}{10})_{i=1}^6)$ & $\frac{32}{5}$ \\
\hline
61 & $\frac{1}{21}(1,4)$ & $(6,2,2,2)$ & $(\frac{5}{21},\frac{3}{7},\frac{13}{21},\frac{17}{21})$ & $\frac{20}{21}$ \\
\hline
62 & $\frac{1}{21}(1,5)$ & $(5,2,2,2,2)$ & $(\frac{1+i}{7})_{i=1}^5$ & $\frac{20}{7}$ \\
\hline
63 & $\frac{1}{21}(1,8)$ & $(3,3,3)$ & $(\frac{3}{7},\frac{2}{7},\frac{3}{7})$ & $\frac{8}{7}$ \\
\hline
64 & $\frac{1}{21}(1,10)$ & $(3,2^9)$ & $(\frac{10+i}{21})_{i=1}^{10}$ & $\frac{200}{21}$ \\
\hline
65 & $\frac{1}{21}(1,13)$ & $(2,3,3,2)$ & $(\frac{2}{3},\frac{1}{3},\frac{1}{3},\frac{2}{3})$ & $\frac{8}{3}$ \\
\hline
66 & $\frac{1}{22}(1,3)$ & $(8,2,2)$ & $(\frac{2}{11},\frac{5}{11},\frac{8}{11})$ & $-\frac{21}{11}$ \\
\hline
67 & $\frac{1}{22}(1,5)$ & $(5,2,3)$ & $(\frac{3}{11},\frac{4}{11},\frac{5}{11})$ & $\frac{3}{11}$ \\
\hline
68 & $\frac{1}{22}(1,7)$ & $(4,2^6)$ & $(\frac{3+i}{11})_{i=1}^7$ & $\frac{63}{11}$ \\
\hline
69 & $\frac{1}{22}(1,13)$ & $(2,4,2,2,2)$ & $(\frac{7}{11},(\frac{1+2i}{11})_{i=1}^4)$ & $\frac{39}{11}$ \\
\hline
70 & $\frac{1}{23}(1,3)$ & $(8,3)$ & $(\frac{4}{23},\frac{9}{23})$ & $-\frac{82}{23}$ \\
\hline
71 & $\frac{1}{23}(1,4)$ & $(6,4)$ & $(\frac{5}{23},\frac{7}{23})$ & $-\frac{58}{23}$ \\
\hline
72 & $\frac{1}{23}(1,5)$ & $(5,3,2)$ & $(\frac{6}{23},\frac{7}{23},\frac{15}{23})$ & $\frac{2}{23}$ \\
\hline
73 & $\frac{1}{23}(1,7)$ & $(4,2,2,3)$ & $(\frac{8}{23},\frac{9}{23},\frac{10}{23},\frac{11}{23})$ & $\frac{50}{23}$ \\
\hline
74 & $\frac{1}{23}(1,9)$ & $(3,3,2,2,2)$ & $(\frac{10}{23},\frac{7}{23},\frac{11}{23},\frac{15}{23},\frac{19}{23})$ & $\frac{86}{23}$ \\
\hline
75 & $\frac{1}{23}(1,11)$ & $(3,2^{10})$ & $(\frac{11+i}{23})_{i=1}^{11}$ & $\frac{242}{23}$ \\
\hline
76 & $\frac{1}{23}(1,13)$ & $(2,5,2,2)$ & $(\frac{14}{23},\frac{5}{23},\frac{11}{23},\frac{17}{23})$ & $\frac{38}{23}$ \\
\hline
77 & $\frac{1}{23}(1,15)$ & $(2,3,2^6)$ & $(\frac{16}{23},(\frac{7+2i}{23})_{i=1}^7)$ & $\frac{170}{23}$ \\
\hline
78 & $\frac{1}{23}(1,17)$ & $(2,2,3,2^4)$ & $(\frac{18}{23},\frac{13}{23},(\frac{5+3i}{23})_{i=1}^5)$ & $\frac{146}{23}$ \\
\hline
79 & $\frac{1}{24}(1,5)$ & $(5,5)$ & $(\frac{1}{4},\frac{1}{4})$ & $-\frac{5}{2}$ \\
\hline
80 & $\frac{1}{24}(1,7)$ & $(4,2,4)$ & $(\frac{1}{3},\frac{1}{3},\frac{1}{3})$ & $\frac{1}{3}$ \\
\hline
81 & $\frac{1}{24}(1,11)$ & $(3,2^4,3)$ & $(\frac{1}{2})_{i=1}^6$ & $5$ \\
\hline
82 & $\frac{1}{24}(1,13)$ & $(2,7,2)$ & $(\frac{7}{12},\frac{1}{6},\frac{7}{12})$ & $-\frac{7}{6}$ \\
\hline
83 & $\frac{1}{24}(1,17)$ & $(2,2,4,2,2)$ & $(\frac{3}{4},\frac{1}{2},\frac{1}{4},\frac{1}{2},\frac{3}{4})$ & $\frac{7}{2}$ \\
\hline
84 & $\frac{1}{24}(1,19)$ & $(2,2,2,3,2,2,2)$ & $(\frac{5}{6},\frac{2}{3},\frac{1}{2},\frac{1}{3},\frac{1}{2},\frac{2}{3},\frac{5}{6})$ & $\frac{19}{3}$ \\
\hline
\end{longtable}
\end{center}
\end{lem}

\begin{lem}
Conditions as in Set-up \ref{setup: rho=1 set-up} and assume that $x_4\not\in S$. Suppose that
\begin{enumerate}
    \item the order of $X\ni x_4$ is $\geq 25$,
    \item $\gamma(X\ni x_4)<\frac{887}{77}$, and
    \item $X\ni x_4$ is not a Du Val singularity.
\end{enumerate}
Then the order of $X\ni x_4$ belongs to  $\{25,26,29,31,34,37,38,39,41\}$, and $(X\ni x_4)$ is of one of the types listed in Table \ref{table order>=25 A type}.
\begin{center}
    \begin{longtable}{|l|l|l|l|l|}
 			\caption{A type, II}\label{table order>=25 A type}\\
		\hline
		No.&Singularity& Weights & L.d. & $\gamma(X\ni x_4)$\\ 
		\hline
1 & $\frac{1}{25}(1,3)$ & $(9,2,2)$ & $(\frac{4}{25},\frac{11}{25},\frac{18}{25})$ & $-\frac{72}{25}$ \\
\hline
2 & $\frac{1}{25}(1,4)$ & $(7,2,2,2)$ & $(\frac{i}{5})_{i=1}^4$ & $0$ \\
\hline
3 & $\frac{1}{25}(1,6)$ & $(5,2^5)$ & $(\frac{4+3i}{25})_{i=1}^6$ & $\frac{96}{25}$ \\
\hline
4 & $\frac{1}{25}(1,7)$ & $(4,3,2,2)$ & $(\frac{8}{25},\frac{7}{25},\frac{13}{25},\frac{19}{25})$ & $\frac{48}{25}$ \\
\hline
5 & $\frac{1}{25}(1,8)$ & $(4,2^7)$ & $(\frac{7+2i}{25})_{i=1}^8$ & $\frac{168}{25}$ \\
\hline
6 & $\frac{1}{25}(1,9)$ & $(3,5,2)$ & $(\frac{2}{5},\frac{1}{5},\frac{3}{5})$ & $0$ \\
\hline
7 & $\frac{1}{25}(1,11)$ & $(3,2,2,3,2)$ & $(\frac{12}{25},\frac{11}{25},\frac{2}{5},\frac{9}{25},\frac{17}{25})$ & $\frac{96}{25}$ \\
\hline
8 & $\frac{1}{25}(1,12)$ & $(3,2^{11})$ & $(\frac{12+i}{25})_{i=1}^{12}$ & $\frac{288}{25}$ \\
\hline
9 & $\frac{1}{26}(1,3)$ & $(9,3)$ & $(\frac{2}{13},\frac{5}{13})$ & $-\frac{59}{13}$ \\
\hline
10 & $\frac{1}{26}(1,5)$ & $(6,2,2,2,2)$ & $(\frac{1+2i}{13})_{i=1}^5$ & $\frac{25}{13}$ \\
\hline
11 & $\frac{1}{26}(1,7)$ & $(4,4,2)$ & $(\frac{4}{13},\frac{3}{13},\frac{8}{13})$ & $\frac{1}{13}$ \\
\hline
12 & $\frac{1}{26}(1,11)$ & $(3,2,3,2,2)$ & $(\frac{6}{13},\frac{5}{13},\frac{4}{13},\frac{7}{13},\frac{10}{13})$ & $\frac{49}{13}$ \\
\hline
13 & $\frac{1}{26}(1,17)$ & $(2,3,2^7)$ & $(\frac{9}{13},(\frac{4+i}{13})_{i=1}^8)$ & $\frac{109}{13}$ \\
\hline
14 & $\frac{1}{29}(1,4)$ & $(8,2,2,2)$ & $(\frac{5}{29},\frac{11}{29},\frac{17}{29},\frac{23}{29})$ & $-\frac{28}{29}$\\
\hline
15 & $\frac{1}{29}(1,5)$ & $(6,5)$ & $(\frac{6}{29},\frac{7}{29})$ & $-\frac{100}{29}$\\
\hline
16 & $\frac{1}{29}(1,7)$ & $(5,2^6)$ & $(\frac{5+3i}{29})_{i=1}^7$ & $\frac{140}{29}$\\
\hline
17 & $\frac{1}{29}(1,8)$ & $(4,3,3)$ & $(\frac{9}{29},\frac{7}{29},\frac{12}{29})$ & $\frac{8}{29}$\\
\hline
18 & $\frac{1}{29}(1,9)$ & $(4,2,2,2,3)$ & $(\frac{10}{29},\frac{11}{29},\frac{12}{29},\frac{13}{29},\frac{14}{29})$ & $\frac{92}{29}$\\
\hline
19 & $\frac{1}{29}(1,12)$ & $(3,2,4,2)$ & $(\frac{13}{29},\frac{10}{29},\frac{7}{29},\frac{18}{29})$ & $\frac{56}{29}$\\
\hline
20 & $\frac{1}{29}(1,16)$ & $(2,6,2,2)$ & $(\frac{17}{29},\frac{5}{29},\frac{13}{29},\frac{21}{29})$ & $\frac{20}{29}$\\
\hline
21 & $\frac{1}{29}(1,18)$ & $(2,3,3,2,2)$ & $(\frac{19}{29},\frac{9}{29},\frac{8}{29},\frac{15}{29},\frac{22}{29})$ & $\frac{104}{29}$\\
\hline
22 & $\frac{1}{29}(1,19)$ & $(2,3,2^8)$ & $(\frac{20}{29},(\frac{9+2i}{29})_{i=1}^9)$ & $\frac{272}{29}$\\
\hline
23 & $\frac{1}{29}(1,23)$ & $(2,2,2,3,2^4)$ & $(\frac{24}{29},\frac{19}{29},\frac{14}{29},(\frac{5+4i}{29})_{i=1}^5)$ & $\frac{212}{29}$\\
\hline
24 & $\frac{1}{31}(1,4)$ & $(8,4)$ & $(\frac{5}{31},\frac{9}{31})$ & $-\frac{86}{31}$ \\
\hline
25 & $\frac{1}{31}(1,5)$ & $(7,2^4)$ & $(\frac{1+5i}{31})_{i=1}^5$ & $\frac{30}{31}$ \\
\hline
26 & $\frac{1}{31}(1,6)$ & $(6,2^5)$ & $(\frac{3+4i}{31})_{i=1}^6$ & $\frac{90}{31}$ \\
\hline
27 & $\frac{1}{31}(1,7)$ & $(5,2,4)$ & $(\frac{8}{31},\frac{9}{31},\frac{10}{31})$ & $-\frac{18}{31}$ \\
\hline
28 & $\frac{1}{31}(1,10)$ & $(4,2^9)$ & $(\frac{9+2i}{31})_{i=1}^{10}$ & $\frac{270}{31}$ \\
\hline
29 & $\frac{1}{31}(1,11)$ & $(3,6,2)$ & $(\frac{12}{31},\frac{5}{31},\frac{18}{31})$ & $-\frac{30}{31}$ \\
\hline
30 & $\frac{1}{31}(1,12)$ & $(3,3,2,3)$ & $(\frac{13}{31},\frac{8}{31},\frac{11}{31},\frac{14}{31})$ & $\frac{66}{31}$ \\
\hline
31 & $\frac{1}{31}(1,14)$ & $(3,2,2,2,3,2)$ & $(\frac{15}{31},\frac{14}{31},\frac{13}{31},\frac{12}{31},\frac{11}{31},\frac{21}{31})$ & $\frac{150}{31}$ \\
\hline
32 & $\frac{1}{31}(1,18)$ & $(2,4,3,2)$ & $(\frac{19}{31},\frac{7}{31},\frac{9}{31},\frac{20}{31})$ & $\frac{54}{31}$ \\
\hline
33 & $\frac{1}{31}(1,22)$ & $(2,2,4,2,2,2)$ & $(\frac{23}{31},\frac{15}{31},\frac{7}{31},\frac{13}{31},\frac{19}{31},\frac{25}{31})$ & $\frac{138}{31}$ \\
\hline
34 & $\frac{1}{31}(1,23)$ & $(2,2,3,2^6)$ & $(\frac{24}{31},\frac{17}{31},(\frac{7+3i}{31})_{i=1}^7)$ & $\frac{258}{31}$\\
\hline
35 & $\frac{1}{34}(1,5)$ & $(7,5)$ & $(\frac{3}{17},\frac{4}{17})$ & $-\frac{101}{17}$ \\
\hline
36 & $\frac{1}{34}(1,9)$ & $(4,5,2)$ & $(\frac{5}{17},\frac{3}{17},\frac{10}{17})$ & $-\frac{15}{17}$ \\
\hline
37 & $\frac{1}{34}(1,11)$ & $(4,2^{10})$ & $(\frac{5+i}{17})_{i=1}^{11}$ & $\frac{165}{17}$ \\
\hline
38 & $\frac{1}{34}(1,13)$ & $(3,3,3,2)$ & $(\frac{7}{17},\frac{4}{17},\frac{5}{17},\frac{11}{17})$ & $\frac{33}{17}$ \\
\hline
39 & $\frac{1}{34}(1,15)$ & $(3,2,2,3,2,2)$ & $(\frac{8}{17},\frac{7}{17},\frac{6}{17},\frac{5}{17},\frac{9}{17},\frac{13}{17})$ & $\frac{81}{17}$ \\
\hline
40 & $\frac{1}{34}(1,27)$ & $(2,2,2,3,2^5)$ & $(\frac{14}{17},\frac{11}{17},\frac{8}{17},(\frac{3+2i}{17})_{i=1}^6)$ & $\frac{141}{17}$ \\
\hline
41 & $\frac{1}{37}(1,5)$ & $(8,2,3)$ & $(\frac{6}{37},\frac{11}{37},\frac{16}{37})$ & $-\frac{96}{37}$ \\
\hline
42 & $\frac{1}{37}(1,6)$ & $(7,2^5)$ & $(\frac{2+5i}{37})_{i=1}^6$ & $\frac{72}{37}$ \\
\hline
43 & $\frac{1}{37}(1,7)$ & $(6,2,2,3)$ & $(\frac{5+3i}{37})_{i=1}^4$ & $\frac{12}{37}$ \\
\hline
44 & $\frac{1}{37}(1,8)$ & $(5,3,3)$ & $(\frac{9}{37},\frac{8}{37},\frac{15}{37})$ & $-\frac{80}{37}$ \\
\hline
45 & $\frac{1}{37}(1,9)$ & $(5,2^8)$ & $(\frac{7+3i}{37})_{i=1}^9$ & $\frac{252}{37}$ \\
\hline
46 & $\frac{1}{37}(1,10)$ & $(4,4,2,2)$ & $(\frac{11}{37},\frac{7}{37},\frac{17}{37},\frac{27}{37})$ & $\frac{36}{37}$ \\
\hline
47 & $\frac{1}{37}(1,11)$ & $(4,2,3,2,2)$ & $(\frac{12}{37},\frac{11}{37},\frac{10}{37},\frac{19}{37},\frac{28}{37})$ & $\frac{108}{37}$ \\
\hline
48 & $\frac{1}{37}(1,12)$ & $(4,2^{11})$ & $(\frac{11+2i}{37})_{i=1}^{12}$ & $\frac{396}{37}$ \\
\hline
49 & $\frac{1}{37}(1,17)$ & $(3,2,2,2,2,3,2)$ & $((\frac{19-i}{37})_{i=1}^6,\frac{25}{37})$ & $\frac{216}{37}$ \\
\hline
50 & $\frac{1}{37}(1,21)$ & $(2,5,2^4)$ & $(\frac{22}{37},(\frac{1+6i}{37})_{i=1}^5)$ & $\frac{132}{37}$ \\
\hline
51 & $\frac{1}{37}(1,22)$ & $(2,4,2^6)$ & $(\frac{23}{37},(\frac{5+4i}{37})_{i=1}^7)$ & $\frac{240}{37}$ \\
\hline
52 & $\frac{1}{37}(1,23)$ & $(2,3,3,2,2,2)$ & $(\frac{24}{37},\frac{11}{37},(\frac{2+7i}{37})_{i=1}^4)$ & $\frac{168}{37}$ \\
\hline
53 & $\frac{1}{38}(1,5)$ & $(8,3,2)$ & $(\frac{3}{19},\frac{5}{19},\frac{12}{19})$ & $-\frac{53}{19}$ \\
\hline
54 & $\frac{1}{38}(1,7)$ & $(6,2,4)$ & $(\frac{4}{19},\frac{5}{19},\frac{6}{19})$ & $-\frac{29}{19}$ \\
\hline
55 & $\frac{1}{38}(1,9)$ & $(5,2,2,2,3)$ & $(\frac{4+i}{19})_{i=1}^5$ & $\frac{43}{19}$ \\
\hline
56 & $\frac{1}{38}(1,15)$ & $(3,3,2^6)$ & $(\frac{8}{19},(\frac{3+2i}{19})_{i=1}^7)$ & $\frac{127}{19}$ \\
\hline
57 & $\frac{1}{38}(1,21)$ & $(2,6,2,2,2)$ & $(\frac{11}{19},(\frac{-1+4i}{19})_{i=1}^4)$ & $\frac{31}{19}$ \\
\hline
58 & $\frac{1}{38}(1,27)$ & $(2,2,4,2^4)$ & $(\frac{14}{19},\frac{9}{19},(\frac{1+3i}{19})_{i=1}^5)$ & $\frac{103}{19}$ \\
\hline
59 & $\frac{1}{39}(1,5)$ & $(8,5)$ & $(\frac{2}{13},\frac{3}{13})$ & $-\frac{70}{13}$ \\
\hline
60 & $\frac{1}{39}(1,7)$ & $(6,3,2,2)$ & $(\frac{8}{39},\frac{3}{13},\frac{19}{39},\frac{29}{39})$ & $\frac{2}{39}$ \\
\hline
61 & $\frac{1}{39}(1,11)$ & $(4,3,2^4)$ & $(\frac{4}{13},(\frac{1+2i}{13})_{i=1}^5)$ & $\frac{50}{13}$ \\
\hline
62 & $\frac{1}{39}(1,14)$ & $(3,5,3)$ & $(\frac{5}{13},\frac{2}{13},\frac{5}{13})$ & $-\frac{10}{13}$ \\
\hline
63 & $\frac{1}{39}(1,16)$ & $(3,2,5,2)$ & $(\frac{17}{39},\frac{4}{13},\frac{7}{39},\frac{23}{39})$ & $\frac{38}{39}$ \\
\hline
64 & $\frac{1}{39}(1,17)$ & $(3,2,2,4,2)$ & $(\frac{6}{13},\frac{5}{13},\frac{4}{13},\frac{3}{13},\frac{8}{13})$ & $\frac{38}{13}$ \\
\hline
65 & $\frac{1}{39}(1,25)$ & $(2,3,2,2,3,2)$ & $(\frac{2}{3},\frac{1}{3},\frac{1}{3},\frac{1}{3},\frac{1}{3},\frac{2}{3})$ & $\frac{14}{3}$ \\
\hline
66 & $\frac{1}{39}(1,29)$ & $(2,2,3,2^8)$ & $(\frac{10}{13},\frac{7}{13},(\frac{3+i}{13})_{i=1}^9)$ & $\frac{134}{13}$ \\
\hline
67 & $\frac{1}{39}(1,31)$ & $(2,2,2,3,2^6)$ & $(\frac{32}{39},\frac{25}{39},\frac{6}{13},(\frac{7+4i}{39})_{i=1}^7)$ & $\frac{362}{39}$ \\
\hline
68 & $\frac{1}{41}(1,5)$ & $(9,2^4)$ & $(\frac{-1+7i}{41})_{i=1}^5$ & $-\frac{40}{41}$ \\
\hline
69 & $\frac{1}{41}(1,6)$ & $(7,6)$ & $(\frac{7}{41},\frac{8}{41})$ & $-\frac{220}{41}$ \\
\hline
70 & $\frac{1}{41}(1,8)$ & $(6,2^7)$ & $(\frac{5+4i}{41})_{i=1}^8$ & $\frac{200}{41}$ \\
\hline
71 & $\frac{1}{41}(1,9)$ & $(5,3,2,2,2)$ & $(\frac{10}{41},(\frac{1+8i}{41})_{i=1}^4)$ & $\frac{80}{41}$ \\
\hline
72 & $\frac{1}{41}(1,10)$ & $(5,2^9)$ & $(\frac{8+3i}{41})_{i=1}^{10}$ & $\frac{320}{41}$ \\
\hline
73 & $\frac{1}{41}(1,11)$ & $(4,4,3)$ & $(\frac{12}{41},\frac{7}{41},\frac{16}{41})$ & $-\frac{28}{41}$ \\
\hline
74 & $\frac{1}{41}(1,12)$ & $(4,2,4,2)$ & $(\frac{13}{41},\frac{11}{41},\frac{9}{41},\frac{25}{41})$ & $\frac{44}{41}$ \\
\hline
75 & $\frac{1}{41}(1,13)$ & $(4,2^5,3)$ & $(\frac{13+i}{41})_{i=1}^7$ & $\frac{212}{41}$ \\
\hline
76 & $\frac{1}{41}(1,16)$ & $(3,3,2,2,3)$ & $(\frac{17}{41},\frac{10}{41},\frac{13}{41},\frac{16}{41},\frac{19}{41})$ & $\frac{128}{41}$ \\
\hline
77 & $\frac{1}{41}(1,17)$ & $(3,2,4,2,2)$ & $(\frac{18}{41},\frac{13}{41},\frac{8}{41},\frac{19}{41},\frac{30}{41})$ & $\frac{116}{41}$ \\
\hline
78 & $\frac{1}{41}(1,23)$ & $(2,5,3,2)$ & $(\frac{24}{41},\frac{7}{41},\frac{11}{41},\frac{26}{41})$ & $\frac{32}{41}$ \\
\hline
79 & $\frac{1}{41}(1,26)$ & $(2,3,2,3,2,2)$ & $(\frac{27}{41},\frac{13}{41},\frac{12}{41},\frac{11}{41},\frac{21}{41},\frac{31}{41})$ & $\frac{188}{41}$ \\
\hline
80 & $\frac{1}{41}(1,34)$ & $(2^4,3,2^5)$ & $((\frac{41-6i}{41})_{i=1}^5,(\frac{11+5i}{41})_{i=1}^5)$ & $\frac{380}{41}$ \\
\hline
\end{longtable}
\end{center}
\end{lem}
\begin{proof}
By Lemma \ref{lem: summary four singular points condition, II}(4)(5), the order of $X\ni x_4$ belongs to  $\{25,26,29,31,34,37,38,39,41\}$. The rest is enumeration.
\end{proof}

\end{document}